\documentclass[a4paper, 11pt]{article}
\usepackage{amsmath}
\usepackage{amssymb,amsmath,amsthm,color}
\usepackage{lineno,hyperref}
\usepackage[table,x11names,dvipsnames,table]{xcolor}
\usepackage{authblk}
\usepackage{subcaption,booktabs}
\usepackage{graphicx}
\usepackage{booktabs} 
\usepackage{multirow} 
\usepackage[nolist,nohyperlinks]{acronym}
\usepackage[superscript]{cite}
\usepackage{tabularx}
\usepackage{float}
\usepackage[group-digits=none,table-text-alignment=left]{siunitx}
\usepackage{geometry}
\usepackage{amsmath}
\usepackage{algorithm}
\usepackage{algpseudocode}
\usepackage[utf8]{inputenc}
\usepackage[T1]{fontenc}
\usepackage{lmodern}
\usepackage{amsmath}
\usepackage{algorithm}
\usepackage{algpseudocode}
\usepackage{amsfonts}
\usepackage{colortbl} 
\usepackage{booktabs} 
\usepackage{enumitem}
\usepackage{adjustbox}
\newlist{dasheditemize}{itemize}{1}
\setlist[dasheditemize,1]{label=--, left=0em, itemsep=0.5em}
 \usepackage{hyperref}		
\hypersetup{
  colorlinks=true,
  citecolor=blue,  
  linkcolor=blue,
  urlcolor=blue
}
\usepackage{url}
\usepackage[square,numbers]{natbib}

\usepackage{titlesec}
\usepackage{pifont}

\theoremstyle{remark}
\newtheorem{remark}{Remark}
\newtheorem{example}{Example}

\titleformat{\section}
  {\normalfont\Large\bfseries} 
  {\thesection.} 
  {0.5em} 
  {} 

\titleformat{\subsection}
  {\normalfont\large\itshape} 
  {\thesubsection.} 
  {0.5em} 
  {} 

\titleformat{\subsubsection}
  {\normalfont\large\itshape} 
  {\thesubsubsection.} 
  {0.5em} 
  {} 

\begin{document}

\makeatletter
\renewcommand*{\@fnsymbol}[1]{\ensuremath{\ifcase#1\or *\or \text{\ding{86}} \or
   \mathsection\or \mathparagraph\or \|\or **\or \dagger\dagger
   \or \ddagger\ddagger \else\@ctrerr\fi}}
\makeatother
\title{\LARGE \textbf{An Efficient Exponential Sum Approximation of Power-Law Kernels for Solving Fractional Differential Equations}\footnote{This work was supported by the German Federal Ministry of Research, Technology and Space under Grant No.~05M22WHA.}}

\author{\Large Renu Chaudhary, Kai Diethelm, Afshin Farhadi\thanks{Corresponding author\newline \textit{Email addresses:} renu.chaudhary@thws.de (R. Chaudhary), kai.diethelm@thws.de (K. Diethelm), \\ afshin.farhadi@thws.de (A. Farhadi), fredalbert.fuchs@thws.de (F.A. Fuchs). } , Fred A. Fuchs}
\affil{\small \textit{Faculty of Applied Natural Sciences and Humanities (FANG), Technical University of Applied Sciences
W\"{u}rzburg-Schweinfurt, Ignaz-Sch\"{o}n-Str. 11, 97421 Schweinfurt, Germany}}
\vspace{0.3cm}
\date{}
\maketitle
\vspace{-0.5cm}
\hrule
\vspace{0.3cm}
\textbf{Abstract.} In this work, we present a comprehensive framework for approximating the weakly singular power-law kernel $t^{\alpha-1}$ of fractional integral and differential operators, where $\alpha \in (0,1)$ and $t \in [\delta,T]$ with $0 < \delta < T < \infty$, using a finite sum of exponentials. This approximation method begins by substituting an exponential function into the Laplace transform of the power function, followed by the application of the trapezoidal rule to approximate the resulting integral. To ensure computational feasibility, the integral limits are truncated, leading to a finite exponential sum representation of the kernel.  In contrast to earlier approaches, we pre-specify the admitted computational cost (measured in terms of the number of exponentials) and minimize the approximation error. Furthermore, to reduce computational cost while maintaining accuracy, we present a two-stage algorithm based on Prony’s method that compresses the exponential sum. The compressed kernel is then embedded into the Riemann-Liouville fractional integral and applied to solve fractional differential equations. To this end, we discuss two solution strategies, namely (a) method based on piecewise constant interpolation and (b) a transformation of the original fractional differential equation into a system of first-order ordinary differential equations (ODEs). This reformulation makes the problem solvable by standard ODE solvers with low computational cost while retaining the accuracy benefits of the exponential-sum approximation. Finally, we apply the proposed strategies to solve some well-known fractional differential equations and demonstrate the advantages, accuracy and the experimental order of convergence of the methods through numerical results.

\vspace{0.3cm}
\hrule

\section{Introduction}

The Riemann-Liouville fractional integral is a fundamental concept in fractional calculus, characterized (up to a multiplicative constant) by having the form of a Laplace-type convolution with kernel function $w(t)= t^{\alpha-1}$ where $\alpha > 0$ and $t \in (0,T]$. Therefore, the accurate and efficient approximation of power-law kernels plays a crucial role in the computation of fractional integrals as required, e.g., in typical algorithms for solving fractional differential equations. In the vast majority of applications \cite{BDST,HFCA7,HFCA8,Diethelm2010,HFCA6,HFCA4,HFCA5}, one needs to deal with the case that $0 < \alpha < 1$ . In this case, the  singularity at $t=0$ and the slow decay of the kernel make direct numerical computations challenging. The sum-of-exponentials method addresses this challenge by approximating the kernel $t^{\alpha-1}$ with a sum of exponential functions which are easier to handle numerically and lead to more efficient algorithms. 

A variety of sum-of-exponentials methods has been introduced for this purpose. For example, Baffet \cite{Baffet2019} developed a composite Gauss–Jacobi quadrature on the Laplace-domain integral representation of $t^{\alpha-1}$. By subdividing the integral into finite intervals and truncating the infinite tails, he obtained a rapidly convergent approximation of the form
\begin{equation*}
w(t) \approx \sum_{l=1}^{P}b_l\mathrm e^{-a_lt},
\end{equation*}
valid for $t\in[\delta,T]$, with real coefficients $a_l,b_l$ chosen by the quadrature. This scheme delivers high accuracy with a modest number of terms $P$, uniformly in $\alpha\in(0,1)$, and provides explicit error bounds. He also incorporated the compressed kernel into a time-stepping algorithm to solve fractional differential equations. Baffet and Hesthaven \cite{baffet2017kernel} split the fractional integral of a function $f$ into local and history parts and proposed a kernel compression scheme based on a multipole approximation of the Laplace transform of the kernel function. This method transforms the history part into a linear combination of auxiliary variables, each defined by a standard ordinary differential equation, so that the fractional problem can be converted to a system of ODEs. They  also derived estimates for the number of auxiliary variables required to attain a prescribed error tolerance. 
In a different approach, Li \cite{rebecca2010} exploited the representation of the kernel as a Laplace transform,
\begin{equation*}
	t^{\alpha-1} =\mathcal{L}\{\xi^{-\alpha}/ \Gamma(1-\alpha)\} (t),
\end{equation*}
and constructed a fixed $Q$-point quadrature rule for the integral $\int_0^{\infty}\mathrm e^{-\xi t}\xi^{-\alpha}\,\mathrm d\xi$ that is accurate to a tolerance $\epsilon$ for all $t$. The resulting fast time-stepping algorithm uses this quadrature with time-independent nodes and weights, avoiding repeated re-computation. This algorithm does not impose any constraints on the final simulation time, making it suitable for modeling systems over very long durations. The storage requirement of Li's algorithm is $O(Q)$ and its complexity is $O(NQ)$, where $N$ is the number of time steps.

Beylkin and Monzón \cite{beylkin2005} introduced an algorithmic framework based on Prony’s method and Hankel matrices for approximating the power function $t^{-\beta}$, $\beta> 0$, on a finite interval by a short sum of exponentials. Their method produces, for any prescribed accuracy, a representation with far fewer terms than the standard Fourier representations would require. In a follow-up study \cite{beylkin2010approximation}, they achieved approximations for such functions with uniform relative error on the whole real line. In their method, they applied Poisson summation to discretize the integral representations of these functions. They also introduced a straightforward algorithm based on Prony's method to remove those terms with small exponents from the initial approximations.

Building on the Beylkin-Monzón approach, McLean \cite{mclean2018exponential} revisited the exponential-sum approximation of $t^{-\beta}$, $\beta > 0$.  In order to achieve a desired accuracy, he compared the performance of the exponential-sum approximation of the Laplace integral representation using the substitutions $p =\mathrm e^{x}$ and $p = \exp(x-\mathrm e^{-x})$ both before and after applying Prony's method. He demonstrated that the second substitution results in significantly fewer terms prior to the application of Prony's method, while both approximations yield approximately the same number of exponential terms once Prony's method is applied.
 
In this paper, we develop an efficient exponential-sum approximation for the power-law kernel $t^{\alpha-1}$  in the given interval with given fractional order $\alpha$, and $L$ exponential terms. Specifically, we pre-specify the computational cost (i.e.\ the number of terms $L$ in the sum of exponentials) and then minimize the approximation error. This approach differs from existing relative error tolerance-driven approaches that compute quadrature ranges and hence the number of exponential terms $L$ to meet a prescribed error and accept the resulting computational cost. Starting from the Laplace-transform integral for $t^{\alpha-1}$, we apply a composite trapezoidal rule to approximate it on $[\delta,T]$, truncating the infinite tails to render the computation finite. We then introduce Prony's method to reduce the number of exponential terms, thereby decreasing memory usage and computational time. The result is a compact exponential-sum representation of the power-law kernel, matching or improving upon the initial maximum error before applying Prony's method. Further, we embed this compressed kernel into Riemann-Liouville fractional integral for solving fractional differential equations, and present two solution strategies based on constant interpolation and first-order ODE formulations. The reduced ODE system is solved with the backward Euler method or the trapezoidal rule. Numerical experiments demonstrate that the proposed methods are practical and advantageous with significantly reduced computational cost for solving a wide variety of problems in fractional calculus.  We also show the accuracy and effectiveness of the different solution schemes depends upon fractional order $\alpha$ and number of exponentials $L$ and discuss experimental order of convergence (EOC) for the tested schemes.

\section{Kernel Approximation by Exponential Sums}
\label{sec:approx-exp-sum}

\subsection{The General Idea}
Our aim is to construct an efficient numerical approximation by a sum of exponentials to the kernel function $ t^{\alpha-1}$ of Riemann-Liouville fractional integrals of order $\alpha\in (0,1)$ for functions $f\in C[0,T]$ (where $T$ is a real positive number) with a starting point at $0$ represented as
\begin{equation}\label{RLIntegral}
	I^{\alpha} f(t) = \frac{1}{\Gamma(\alpha)}\int_0^t(t-\tau)^{\alpha-1}f(\tau) \, \mathrm d\tau.
\end{equation}

The convolution kernel function $t^{\alpha-1}$ of this operator can be written as the Laplace transform of a power function in the form
\begin{equation}\label{LapTrans}
	t^{\alpha-1} = \frac{1}{\Gamma(1-\alpha)} \int_{0}^{\infty}\mathrm e^{-pt}p^{-\alpha}\, \mathrm dp.
\end{equation}
Choosing $ p = \mathrm e^{\omega}$, eq. (\ref{LapTrans}) reads
\begin{equation}\label{KerExp1}
	t^{\alpha-1} = \frac{1}{\Gamma(1-\alpha)} \int_{-\infty}^{\infty}\mathrm e^{(1-\alpha)\omega}\mathrm e^{-t\mathrm e^{\omega}}\, \mathrm d\omega.
\end{equation}

Since the integral $(\ref{KerExp1})$ is over the unbounded domain $(-\infty,\infty)$, we truncate the domain to a finite interval $[l_{\min},l_{\max}]$, where $l_{\min} < 0$ and $l_{\max} > 0$ are real numbers. Then for some given $L \in \mathbb{N}$, by generating the nodes $\omega_{l} = l_{\min} + (l-1)h$ for $l=1,\ldots,L$, we divide this interval into $L-1$ subintervals of width $ h =(l_{\max}-l_{\min})/(L-1)$, and finally we apply the trapezoidal rule to obtain
\begin{equation}\label{KerExp2}
	t^{\alpha-1} \approx \frac{1}{\Gamma(1-\alpha)} \int_{l_{\min}}^{l_{\max}}g(\omega)\, \mathrm d\omega 
	\approx \frac{h}{\Gamma(1-\alpha)} \bigg(\frac{g(\omega_{1})+g(\omega_{L})}{2}+\sum_{l=2}^{L-1}g(\omega_{l})\bigg)
\end{equation}
where 
\begin{equation}
	\label{eq:def-g}
	g(\omega) =\mathrm e^{(1-\alpha)\omega}\mathrm e^{-t\mathrm e^{\omega}}. 
\end{equation}
A finite exponential sum approximation for the kernel function can then be obtained as
\begin{equation}\label{TrapExpFin}
 	t^{\alpha-1} \approx f(t,\theta) := \frac{1}{\Gamma(1-\alpha)}\sum_{l=1}^{L} w_{l}\mathrm e^{b_{l}t} \qquad \text{for} \quad t \in [\delta,T],
\end{equation}
where $0 < \delta < T$, $ \theta = (w_1, \dots, w_L, b_1, \dots, b_L) \in \mathbb{R}^{2L}$ is the vector of parameters $w_l$ and $b_{l}$ called the weights and exponents, respectively, and 
\begin{equation}\label{weightandnode}
\begin{aligned}
	w_{1} &=\frac{h}{2}\mathrm e^{(1-\alpha)l_{\min}}, & \quad & w_{L} =\frac{h}{2}\mathrm e^{(1-\alpha)l_{\max}},\\
	w_{l} &=h\mathrm e^{(1-\alpha)\omega_{l}} & \qquad &\text{for} \quad l = 2,\ldots,L-1, \\
	b_{l} &=-\mathrm e^{\omega_{l}} &\qquad &\text{for} \quad l = 1,\ldots,L.
\end{aligned}
\end{equation}
With this approach, our goal is to find a set of values $ \theta = (w_1, \dots, w_L, b_1, \dots, b_L) \in \mathbb{R}^{2L} $ such that the absolute error
\begin{equation}\label{Abserror}
	e(t,\theta) =  t^{\alpha-1} - \frac{1}{\Gamma(1-\alpha)}\sum_{l=1}^{L} w_{l} \mathrm e^{b_{l}t} 
\end{equation}
is small in modulus for $t  \in [\delta, T]$, where data $\alpha, \delta, T$ and $L$ are given. As such, our course of action differs from the approach of McLean \cite{mclean2018exponential} who computed suitable values for $l_{\min}$ and $l_{\max}$ and hence, indirectly, $L$ based on a given relative error tolerance for the finite sum approximation $(\ref{TrapExpFin})$. Thus, while McLean intends to achieve a certain given accuracy with whatever computational cost is necessary, we pre-specify the computational cost that we are prepared to invest (which essentially depends on $L$) and want to minimize the error.

\begin{remark}
	The case of interest for us here is when $\alpha \in (0,1)$. Therefore, the function $t^{\alpha-1}$ is unbounded as $t \to 0$, a feature that the approximation $f(t, \theta)$ cannot reproduce for any $\theta$.
	This is the reason why, as indicated above, we aim for an approximation of $t^{\alpha-1}$ only on the interval $[\delta, T]$ with some $\delta > 0$, not on the entire interval $(0,T]$.
	We will explain in Section \ref{Sec3} how we deal with the interval $(0, \delta]$ that is not covered by this approach.
\end{remark}

\begin{remark}
	\label{rem:rescaling}
	Note also that we can apply a simple rescaling technique to transfer our results from a general interval $[\delta, T]$ to the interval $[\delta / T, 1]$ which is normalized in the sense that its right end point is fixed at $1$ (or vice versa): 
	If $\theta = (w_1, \dots, w_L, b_1, \dots, b_L) $ is a choice of the weights and exponents for the interval $[\delta, T]$ and
	\begin{equation}\label{rescale}
		\tilde w_l = T^{1-\alpha} w_l, \quad 
		\tilde b_l = T b_l,
		\qquad (l=1,\dots,L)
	\end{equation}
	then $\tilde \theta = (\tilde w_1, \dots, \tilde w_L, \tilde b_1, \dots, \tilde b_L)$ can be used on $[\delta/T, 1]$, 
	and the associated errors are connected to each other via
	\begin{equation}
		e(t, \theta) = T^{\alpha-1} e(\tilde t, \tilde \theta).
	\end{equation}
	where $t \in [\delta, T]$ and $\tilde t = t / T \in [\delta / T, 1]$.
\end{remark}

The first technical step for us is to determine \( l_{\min} \) and \( l_{\max} \) for truncating the integral  $(\ref{KerExp1})$ such that the integral's contribution outside the interval $[l_{\min},l_{\max}]$ is negligible. We need the integrals 
\begin{equation}\label{3}
I_{\min} =\int_{-\infty}^{l_{\min}} g(\omega) \, \mathrm d\omega \qquad  \text{and} \qquad I_{\max} =\int_{l_{\max}}^{\infty} g(\omega) \, \mathrm d\omega
\end{equation}
to be less than a small threshold $\epsilon$.

\subsection{Condition for  $I_{\min} \le \epsilon $}
To find a suitable choice for $l_{\min}$, we consider the integral $I_{\min}$ in $(\ref{3})$. By substituting $u = \mathrm e^{\omega}$, and  $v = t u$, we get 
\[
I_{\min} = \frac{1}{t^{1-\alpha}}\int_{0}^{t\mathrm e^{l_{\min}}} v^{-\alpha} \mathrm e^{-v} \, \mathrm dv = \frac{1}{t^{1-\alpha}} \gamma(1-\alpha, t \mathrm e^{l_{\min}})
\]
where $\gamma$ is the lower incomplete Gamma function defined as 
\begin{equation}\label{6}
	\gamma(\beta,q) = \int_{0}^{q}  p^{\beta-1}  \mathrm e^{-p}  \, \mathrm dp \quad \text{for} \quad \beta > 0 \quad \text{and} \quad q >0,
\end{equation}
see \cite[\S 8.2]{DLMF}. For small $s = t \mathrm e^{l_{\min}}$, the function $\gamma$ has the asymptotic behavior $\gamma(1-\alpha, s) \sim s^{1-\alpha}/(1-\alpha)$ by \cite[eqs.~(8.2.6) and (8.7.1)]{DLMF}. Thus, $I_{\min} \approx \frac{\exp((1-\alpha) l_{\min})}{1-\alpha}$. Since we need $I_{\min} \le \epsilon$, we obtain
\begin{equation}\label{7}
	l_{\min} \le \frac{\ln (\epsilon (1-\alpha))}{1-\alpha}.
\end{equation}
Additionally, we must make sure that $s = t \mathrm e^{l_{\min}}$ is sufficiently small. Hence, we must suppose that $t \mathrm e^{l_{\min}} < \epsilon$ for all $t \in [\delta, T]$, i.e.\ $T \mathrm e^{l_{\min}} < \epsilon$. We thus also require
\[
	l_{\min} \le \ln \left(\frac{\epsilon}{T}\right).
\]
Combining these two inequalities, we see that we can choose a lower boundary for truncating the integral $(\ref{KerExp1})$ in the form
\begin{equation}\label{8}
	l_{\min} = \min \left\{ \ln \left(\frac{\epsilon}{T}\right), \frac{\ln (\epsilon (1-\alpha))}{1-\alpha}\right\}.
\end{equation}

\begin{figure}[h!t]
\centering
\begin{subfigure}[b]{0.49\textwidth}
    \includegraphics[width=\textwidth]{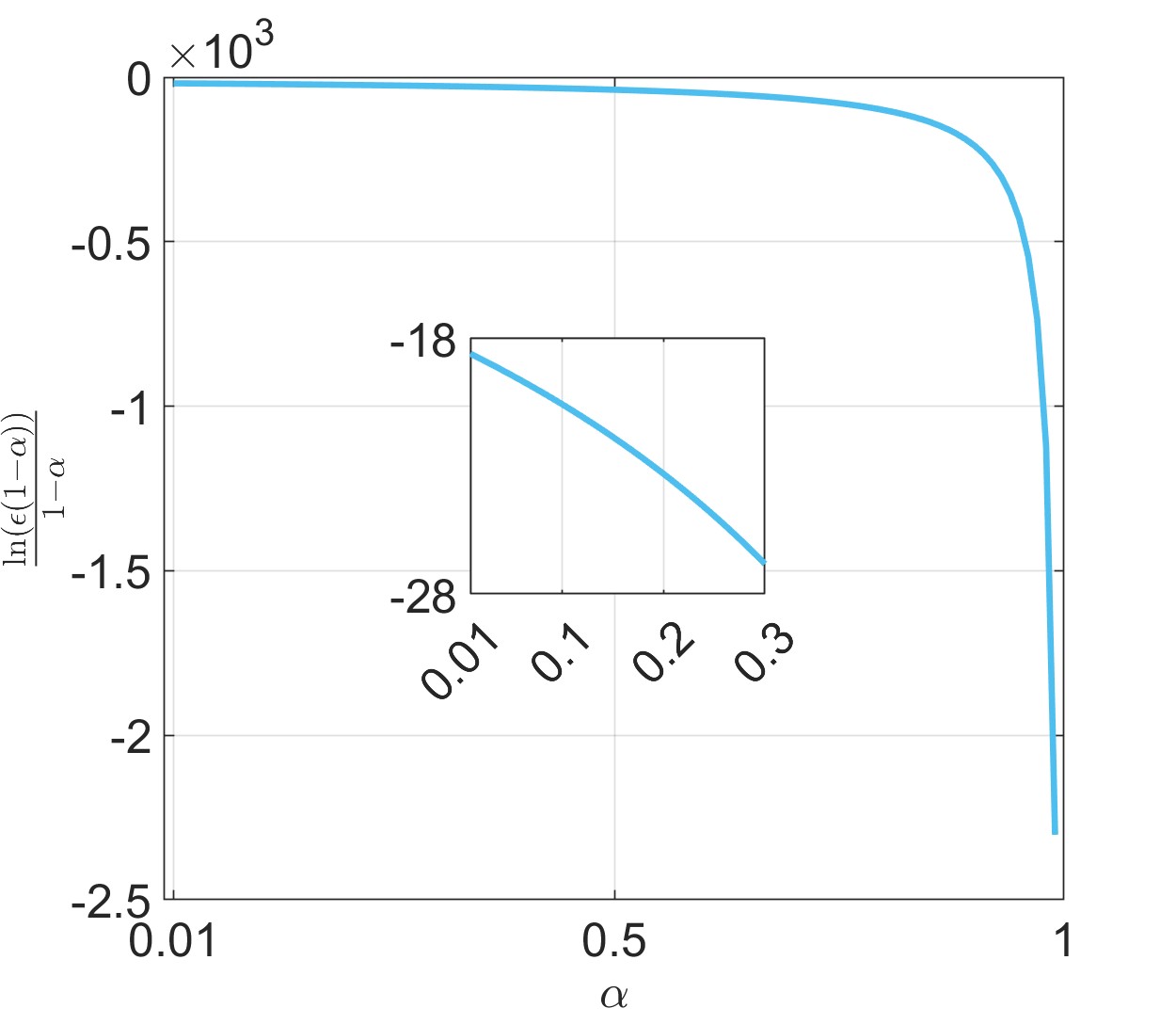}
    \caption{}
    \label{subfig:1-a}
\end{subfigure}
\hfill
\begin{subfigure}[b]{0.49\textwidth}
    \includegraphics[width=\textwidth]{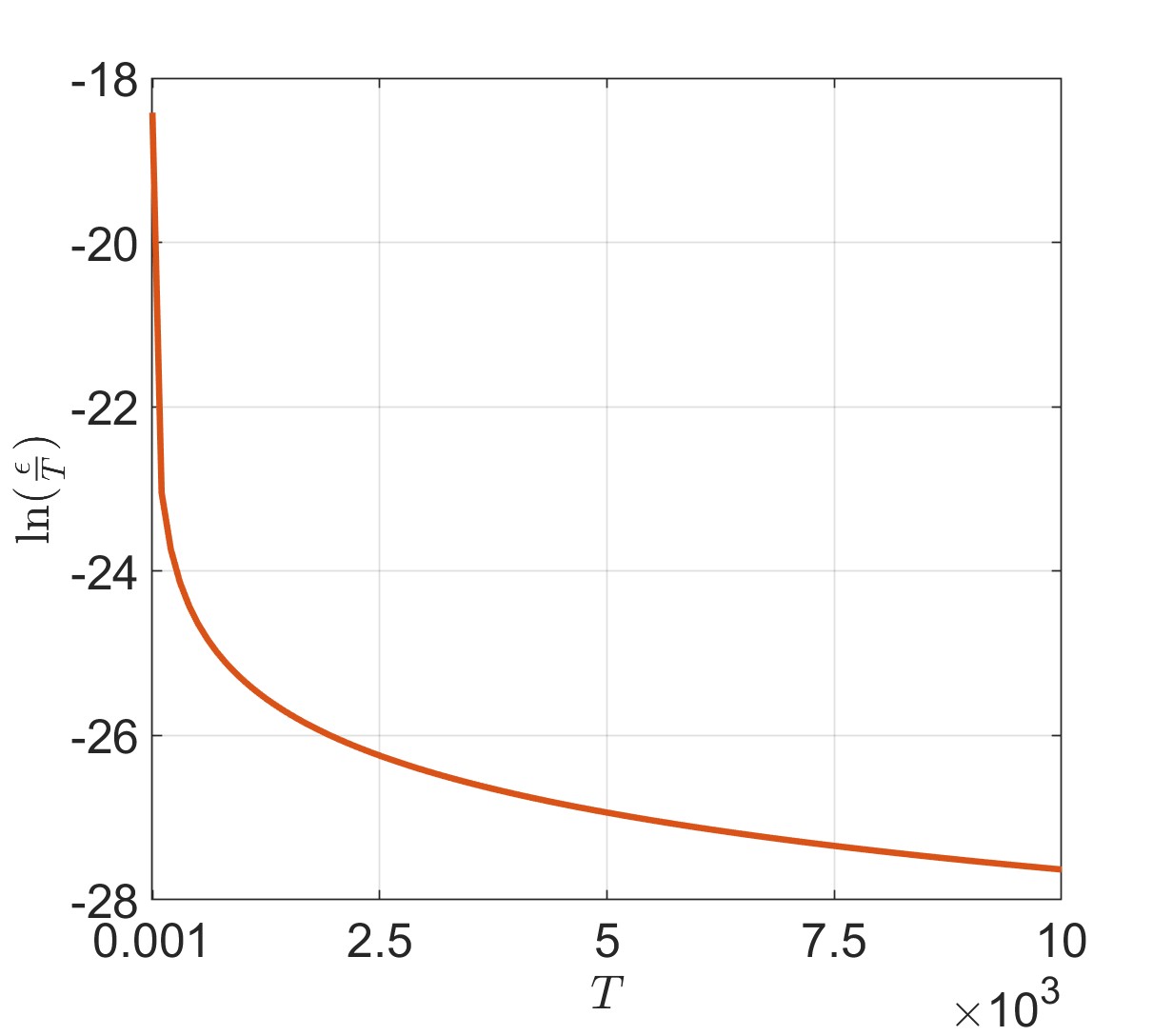}
    \caption{}
    \label{subfig:1-b}
\end{subfigure}
\caption{(a) Behaviour of $\frac{\ln (\epsilon (1-\alpha))}{1-\alpha}$ with respect to the fractional order $\alpha$. (b) Behavior of $\ln \left(\frac{\epsilon}{T}\right)$ with respect to the final time $T$.}
\label{fig:1}
\end{figure}

We show in Figure \ref{fig:1} how different values of $\alpha$ and $T$ can affect  $l_{\min}$. Figure \ref{subfig:1-a} shows the behaviour of $f_{1}(\alpha):=\frac{\ln (\epsilon (1-\alpha))}{1-\alpha}$ for different values of $\alpha$ in $[0.01,0.99]$ and $\epsilon = 10^{-8}$. Figure \ref{subfig:1-b} shows the behaviour of $f_{2}(T):=\ln \left(\epsilon/T\right)$ for different values of $T$ in $[1,10^{4}]$ and the same $\epsilon$. We here see that different values of $T$ from the interval $[1,10^{4}]$ can have an effect on our choice for $l_{\min}$ if $ \alpha$ in $[0.01,0.3]$ because both values of $f_{1}$ and $f_{2}$ change in the interval $[-18,-28]$. But for values of $\alpha$ greater than $0.3$, $T$ does not affect $l_{\min}$. Thus, the value of $\alpha$ is highly significant for approximating $l_{\min}$. It should be noted that values of $\alpha$ closer to $1$ result in more negative $l_{\min}$ and for fixed values of $\delta$ and $\epsilon$ a larger truncated interval $[l_{\min}, l_{\max}]$ is then obtained. 

\subsection{ Condition for $I_{\max} \le \epsilon $}
Next, we consider the integral $I_{\max}$ in $(\ref{3})$. Similarly as above, by substituting $u = \mathrm e^{\omega}$ and $v = t u$, we get
\[
	I_{\max} = \frac{1}{t^{1-\alpha}} \int_{t \mathrm e^{l_{\max}}}^{\infty} v^{-\alpha} \mathrm e^{-v} \, \mathrm dv = \frac{1}{t^{1-\alpha}} \Gamma(1-\alpha, t \mathrm e^{l_{\max}}),
\]
where $\Gamma$ is the upper incomplete Gamma function \cite[\S 8.2]{DLMF} defined as 
\begin{equation}\label{4}
	\Gamma(\beta,q) = \int_{q}^{\infty}  p^{\beta-1}  \mathrm e^{-p} \, \mathrm dp \quad \text{for} \quad \beta > 0 \text{ and } q >0.
\end{equation}
For large $s = t \mathrm e^{l_{\max}}$, the function $\Gamma$ asymptotically behaves as $\Gamma(1-\alpha, s) \sim s^{-\alpha} \mathrm e^{-s}$, see \cite[\S~8.11(i)]{DLMF}. Thus
\begin{equation*}
	I_{\max} 
	\approx \frac{1}{t^{1-\alpha}} (t \mathrm e^{l_{\max}})^{-\alpha} \exp(-t \mathrm e^{l_{\max}}) 
	= \frac{1}{t} \mathrm e^{-\alpha l_{\max}} \exp(-t \mathrm e^{l_{\max}}).
\end{equation*}
Now, $I_{\max} \le \epsilon$ results in $ -\ln t - \alpha l_{\max} - t \mathrm e^{l_{\max}} \le \ln \epsilon$. Since the dominant term is $t \mathrm e^{l_{\max}}$, an estimate for $l_{\max}$ dependent on $t$ is then obtained as 
\begin{equation}\label{lmin_approx1}
	l_{\max} (t) \approx \ln \left(-\frac{\ln \epsilon}{t}\right) \quad  \text{for} \quad \delta \leq t \leq T.
\end{equation}
Since for $t \geq \delta$, $l_{\max}(t) \leq l_{\max}(\delta)$, we can choose an upper bound for truncating the integral at $t = \delta$, viz.\ 
\begin{equation}\label{lmin_approx2}
	l_{\max} = \ln \left(-\frac{\ln \epsilon}{\delta}\right).
\end{equation}

\begin{remark}\label{lmax}
In the case $\delta > \ln (1 / \epsilon)$, we obtain $l_{\max} < 0$. While this case is theoretically possible, in practical numerical settings one typically has $\delta \ll \ln (1 / \epsilon)$, so such a situation is unlikely to occur. Moreover, even if it does arise, the algorithm remains well defined and may in fact yield an especially strong reduction in the number of terms after applying Prony’s method. (See also Remark \ref{rem:sum2empty}.)
\end{remark}

\begin{figure}[h!t]
\centering
\begin{subfigure}[b]{0.49\textwidth}
    \includegraphics[width=\textwidth]{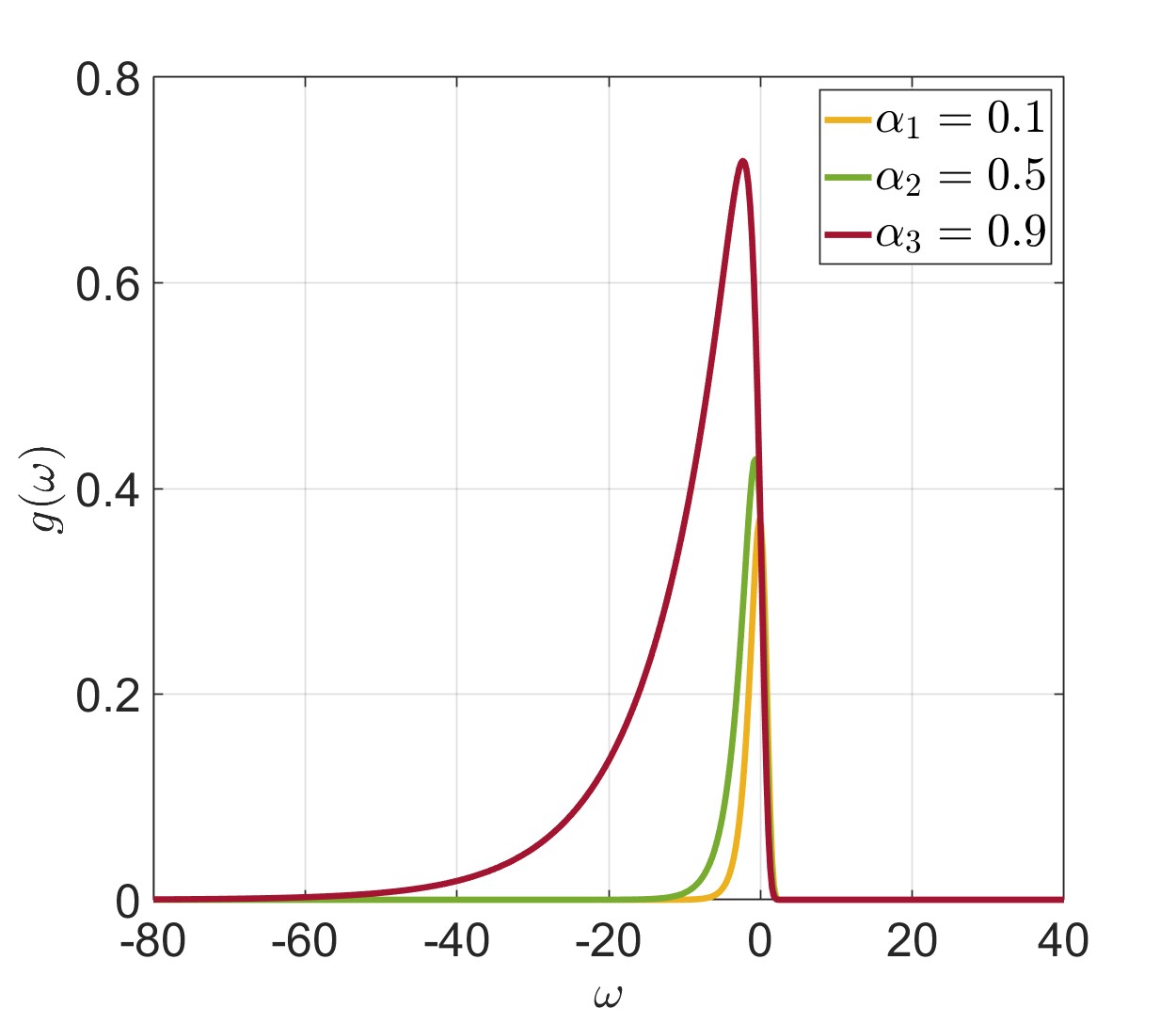}
    \caption{}
    \label{subfig:2-a}
\end{subfigure}
\hfill
\begin{subfigure}[b]{0.49\textwidth}
    \includegraphics[width=\textwidth]{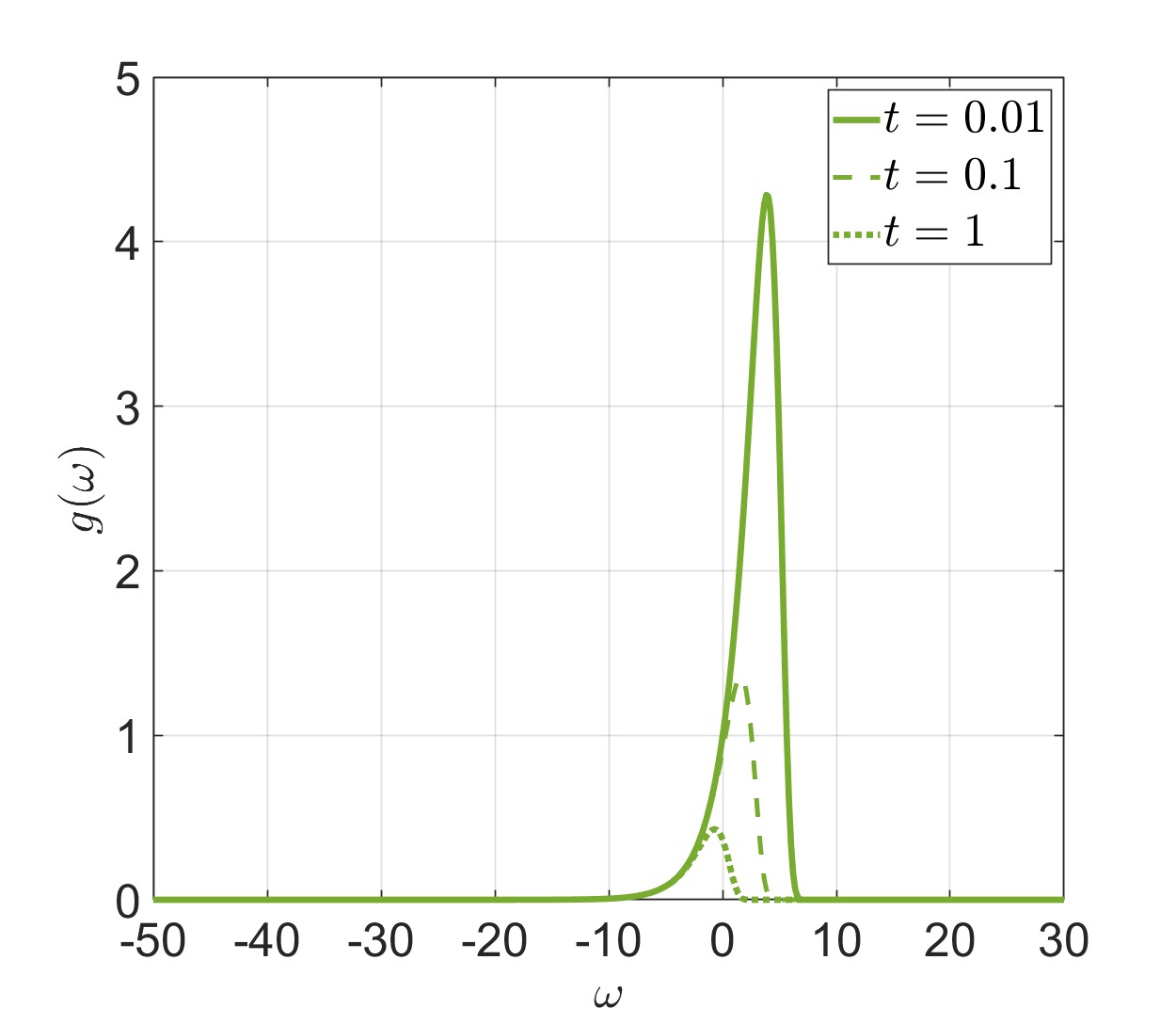}
    \caption{}
    \label{subfig:2-b}
\end{subfigure}
\caption{Behaviour of the integrand from (\ref{KerExp1}) at (a) $t=1$ for different values of $\alpha$, (b) at different values of $t$ for $\alpha = 0.5$.}
\label{fig:2}
\end{figure}

Figure \ref{subfig:2-a} shows the behaviour of $g(\omega)$ as defined in eq.~\eqref{eq:def-g} for different values of $\alpha$ in $(0,1)$ at $t=1$. In all cases, the function decays very rapidly as one moves from the local maximum of $g$ to the right. Thus, we may always choose a rather small positive value for the truncation point $l_{\max}$, independent of $\alpha$. This is because the $\alpha$-dependent term in the asymptotic limit of the upper incomplete Gamma function is subdominant.  Figure \ref{subfig:2-b} shows the behaviour of $g(\omega)$ for different values of $t$ in $[\delta,T]$ and for a fixed value of $\alpha$. The larger values of $l_{\max}$ obtained by decreasing $t$ such that its maximum value is at $t = \delta$. We here choose this value, that is, $ \ln \left(-\frac{\ln \epsilon}{\delta}\right)$  independent of $t$ as an approximation of $l_{\max}$. So, the value of $\delta$ plays a significant role in determining $l_{\max}$. It should be noted that by decreasing $\delta$, $l_{\max}$ become more positive and for fixed values of $\alpha$ and $\epsilon$ a larger truncated interval $[l_{\min}, l_{\max}]$ is then obtained.

\begin{remark}
As another approach, one can make an estimate of $l_{\min}$ and $l_{\max}$ by analyzing the integrand $g(\omega) = \mathrm e^{(1-\alpha)\omega} \mathrm e^{-t \mathrm e^{\omega}} $ to ensure that its contribution outside \([l_{\min}, l_{\max}]\) is negligible, i.e., we require $\mathrm e^{(1-\alpha)\omega-t \mathrm e^{\omega}}  < \epsilon$. For large positive $\omega$, $ \mathrm e^{\omega} $ grows rapidly, and the term $-t \mathrm e^{\omega}$ dominates. Thus $\omega > \ln \left(-\frac{\ln \epsilon}{\delta}\right)$ and this yields an approximation of $l_{\max}$ as it is given by $(\ref{lmin_approx2})$. For large negative $\omega$, $\mathrm e^{\omega} \to 0 $ leads to $ g(\omega) \approx \mathrm e^{(1-\alpha)\omega}$, and then we get $\omega < \frac{\ln \epsilon}{1-\alpha}$. Thus, we obtain an approximation as $l_{\min} \approx \frac{\ln \epsilon}{1-\alpha}$. This estimate for $l_{\min}$ does not depend on $T$ but shows that it is effectively dependent on $\alpha$ similar to the first approach. It should be noted that the estimate $\frac{\ln \epsilon}{1-\alpha}$ for different values of $\alpha$ and a fixed value of $\epsilon$ exhibits similar behavior to that shown in Figure \ref{subfig:1-a}, although with very small differences in their values. However, in the numerical examples presented here, we determine the values of $l_{\min}$ and $l_{\max}$ using the first approach.
\end{remark}

\section{Prony's Method}\label{Sec2}

To improve the level of accuracy in the kernel approximation, a straightforward idea is to increase the value of $L$, i.e.\ to use more terms in the sum of exponentials. However, using larger values of $L$ leads to more evaluations and hence an increase in memory consumption and computation time. Thus, it is essential to equip our numerical algorithm for approximating the kernel with efficient methods that can reduce the number of exponentials in the sum while preserving accuracy.

\begin{figure}[h]
    \centering
    \begin{subfigure}{0.33\textwidth}
        \centering
        \includegraphics[width=\textwidth]{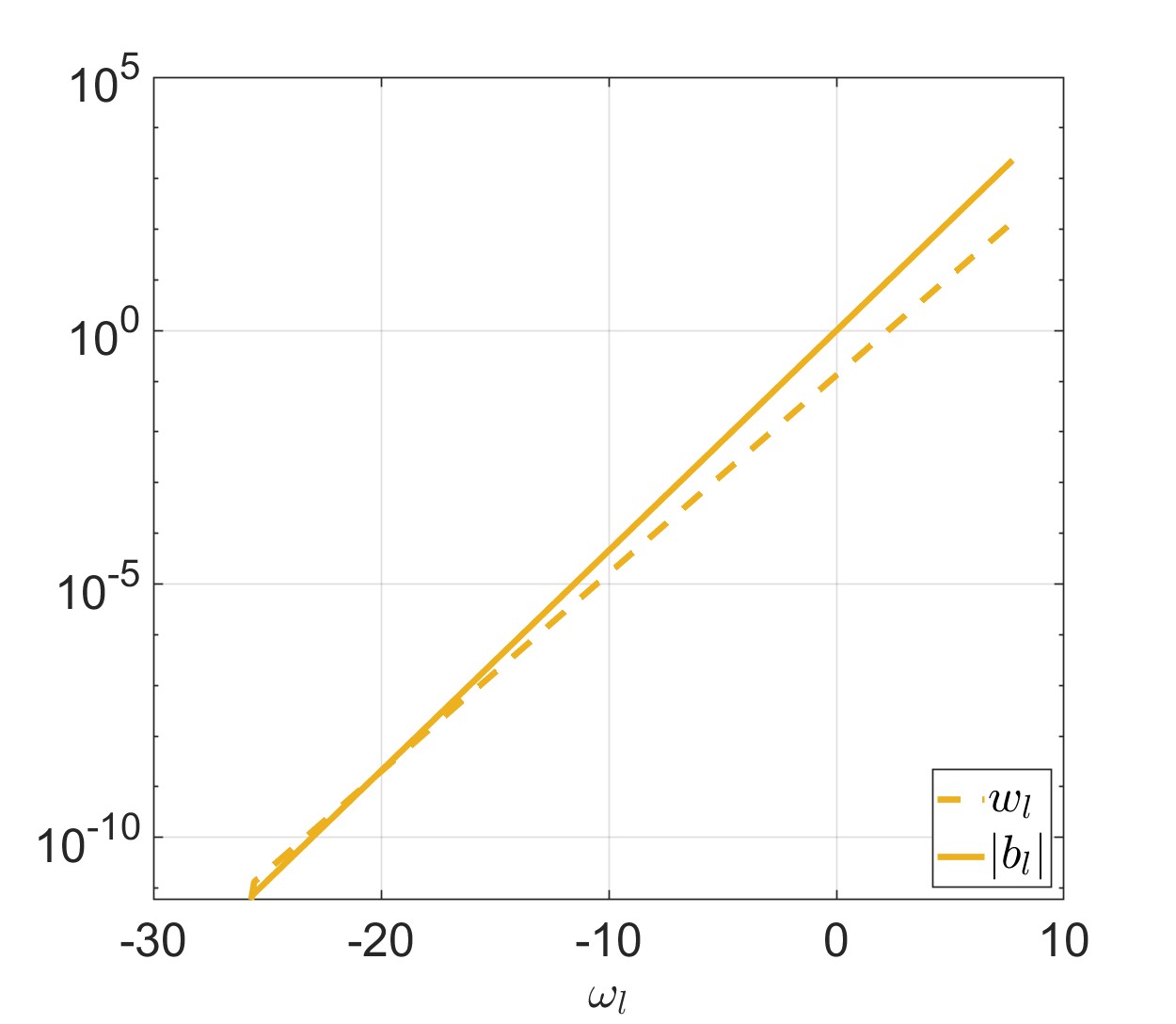}
        \caption{$\alpha = 0.1$}
        \label{fig:image1}
    \end{subfigure}\hfill
    \begin{subfigure}{0.33\textwidth}
        \centering
        \includegraphics[width=\textwidth]{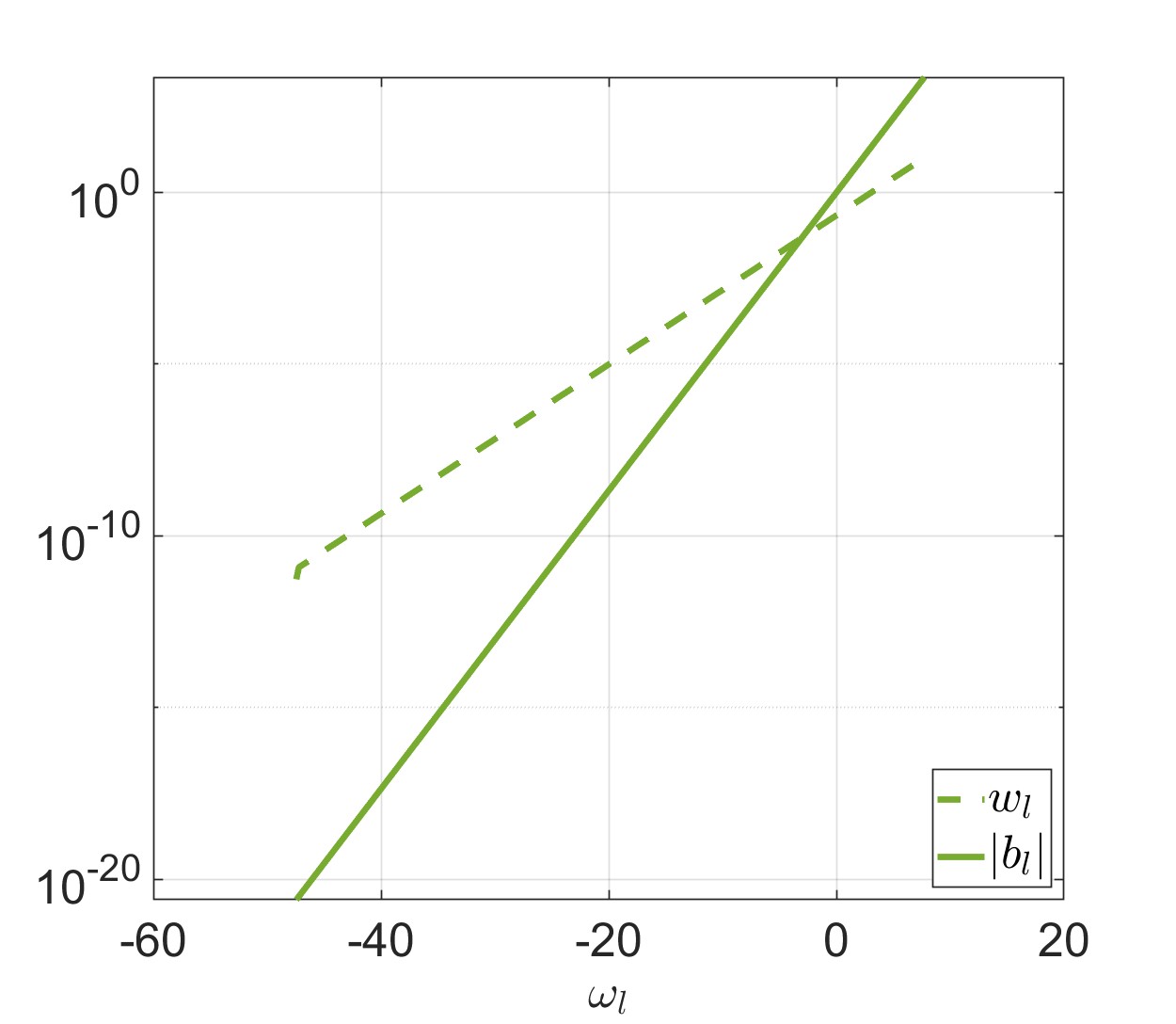}
        \caption{$\alpha = 0.5$}
        \label{fig:image2}
    \end{subfigure}\hfill
    \begin{subfigure}{0.33\textwidth}
        \centering
        \includegraphics[width=\textwidth]{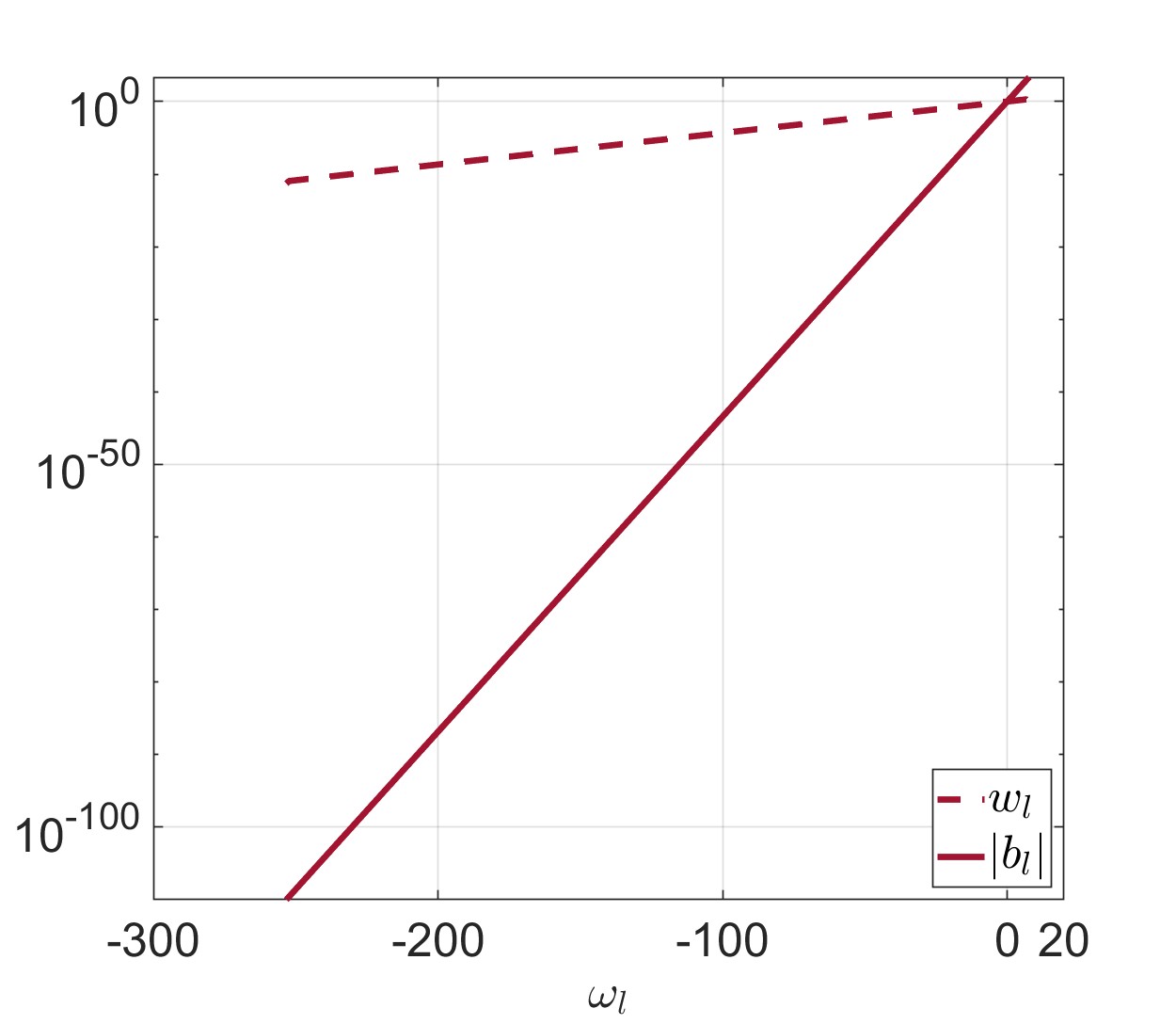}
        \caption{$\alpha = 0.9$}
        \label{fig:image3}
    \end{subfigure}
    \caption{Comparison of weights and exponents obtained by the kernel approximation $t^{\alpha-1}$ over the interval $[10^{-2},1]$ for varying $\alpha$ values $(0.1, 0.5, 0.9)$ with respect to nodes on the truncation domain $[l_{\min},l_{\max}]$. Here,  the value of $\epsilon$ is set to $10^{-10}$. Note that the vertical axis is shown on a logarithmic scale.}
    \label{fig:3}
\end{figure}

Originally introduced in \cite{Prony} and discussed in general in \cite{Sauer} and with an emphasis on our specific application in \cite{beylkin2005, beylkin2010approximation, mclean2018exponential}, the Prony's method allows us to reduce the number of terms in the sum on the right-hand side of eq.~\eqref{TrapExpFin}, i.e.\ our approximation of the kernel of the integral operator, without losing accuracy. Specifically, we split up the sum as
\[
	\sum_{l=1}^L w_l \mathrm e^{b_l t} = \sum_{\omega_l \le 0} w_l \mathrm e^{b_l t} + \sum_{\omega_l > 0} w_l \mathrm e^{b_l t},
\]
and replace the first sum on the right-hand side by an expression with the same structure but a smaller number of summands while not changing the second sum. This results in an approximation of the kernel function using fewer exponential terms. Figure \ref{fig:3} illustrates the behavior of weights and exponents $(w_{l}, b_{l})$ for $l = 1, \ldots, 256$ versus the nodes $\omega_l$ on the truncation domain $[l_{\min},l_{\max}]$ obtained by approximating $t^{\alpha-1}$ over the interval $[10^{-2},1]$ for various values of $\alpha$ with $\epsilon = 10^{-10}$. With $\omega_{l}$ in the interval $[l_{\min}, 0]$ the values of $b_l$ for $\alpha = 0.1$ range from $-6.8873 \times 10^{-12}$ to $-0.8806$; for $\alpha = 0.5$,  from $-2.5 \times 10^{-21}$ to $-0.9528$; and for $\alpha = 0.9$, from $-9.9 \times 10^{-111}$ to $-0.6394$. As we see, the exponents $ b_{l} = -\mathrm {e}^{\omega_{l}} $, where $ \omega_{l} $ is within the interval $[l_{\min}, 0] $, mostly take values that are very near to zero in the interval $(-1,0)$. This allows us to justify the application of Prony's method, resulting in an effective reduction. Furthermore, the weights begin with very small positive values and increase, and their maximum values decrease as $\alpha$ increases. 
For example,
\[
	\max_{l=1, 2, \ldots, 256} w_l = \begin{cases} 
							123.7369 & \text{for } \alpha = 0.1, \\
							\phantom{12}9.3188 & \text{for } \alpha = 0.5,\text{ and} \\
							\phantom{12}2.0041 & \text{for } \alpha = 0.9.
						\end{cases}
\]

For a precise description of Prony's method, we follow McLean \cite{mclean2018exponential}. Recalling first that $\omega_l = l_{\min} + (l-1) (l_{\max} - l_{\min})/(L-1)$ for $l = 1, 2, \ldots, L$, we set
\begin{equation}
	\label{eq:def-M}
	M := \max \{ l \in 1, 2, \ldots, L : \omega_l \le 0 \}
\end{equation}
so that $M \le L$, $\omega_l \le 0$ for $l = 1, 2, \ldots, M$ and $\omega_{M+1} > 0$ if $M < L$. Furthermore, we choose some $L_{\text p} \in \{1, 2, \ldots, M \}$ and a positive integer $K$ such that $2K - 1 \le L_{\text p}$. The goal of Prony's method is to replace the first $L_{\text p}$ summands in the right-hand side of eq.~\eqref{TrapExpFin} by 
\begin{equation}
	\label{eq:pronysum}
	\sum_{k=1}^{K} \rho_{k}\mathrm e^{\eta_{k}t},
\end{equation}
i.e.\ by a sum with at most a bit more than half as many (and, if possible, even fewer) terms as the original one, where the parameters $\rho_k$ and $\eta_k$ are chosen such that we can maintain the desired accuracy, i.e.\ such that
\begin{equation}\label{pronyerror}
	\varepsilon_{\text p}(t) = \frac{1}{\Gamma(1-\alpha)} \left( \sum_{l=1}^{L_{\text p}} w_{l}\mathrm e^{b_{l}t} - \sum_{k=1}^{K} \rho_{k}\mathrm e^{\eta_{k}t} \right)
\end{equation}
is small enough in modulus for all $t \in [\delta,T]$, thus permitting to reduce the total number of terms from $L$ to 
\[
	L_{\text f} = K + L - L_{\text p}.
\] 
Based on Prony's method, in order to find the $2K$ parameters $\rho_{k}$ and $\eta_{k}$, we take the following steps:
\begin{enumerate}
\item Check that $2K-1 \leq L_{\text p}$ and compute 
\begin{equation*}
	g_{j} = \sum_{l=1}^{L_{\text p}}w_{l}b_{l}^{j}   \qquad (j = 0, 1,, \ldots, 2K-1).
\end{equation*}
\item \label{item:hankel} Construct the $K \times K$ Hankel matrix 
\begin{equation*}
	\mathbf{H}=
		\begin{pmatrix}
			g_{0}& g_{1} & \cdots   &   g_{K-1}           \\
			g_{1} & g_{2} & \cdots &   g_{K}       \\
			 \vdots  & \vdots & \ddots & \vdots      \\
			g_{K-1} & g_{K} & \cdots &   g_{2K-2}       
		\end{pmatrix}
\end{equation*}
and the vector $\mathbf{b} = - (g_K, g_{K+1}, \ldots, g_{2K-1})^{\mathrm T}$,
and solve the linear system $\mathbf{H}\mathbf{q} = \mathbf{b}$ with $\mathbf q = (q_0, q_1, \ldots, q_{K-1})^{\mathrm T}$.
\item \label{item:polynomial} Form the monic polynomial $q(z) = q_{0}+q_{1}z+\ldots+q_{K-1}z^{K-1}+z^{K}$, and find its roots $\eta_{k}$ ($k=1,\ldots, K$) to obtain the Prony's exponents for eq.~\eqref{eq:pronysum}.
\item \label{item:vandermonde} Construct the $(2K) \times K$ Vandermonde matrix 
\begin{equation*}
\mathbf{V}=
\begin{pmatrix}
1& 1 & \cdots   &   1           \\
\eta_{1} & \eta_{2} & \cdots &  \eta_{K}       \\
 \vdots  & \vdots & \ddots & \vdots      \\
\eta_{1}^{2K-1} & \eta_{2}^{2K-1} & \cdots &   \eta_{K}^{2K-1}       
\end{pmatrix}
\end{equation*}
and solve the linear system $\mathbf V \boldsymbol \rho = \mathbf{g}$ where $\mathbf{g} =(g_{0}, g_{1},  \ldots, g_{2K-1})^{\mathrm T}$ and
$\boldsymbol \rho = (\rho_1, \rho_2, \ldots, \rho_k)^{\mathrm T}$ to determine the Prony's weights $\rho_{k}$ in eq.~\eqref{eq:pronysum}.
\end{enumerate}

\begin{remark}
The roots of the polynomial $q(z)$ that arises in item \ref{item:polynomial} above might not all be real and positive. In addition, the  $K \times K$ Hankel matrix mentioned in item \ref{item:hankel} might be nearly singular in some cases. In practice, the detailed discussion in  \cite{beylkin2005, beylkin2010approximation, mclean2018exponential} shows that these issues only cause minor difficulties for our algorithm that can easily be circumvented.
\end{remark}

\begin{remark}
Although only the first $K$ rows of the Vandermonde system in item \ref{item:vandermonde} would suffice to uniquely determine the weights $\rho_{k}$, in practice this choice may lead to severe numerical ill-conditioning, especially when the nodes $\eta_{k}$ are very small. By formulating the problem as an overdetermined  $(2K) \times K$ system, we gain additional stability and robustness against numerical errors, while still recovering the same solution in well-conditioned cases.
\end{remark}

To determine suitable values for $K$ and $L_{\text p}$, we employ the method outlined in  \cite[Example~1]{mclean2018exponential} but present it as an algorithm. Assume that the original approximation $t^{\alpha-1} \approx \sum_{l=1}^{L} w_l \mathrm e^{b_l t} $ yields a maximum error of $\epsilon' $ over the interval $ [\delta, T] $. The steps of the method are as follows:

\begin{enumerate}
    \item Initialize with $K = 1$ and $L_{\text p} = M$ to maximize the initial reduction in the number of points.
    \item Compute the parameters of Prony's method:
    \begin{itemize}
        \item Calculate the exponents $\eta_k$ and weights $\rho_k$ for $k = 1, \dots, K$ using Prony's method as described above.
        \item Evaluate the Prony approximation error $\varepsilon_{\text p}(t)$ from $(\ref{pronyerror})$ for $t \in [\delta, T]$.
    \end{itemize}
    \item Check the error condition:
    \begin{itemize}
        \item If $\max_{t \in [\delta, T]} \varepsilon_{\text p}(t) \leq \epsilon'$, accept the current values of $K$ and $L_{\text p}$ as suitable choices for reducing the total number of terms. Terminate the algorithm.
        \item Otherwise, proceed to step $4$.
    \end{itemize}
    \item Update parameters:
    \begin{itemize}
        \item If $2K - 1 < L_{\text p}$, decrement $L_{\text p}$ by 1 and go back to step 2.
        \item If $2K - 1 \geq L_{\text p}$, increment $K$ by 1, reset $L_{\text p} = M$, and go back to step 2, provided $2K - 1 < L_{\text p}$. Repeat until the error condition is satisfied.
    \end{itemize}
\end{enumerate}
With these considerations, we can replace the initial sum-of-exponentials approximation $f(t, \theta_0)$, say, for $t^{\alpha-1}$ from eq.~\eqref{TrapExpFin} with a certain parameter vector $\theta_0$ by the new, computationally cheaper, approximation
\begin{equation}
	\label{TrapExpProny}
	t^{\alpha-1} \approx f_{\text p}(t, \theta_{\text p}) 
		= \frac 1 {\Gamma(1-\alpha)} \left( \sum_{k=1}^{K} \rho_{k}\mathrm e^{\eta_{k}t} + \sum_{l=L_{\text p}+1}^{L} w_{l}\mathrm e^{b_{l}t} \right)
\end{equation}
where $\theta_{\text p} = (\rho_1, \rho_2, \ldots, \rho_K, w_{L_{\text p} + 1}, w_ {L_{\text p} + 2}, \ldots, w_L, \eta_1, \eta_2, \ldots, \eta_K, b_{L_{\text p} + 1}, b_ {L_{\text p} + 2}, \ldots, b_L)$. 
For this final approximation after the Prony reduction step, we can immediately derive an upper bound for the absolute value of the error $e_{\text p}(t) = t^{\alpha-1} -  f_{\text p}(t, \theta_{\text p})$ according to
\begin{equation}
	|e_{\text p}(t)| 
	= | t^{\alpha-1} -  f_{\text p}(t, \theta_{\text p}) | 
	\le | t^{\alpha-1} -  f(t, \theta_0) | + | f(t, \theta_0) -  f_{\text p}(t, \theta_{\text p}) | 
	= |e(t,\theta_{0})| + |\varepsilon_{\text p}(t)|
\end{equation}
This gives us $\max_{t \in [\delta, T]} |e_{\text p}(t)| \leq \epsilon' + \epsilon' = 2 \epsilon'$. Therefore, after reduction, this algorithm allows us to preserve the maximal error close to the initial maximal error obtained before applying Prony's method. 

Table \ref{tab:1} presents a comparison of the maximum error values $ \max_{t \in [10^{-2}, 1]} |e(t, \theta_0)| $ before applying Prony's method with the corresponding values $ \max_{t \in [10^{-2}, 1]} |e_{\text p}(t)| $ after the Prony procedure for various values of $ L $  from $ 32 $ to $ 1024 $ and $ \alpha \in \{0.1, 0.5, 0.9\} $. The results in Table \ref{tab:1} show the following:
\begin{enumerate}
\item   In most cases, the maximum error after applying Prony's method, viz.\ $ \max_{t \in [10^{-2}, 1]} |e_{\text p}(t)| $, either does not change or actually improves upon the initial error $ \max_{t \in [10^{-2}, 1]} |e(t, \theta_0)| $. For example, 
\begin{itemize}
\item in the case  $ \alpha = 0.1$ and $L = 32 $, the error slightly decreases from $ 5.102467 \times 10^{-2} $ to $ 5.102386 \times 10^{-2} $, 
\item for $ \alpha = 0.1$ and $L = 128 $, the error significantly improves from $ 1.320726 \times 10^{-8} $ to $ 1.980379 \times 10^{-10} $, and 
\item for $ \alpha = 0.5$ and $L = 256 $, the error decreases from $ 3.518998 \times 10^{-10} $ to $ 5.593037 \times 10^{-11} $. 
\end{itemize}
This confirms that the algorithm maintains the accuracy within the theoretical bound $ 2\epsilon' $ or even enhances it.
\item The value of $ K $ is always small (ranging from 1 to 5). This indicates an efficient compression of the original exponential sum. For instance, for $ \alpha = 0.9, L = 128 $, only $ K = 2 $ terms are needed, reducing $ L_{\text f} $ to 6, demonstrating that a minimal number of terms can achieve comparable accuracy.
\item The algorithm initializes $L_{\text p} = M $, and the iterative process adjusts $ L_{\text p} $ and $ K $ to satisfy the error criterion $ \max_{t \in [10^{-2}, 1]} |\varepsilon_p(t) | \leq \epsilon' $. The table shows that in all cases, the final values of $L_{\text p}$ are equal to $ M $, indicating that the algorithm can use all available negative nodes initially, optimizing the reduction process.
\item  For all $\alpha$, the algorithm significantly reduces the number of terms ($ L_{\text f }\ll L $). As $\alpha$ increases, $L_{\text f}$ becomes particularly small, i.e.\ especially few terms are required for comparable accuracy. Thus, Prony’s method effectively reduces the number of terms in the exponential sum approximation.
\end{enumerate}
Therefore, Table \ref{tab:1} confirms that the proposed algorithm effectively reduces the number of terms in the exponential sum approximation while keeping the maximum error close to or better than the initial error. The small values of $ K $ and the significant reduction of the total number of terms, and hence the computational cost, from $L$ to $ L_{\text f} $ highlight the algorithm's efficiency, making it highly suitable for applications where computational efficiency and accuracy are critical. 

In all tables, the notation $5.102467(-2)$ with a mantissa and an exponent in parentheses represents the number $5.102467 \times 10^{-2}$.

\begin{table}[H]
\centering
\renewcommand{\arraystretch}{1.2} 
\caption{Max.\ error values obtained before and after the application of Prony's method for various values of $L$ and $\alpha \in \{0.1, 0.5, 0.9\}$ over the interval $[10^{-2}, 1]$. In all cases, the numbers of non-positive nodes $M$, the values of $L_{\text p}$ and $K$ in Prony's method and the final number of exponential terms after reduction $L_{\text f}$ are provided. Here, the value of $\epsilon$ is set to $10^{-10}$.}
\label{tab:1}
\begin{tabular}{cccccccc}
\toprule
$\boldsymbol{\alpha}$ & $\boldsymbol{L}$ & $\boldsymbol{M}$ & $\boldsymbol{L_{\text p}}$ & $\boldsymbol{K}$ & $\boldsymbol{L_{\text f}}$ & 
	$\boldsymbol{\max_{t \in [10^{-2},1]} |e(t,\theta_{0})|}$ & $\boldsymbol{\max_{t \in [10^{-2}, 1]} |e_{\text p}(t)|}$ \\
\midrule
\multirow{4}{*}{$0.1$} & $32$  & 24  & 24  & 1 & 9  & 5.102467$(-2)\phantom{0}$ & 5.102386$(-2)\phantom{0}$ \\
& $64$  & 49   & 49  & 3 & 18  & 6.483821$(-6)\phantom{0}$ & 6.510213$(-6)\phantom{0}$ \\
& $128$  & 98    & 98  & 4 & 34  & 1.320726$(-8)\phantom{0}$ & 1.980379$(-10)$ \\
& $256$  & 196  & 196 & 4 & 64  & 7.207191$(-9)\phantom{0}$ & 6.868319$(-10)$ \\

\midrule
\multirow{4}{*}{$0.5$} & $32$  & 27 & 27  & 1 & 6  & 8.401490$(-2)\phantom{0}$ & 8.401582$(-2)\phantom{0}$ \\
& $64$  & 55   & 55  & 2 & 11  & 3.577202$(-4)\phantom{0}$ &3.577193$(-4)\phantom{0}$ \\
& $128$  & 110   & 110 & 4 & 22  & 3.988015$(-9)\phantom{0}$ & 3.802676$(-9)\phantom{0}$ \\
& $256$  &220   & 220 & 5 & 41 & 3.518998$(-10)$ & 5.593037$(-11)$ \\

\midrule
\multirow{4}{*}{$0.9$} & $128$  & 124   & 124  & 2 & 6  & 4.027975$(-3)\phantom{0}$ & 4.027975$(-3)\phantom{0}$\\
& $256$  & 248    & 248  & 2 & 10   & 2.591330$(-5)\phantom{0}$ & 2.591330$(-5)\phantom{0}$ \\
& $512$  & 496    & 496 & 4 & 20  & 1.240738$(-9)\phantom{0}$ & 1.240076$(-9)\phantom{0}$ \\
& $1024$  & 993  & 993 &5 & 36 & 1.342770$(-11)$ & 1.039657$(-11)$ \\
\bottomrule
\end{tabular}
\end{table}

\begin{remark}
	\label{rem:sum2empty}
	If the parameters introduced in Section \ref{sec:approx-exp-sum} satisfy $\delta \ge \ln(1/\epsilon)$ then, by eq.~\eqref{lmin_approx2}, we have $\omega_L = l_{\max} \le 0$ and hence eq.~\eqref{eq:def-M} implies $M = L$. It is not uncommon for the algorithm described above (see the examples shown in Table \ref{tab:1}) to yield $L_{\text p} = M = L$, so that in this case the second sum on the right-hand side of eq.~\eqref{TrapExpProny} is empty, thus we obtain an especially large reduction in the number of terms.
\end{remark}

\begin{figure}[H]
    \centering
    \begin{subfigure}{0.33\textwidth}
        \centering
        \includegraphics[width=\textwidth]{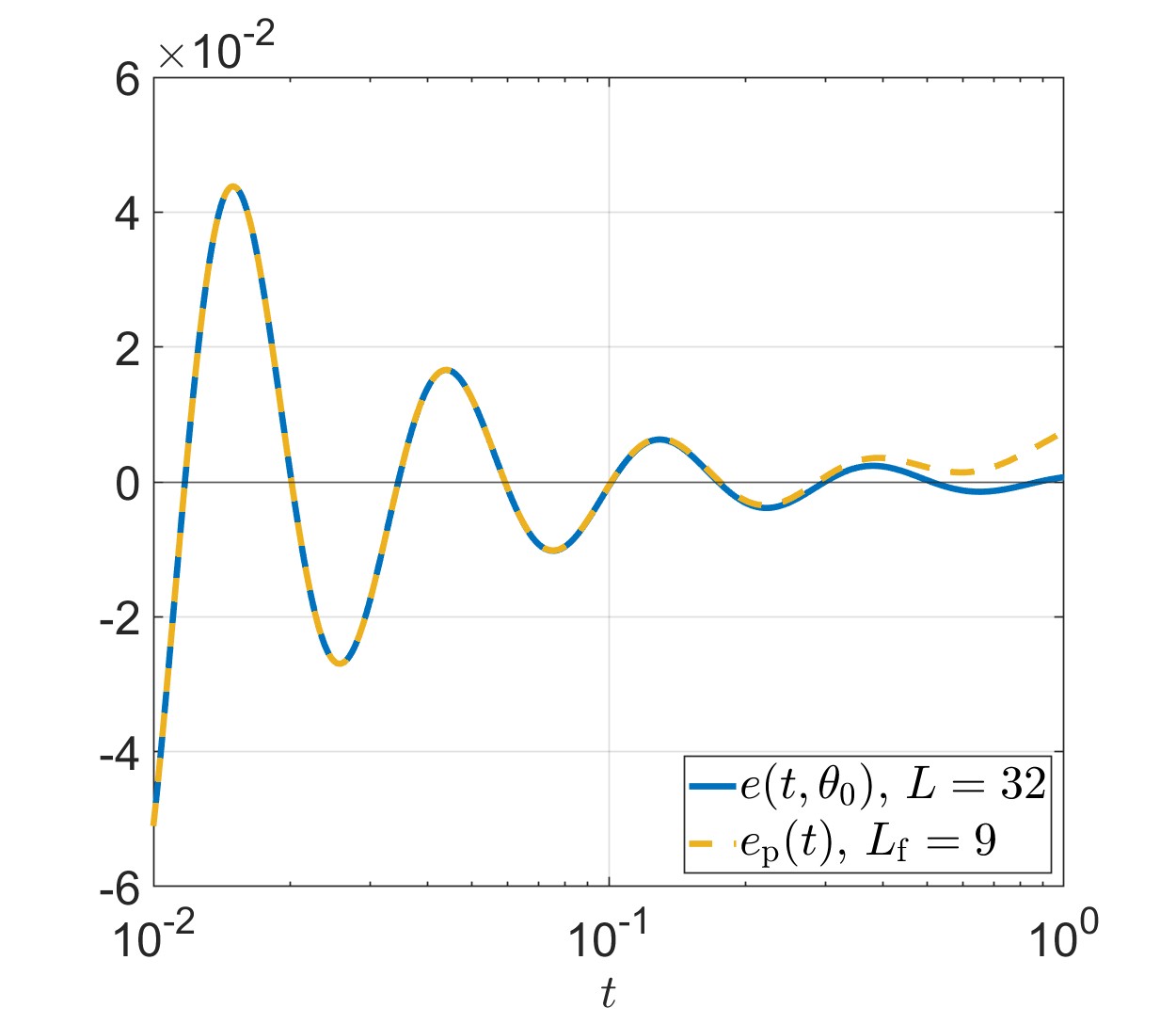}
        \caption{$\alpha = 0.1$}
        \label{fig:4_a}
    \end{subfigure}\hfill
    \begin{subfigure}{0.33\textwidth}
        \centering
        \includegraphics[width=\textwidth]{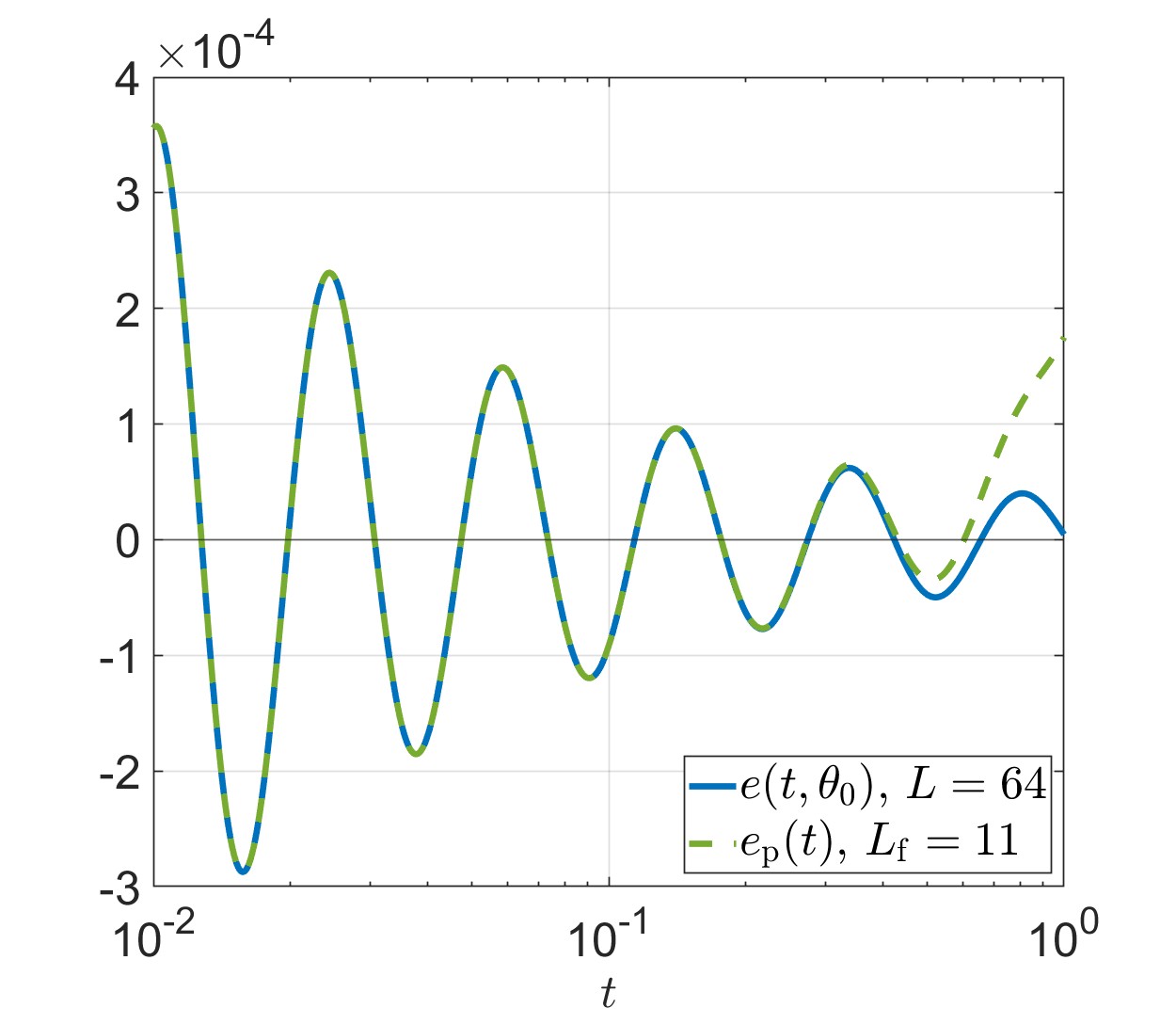}
        \caption{$\alpha = 0.5$}
        \label{fig:4_b}
    \end{subfigure}\hfill
    \begin{subfigure}{0.33\textwidth}
        \centering
        \includegraphics[width=\textwidth]{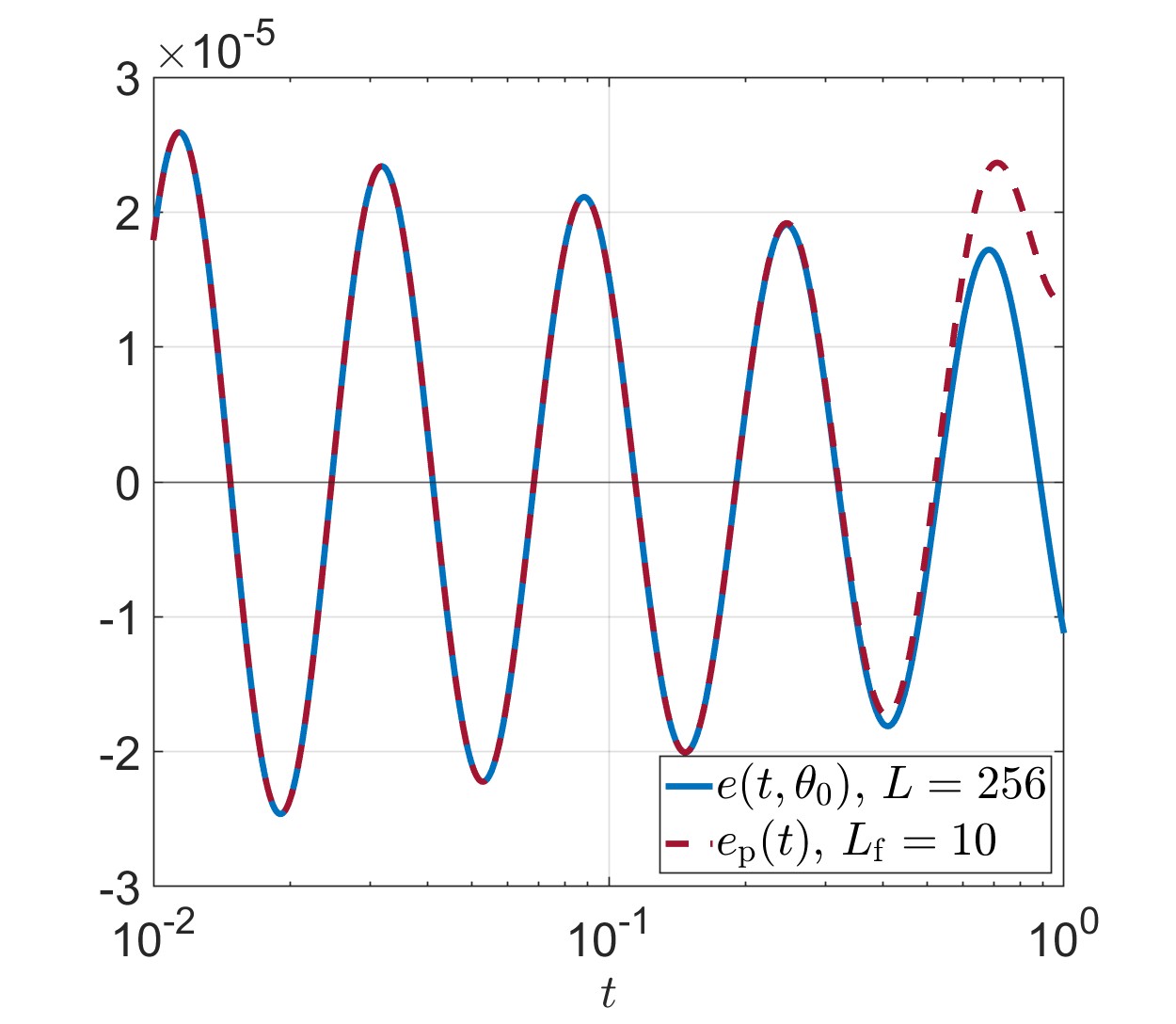}
        \caption{$\alpha = 0.9$}
        \label{fig:4_c}
    \end{subfigure}\hfill
        \begin{subfigure}{0.33\textwidth}
        \centering
        \includegraphics[width=\textwidth]{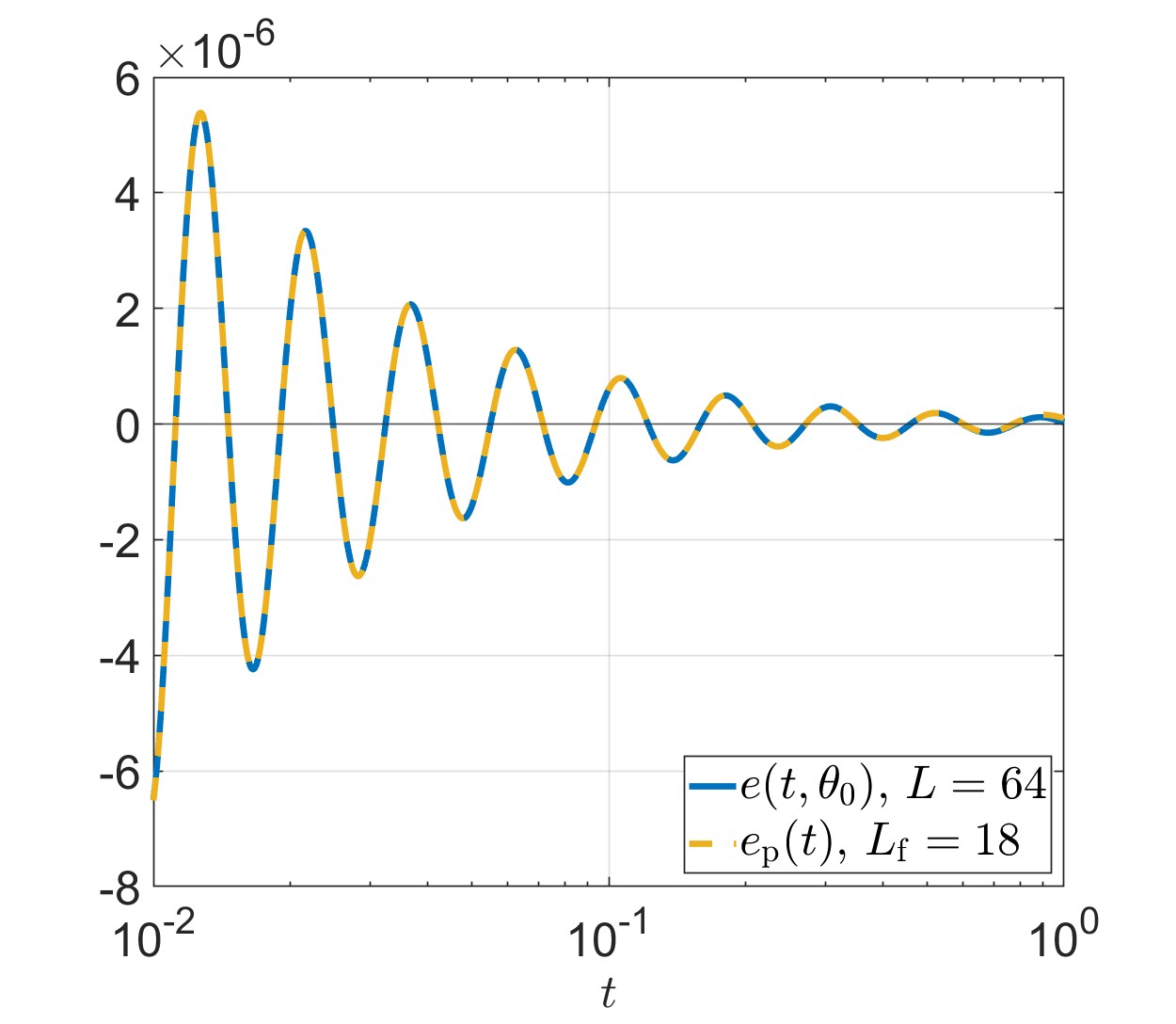}
        \caption{$\alpha = 0.1$}
        \label{fig:4_d}
    \end{subfigure}\hfill
    \begin{subfigure}{0.33\textwidth}
        \centering
        \includegraphics[width=\textwidth]{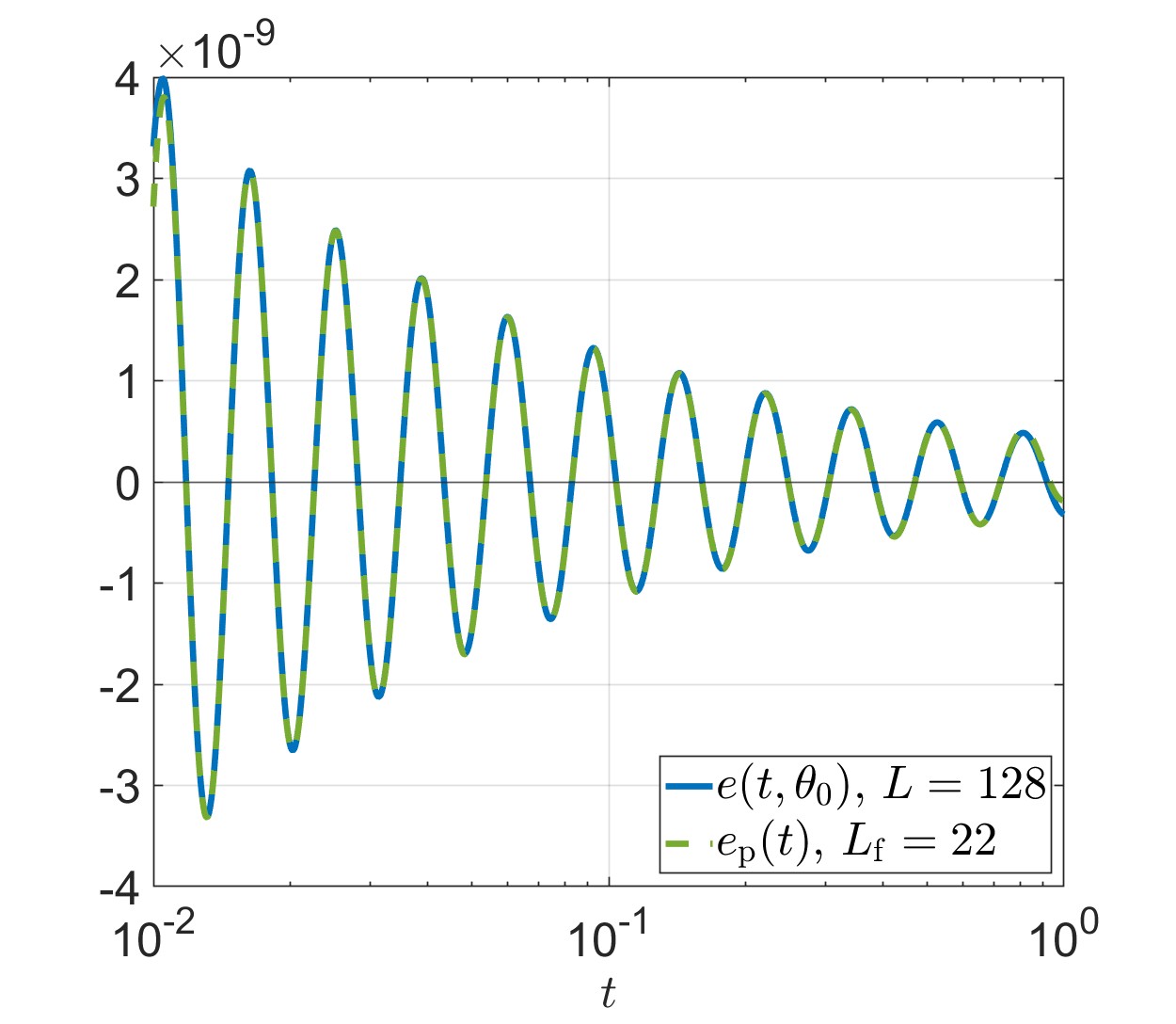}
        \caption{$\alpha = 0.5$}
        \label{fig:4_e}
    \end{subfigure}\hfill
    \begin{subfigure}{0.33\textwidth}
        \centering
        \includegraphics[width=\textwidth]{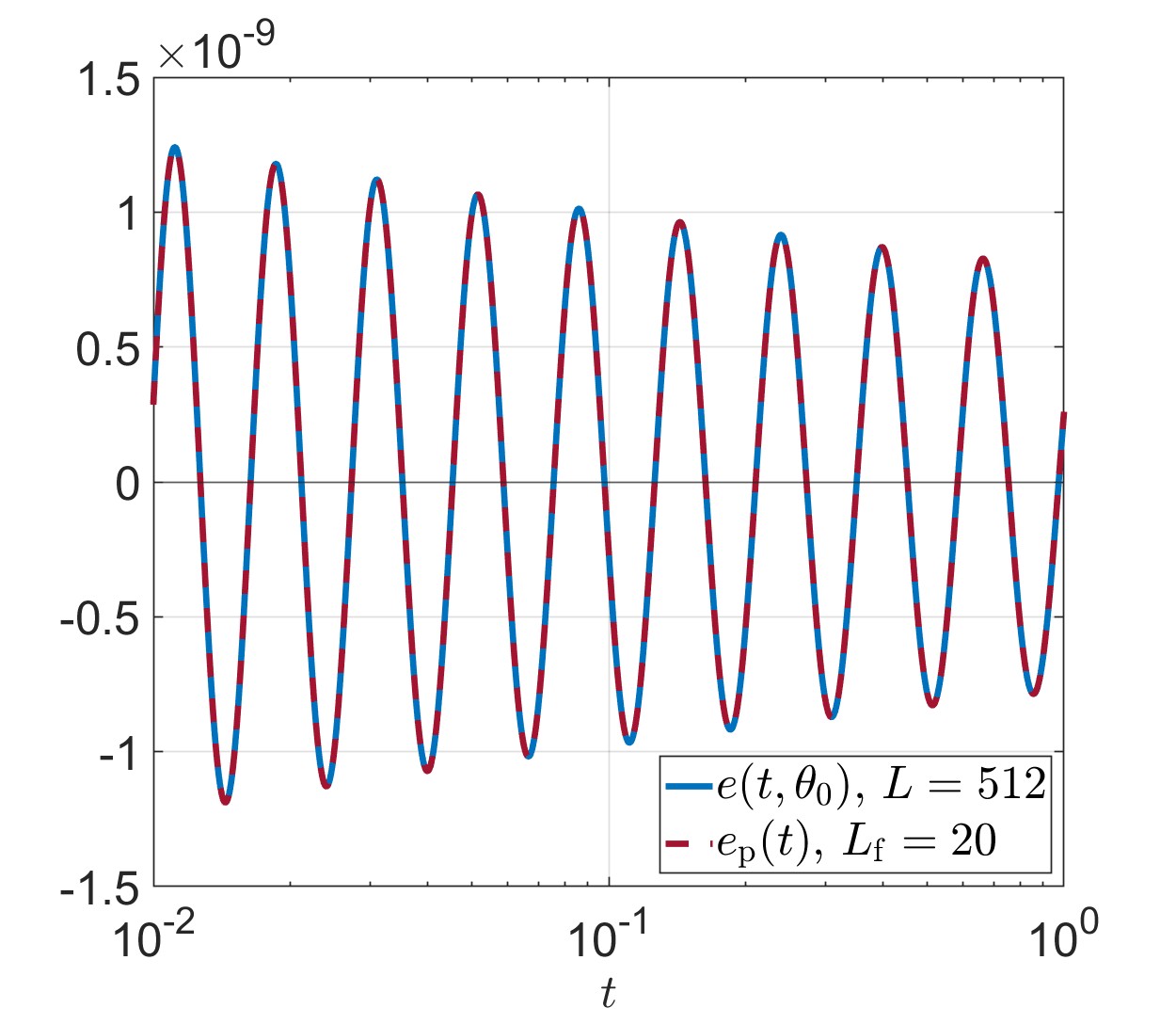}
        \caption{$\alpha = 0.9$}
        \label{fig:4_f}
    \end{subfigure}\hfill
    \caption{Comparison of the difference between the exact and approximation values of the power law kernel $t^{\alpha-1}$ versus $t$ before and after the application of Prony's method for various values of $L$ and $\alpha \in \{0.1, 0.5, 0.9\}$ over the interval $[10^{-2}, 1]$. Here, the value of $\epsilon$ is set to $10^{-10}$ and the horizontal axis is shown on a logarithmic scale.}
    \label{fig:4}
\end{figure}

For a more precise demonstration of the performance of the Prony method, we visualize in Figure~\ref{fig:4} a comparison of the difference between the exact and approximate values of the power law kernel before the application of Prony's method, i.e.\ $e(t, \theta_0)$, vs.\ $t$ (plotted as continuous lines) with the analog error after using Prony's method, viz.\ $e_{\text p}(t)$, vs.\ $t$ (dashed) for various values of $L$ and $\alpha \in \{0.1, 0.5, 0.9\}$, where $t$ ranges over the interval $[10^{-2}, 1]$. As we see, for smaller values of $L$, with Prony's method the difference grows as $t$ increases, but the maximum absolute error remains constant. As $L$ increases, the difference exhibits the same behaviour as observed before Prony's method. 

\begin{remark}
	\label{rem:stability}
	To ensure stability in Prony's reduction algorithm, the Hankel matrix $\mathbf{H}$ must be sufficiently far away from a singular matrix. Based on our numerical observations, for $T=1$, this condition holds, making the algorithm stable. However, as $T$ increases, the Hankel matrix may become singular in some cases, causing algorithm instability. To address this situation in Prony's method, the rescaling approach discussed in Remark \ref{rem:rescaling} can be applied.
	Therefore, we first compute the exponents and weights for the interval $[\delta/T, 1]$, then apply Prony's algorithm, and finally rescale as in eq. (\ref{rescale}) to create an approximation for the interval $[\delta, T]$ with a reduced number of terms. For more details on the maximum error values in this case, see Appendix.
\end{remark}

In Algorithm \ref{alg:prony}, we now present our algorithm of Prony's method formally in a pseudo-code style.
\begin{algorithm}[h!t]
\small
\caption{Prony method for kernel approximation given in a MATLAB-like notation, i.e.\ for an array $X$ we denote its components by $X(1), X(2), \ldots$. \label{alg:prony}}
\begin{algorithmic}[1]
\State \textbf{Input:} $\alpha$, $\delta$, $T$, $N$, $L$, $\epsilon$
\State \textbf{Output:} maximal approximation error after reduction by Prony method
\State Generate time vector $t_j$ for $j = 1, \dots, N$
\State Compute bounds:
\State \quad $l_{\min} \gets \min(\ln(\epsilon/T), \ln(\epsilon(1-\alpha))/(1-\alpha))$
\State \quad $l_{\max} \gets \ln(\ln(\epsilon^{-1})/\delta)$
\State Set step size $h \gets (l_{\max} - l_{\min})/(L-1)$
\State Initialize vectors $w^0$, $b^0$, $\omega$ of length $L$
\For{$i = 1$ to $L$}
    \State $\omega(i) \gets l_{\min} + (i-1)h$
    \State $b^0(i) \gets -\exp(\omega(i))$
    \State $w^0(i) \gets \begin{cases} 
        \frac{h}{2} \exp((1-\alpha)\omega(i)), & \text{if } i = 1 \text{ or } i = L \\
        h \exp((1-\alpha)\omega(i)), & \text{otherwise}
    \end{cases}$
\EndFor
\State Set $M \gets$ number of $\omega(i) \leq 0$
\State Combine parameters: $\theta_0 \gets [w^0, b^0]$
\State Compute exact kernel: $S(t_j) \gets t_j^{\alpha-1}$ for $j = 1, \dots, N$
\State Compute initial approximation: 
\State \quad $f(t_j, \theta_0) \gets \frac{1}{\Gamma(1-\alpha)} \sum_{i=1}^L \theta_0(i) \exp(\theta_0(L+i)t_j)$
\State Compute initial error: $e(t_j, \theta_0) \gets S(t_j) - f(t_j, \theta_0)$
\State Set $\epsilon' \gets \max_{j=1,\dots,N} |e(t_j, \theta_0)|$
\State 
\State Determine suitable $K$ (number of Prony terms) and $L_{\text p}$
\For{$i = 1$ to $K$}
    \State Compute Prony weights $\rho(i)$ and exponents $\eta(i)$ as described at the beginning of this section
\EndFor
\State Set $w \gets [\rho(1), \dots, \rho(K), w^0(L_{\text p}+1), \dots, w^0(L)]$
\State Set $b \gets [\eta(1), \dots, \eta(K), b^0(L_{\text p}+1), \dots, b^0(L)]$
\State Combine parameters: $\theta_{\text p} \gets [w, b]$
\State Set $L_{\text f} \gets L - L_{\text p} + K$
\State Compute approximation error after reduction : $e_{\text p}(t_j) \gets S(t_j) - \frac{1}{\Gamma(1-\alpha)} \sum_{i=1}^{L_{\text f}} \theta_{\text p}(i) \exp(\theta_{\text p}(L_{\text f}+i)t_j)$
\State Compute maximal error after reduction: $\max_{j=1,\dots,N} |e_{\text p}(t_{j})|$
\end{algorithmic}
\end{algorithm}

The algorithm inputs include data \(\alpha\) (fractional order), \(\delta\) (initial time), \(T\) (final time), \(N\) (number of nodes for the discretization of the time interval $[\delta, T]$),  \(L\) (number of exponential terms), and \(\epsilon\) (a very small threshold for neglecting the integrals (\ref{3})). This algorithm is split into two parts:
\begin{itemize}
\item \textbf{Part 1 (Initialization)---lines 1--20:} initializes the input data, then computes the lower and upper bounds of the integral $(\ref{KerExp1})$ based on equations $(\ref{8})$ and $(\ref{lmin_approx2})$ and initial values of the weights and exponents $ \theta_{0}= [w_{0}, b_{0}] \in \mathbb{R}^{2L}$, where \( w_{0}= [w_1^0, \dots, w_L^0] \) and \( b_{0} = [b_1^0, \dots, b_L^0] \in \mathbb{R}^L \) by $(\ref{weightandnode})$, and finally evaluates the initial approximation error of the kernel function, $e(t,\theta_{0})$ using $(\ref{Abserror})$ and its maximum error $\epsilon'$.
\item \textbf{Part 2 (Application of Prony's method)---lines 22--31:} applies the Prony's method to the initial approximation $ \theta_{0}$  to reduce the total number of exponential terms while preserving the accuracy $\epsilon'$.
\end{itemize}

\begin{remark}
One can choose a uniform discretization of the interval $[\delta,T]$ on which the approximation errors are computed. However, to reduce the computation time, especially for large $T$, we can use a geometric grid in the interval $[\delta,T]$ given by $t_{j} = \delta(T/\delta)^{(j-1)/(N-1)}$ for $j=1,\ldots,N$. In that case, the common ratio between consecutive points is $r = \left( T/\delta \right)^{1/(N-1)}$.
\end{remark}

\section{Application to Fractional Ordinary Differential Equations}\label{Sec3}

Next, our aim is to present two algorithms for the numerical solution of the fractional ordinary differential equation (FODE)
\begin{equation}\label{FODE1}
D^{\alpha}y(t) = f(t,y(t)), \quad t \in [0,T],
\end{equation}
subject to the initial condition $y(0) = y_{0}$,
where $ 0<\alpha <1$,  $T$ is a real positive number, and $D^{\alpha}$ is the Caputo type fractional derivative \cite[Chapter 3]{Diethelm2010} of order $\alpha$, which is defined by
\begin{equation}\label{FODE2}
D^{\alpha}y(t) = I^{1-\alpha} y'(t),
\end{equation}
where $I^{1-\alpha}$ is the Riemann-Liouville (RL) fractional integral of order $1-\alpha$ defined by $(\ref{RLIntegral})$ which leads to 
\begin{equation}\label{FODE4}
D^{\alpha}y(t) =\frac{1}{\Gamma(1-\alpha)}\int_{0}^{t}(t-\tau)^{-\alpha}y'(\tau)\, \mathrm d\tau. 
\end{equation}
The Caputo type fractional derivative $D^\alpha$ of order $\alpha$ is a left inverse of the RL integral $I^\alpha$ of order $\alpha$, that is, $D^{\alpha}I^{\alpha}y = y$, but it is not a right inverse since
\begin{equation}\label{FODE5}
I^{\alpha}D^{\alpha}y(t) = y(t) - y(0).
\end{equation}
By applying the RL fractional integral of order $\alpha$ to both sides of Equation (\ref{FODE1}), and considering Equation (\ref{FODE5}), we get the Volterra integral equation
\begin{equation}\label{FODE6}
y(t)=y(0)+\frac{1}{\Gamma(\alpha)}\int_{0}^{t}(t-\tau)^{\alpha-1}f(\tau,y(\tau))\, \mathrm d\tau,
\end{equation}
see \cite[\S~6.1]{Diethelm2010}.
To numerically solve (\ref{FODE6}), we first introduce an arbitrary discretization (e.g., a uniform grid), which is defined as $0 = t_{0} < t_{1}< \ldots < t_{N}=T$,  where the step sizes are \( h_j = t_{j }- t_{j-1} \) for \( j = 1, \ldots, N \), and \( h_j \) may vary. Here, we assume that the approximations $y_{j} \approx y(t_{j})$ for $j=1,2,\ldots,n-1$ and some $n \geq 1$ have been already computed. Then, we compute the approximation $y_{n} \approx y(t_{n})$ using the equation
\begin{equation}\label{FODE7}
	y_{n} = y(0) + \tilde I^{\alpha} f(t_{n},y_{n}),
\end{equation}
where $\tilde I^{\alpha} f(t_{n},y_{n})$ is a suitable approximation for the RL integral
\begin{equation}\label{FODE8}
	I^\alpha [f(\cdot, y(\cdot))] (t_n) = \frac{1}{\Gamma(\alpha)} \int_{0}^{t_{n}}(t_{n}-\tau)^{\alpha-1}f(\tau,y(\tau))\, \mathrm d\tau
\end{equation}
computed with the help of the algorithm described in Sections \ref{sec:approx-exp-sum} and \ref{Sec2}.
Following the approach outlined in \cite{rebecca2010, chaudhary2024novel}, we split up the integral in equation (\ref{FODE8}) into a local part $L^{\alpha}_n$ and a history part $H^{\alpha}_n$ according to
\begin{equation}\label{HL1}
	I^\alpha [f(\cdot, y(\cdot))] (t_n) = H^{\alpha}_n [f(\cdot, y(\cdot))]  + L^{\alpha}_n [f(\cdot, y(\cdot))],
\end{equation}
where the history part is given by
\begin{equation}\label{HL3}
	H^{\alpha}_n [f(\cdot, y(\cdot))] = \frac{1}{\Gamma(\alpha)} \int_{0}^{t_{n-1}} (t_{n}-\tau)^{\alpha-1} f(\tau, y(\tau)) \, \mathrm{d}\tau
\end{equation}
and the local part is defined by
\begin{equation}\label{HL2}
	L^{\alpha}_n [f(\cdot, y(\cdot))]  = \frac{1}{\Gamma(\alpha)} \int_{t_{n-1}}^{t_{n}} (t_{n}-\tau)^{\alpha-1} f(\tau, y(\tau)) \, \mathrm{d}\tau.
\end{equation}
Here, to approximate (\ref{HL3}) and (\ref{HL2}), we present two algorithms.

\subsection{Method of Constant Interpolation} 
In the first approach, to approximate the integrand $f$ in (\ref{HL3}) and (\ref{HL2}), we use piecewise constant interpolation $f (\tau,y(\tau)) \approx f(t_{j},y_{j})$ for all $\tau \in (t_{j-1},t_{j})$ and all $j$.  

In the local part $L^{\alpha}_n$ the integration ranges from $t_{n-1}$ to $t_{n}$, i.e., it does not involve contributions from the process history. Due to the singularity of the kernel function $(t_n - \tau)^{\alpha-1}\) at \(\tau = t_n$, there is a region near $\tau \approx t_n$ where the exponential sum approximation loses accuracy. Therefore, to compute the local contribution,  we need to use standard numerical methods that avoid relying on the exponential sum approximation, achieving accurate results with $O(1)$ time complexity and $O(1)$ memory usage. Thus, we compute $L^{\alpha}f(t_{n},y_{n})$ as
\begin{equation}\label{Local1}
	L^\alpha_n [f(\cdot, y(\cdot))] \approx \frac{1}{\Gamma(\alpha)} \int_{t_{n-1}}^{t_{n}} (t_{n}-\tau)^{\alpha-1} f(t_n, y_n) \, \mathrm{d}\tau =  \frac{h_{n}^{\alpha}}{\Gamma(\alpha+1)}f_{n}, \qquad n=1, \ldots, N,
\end{equation}
where 
\[
	f_{j} = f(t_{j},y_{j}) \quad \text{and} \quad h_{j} = t_{j}-t_{j-1} \qquad \text{ for all } j.
\]
Here, we have used the functional equation $\alpha \Gamma(\alpha) = \Gamma(\alpha+1)$ \cite{artin2015gamma}.

In the history part $H^{\alpha}_n$, the integration ranges from $0$ to $t_{n-1}$, meaning the history part is confined to times up to at most $t_{n-1}$ and thus it (a) avoids the kernel function's singularity at $\tau = t_{n}$, and (b) becomes zero for $n=1$. As a result, the exponential sum technique can be used to efficiently compute an approximation for $H^{\alpha}_n [f(\cdot, y(\cdot))]$. By considering the time discretization, we obtain
\[
	H^{\alpha}_n [f(\cdot, y(\cdot))] = \sum_{j=1}^{n-1} \frac{1}{\Gamma(\alpha)} \int_{t_{j-1}}^{t_j} (t_n - \tau)^{\alpha - 1} f(\tau, y(\tau))\, \mathrm d\tau, \quad n=2, \ldots, N.
\]
Using the piecewise constant interpolant $ f(\tau, y(\tau)) \approx f(t_{j},y_{j}) $ over each interval \( [t_{j-1}, t_j] \) for \( j=1,\dots,n-1 \), and considering the approximation of the kernel function $(t_n - \tau)^{\alpha - 1}$ for $ 0 \leq \tau \leq t_{n-1}$ as given in (\ref{TrapExpFin}) and applying it to the sum, we get
\begin{eqnarray}\label{FODE13}
	H^{\alpha}_n [f(\cdot, y(\cdot))] &\approx& \sum_{j=1}^{n-1}\frac{f_{j}}{\Gamma(\alpha)}  \int_{t_{j-1}}^{t_{j}}(t_{n}-\tau)^{\alpha-1}\, \mathrm{d}\tau  \approx \sum_{j=1}^{n-1}\frac{f_{j}}{\Gamma(\alpha)} \int_{t_{j-1}}^{t_{j}}\sum_{l = 1}^{L} \frac{1}{\Gamma(1-\alpha)}w_{l}\mathrm e^{b_{l}(t_{n}-\tau)}\, \mathrm{d}\tau \nonumber \\
	&=& \sum_{j=1}^{n-1}c_{\alpha} w_{l}f_{j}\sum_{l = 1}^{L} \int_{t_{j-1}}^{t_{j}}\mathrm e^{b_{l}(t_{n}-\tau)}\, \mathrm{d}\tau = \sum_{l = 1}^{L}\Phi_{l}^{n},
\end{eqnarray}
where  $ c_{\alpha} = \frac 1 {\Gamma(\alpha) \Gamma(1-\alpha)} = \frac{\sin \pi\alpha}{\pi}$ in view of the reflection formula (one of the fundamental properties of the Gamma function \cite{artin2015gamma}) and 
\begin{equation}\label{FODE14}
	\Phi_{l}^{n} = \sum_{j=1}^{n-1}K_{l,n,j}f_{j}
\end{equation}
with
\begin{equation}\label{FODE15}
	K_{l,n,j} =c_{\alpha}w_{l}\int_{t_{j-1}}^{t_{j}}\mathrm e^{b_{l}(t_{n}-\tau)}\, \mathrm{d}\tau
\end{equation}
It should be noted that the value of $\delta$ used in the construction of the sum-of-exponentials approximation for the kernel function must be less than or equal to the smallest time step $\min_j h_j = \min_j (t_j - t_{j-1})$. This ensures that $\delta \leq h_n = t_n - t_{n-1} \leq t_n - \tau \leq T$ holds for all $0 \leq \tau \leq t_{n-1}$, which implies $\delta \leq t_n - \tau \leq T$.

The crucial property of this approach is that we can find a relation between $K_{l,n,j}$ and $K_{l,n-1,j}$, namely
\begin{eqnarray}\label{FODE17}
	K_{l,n,j} &=& c_{\alpha}w_{l} \int_{t_{j-1}}^{t_{j}}\mathrm e^{b_{l}(t_{n}-\tau)}\, \mathrm{d}\tau = c_{\alpha}w_{l} \int_{t_{j-1}}^{t_{j}}\mathrm e^{b_{l}(t_{n-1}-\tau+t_{n}-t_{n-1})}\, \mathrm{d}\tau \nonumber \\ 
	&=& \mathrm e^{b_{l}h_{n}}c_{\alpha}w_{l}\int_{t_{j-1}}^{t_{j}}\mathrm e^{b_{l}(t_{n-1}-\tau)}\, \mathrm{d}\tau =  \mathrm e^{b_{l}h_{n}} K_{l,n-1,j}.
\end{eqnarray}
Thus, in this evaluation scheme, we derive the recursive relation
 \begin{eqnarray}\label{FODE18}
\Phi_{l}^{n} = \sum_{j=1}^{n-1}K_{l,n,j}f_{j}&=&K_{l,n,n-1}f_{n-1}+\sum_{j=1}^{n-2}K_{l,n,j}f_{j} = K_{l,n,n-1}f_{n-1}+ \mathrm e^{b_{l}h_{n}} \sum_{j=1}^{n-2}K_{l,n-1,j} f_j \nonumber \\ 
&=&  K_{l,n,n-1}f_{n-1}+\mathrm e^{b_{l}h_{n}}\Phi_{l}^{n-1}, \qquad (n \ge 2) 
\end{eqnarray}
with the initialization $\Phi_{l}^{1} =0$. Using eq.~\eqref{FODE18} instead of eq.~\eqref{FODE14}, we can reduce the computational complexity of the evaluation of $\Phi_l^n$ (required in the $n$th time step) from $O(n)$ to $O(1)$, and therefore the overall cost for taking steps $n = 1, 2, \ldots, N$ goes down from $O(N^2)$ to $O(N)$. The approximate solution to the fractional differential equation at \( t_n \) is thus computed as
\begin{equation}\label{FODE19}
	y_n = y_0 + \frac{h_{n}^{\alpha}}{\Gamma(\alpha+1)}f_{n} + \sum_{l=1}^L \Phi_l^n
	= y_0 + \frac{h_{n}^{\alpha}}{\Gamma(\alpha+1)}f(t_{n}, y_n) + \sum_{l=1}^L \Phi_l^n .
\end{equation}
Here, we estimate the value of $y(t)$ at time $t_n$ using only the information from the immediately preceding time point $t_{n-1}$, without relying on any data from the earlier history. Thus, this approach can also efficiently reduce the active storage requirements because we do not need to keep any history data in memory for the later time steps..

Next we need to solve eq.~\eqref{FODE19} for the unknown variable $ y_n $ that appears on both sides. But the nonlinearity of the function $ f $ makes it generally impossible to find $y_{n}$ directly. Thus, we have to resort to an iterative method for solving this equation. Here we use the iterative Newton-Raphson method, which starts with an initial estimate $y_{n}^{(0)}$ for $y_{n}$ and then computes the successive approximations of $y_{n}$ using the following equation
\[
	y_n^{(s+1)} = y_n^{(s)} - \frac{G(y_n^{(s)})}{G'(y_n^{(s)})}, \quad s = 0,1,2,\dots,
\]
where
\begin{equation}
	G(x) = x - y_0 -\frac{h_{n}^{\alpha}}{\Gamma(\alpha+1)}f(t_n, x) - \sum_{l=1}^L \Phi_l^n 
	\quad \text{such that} \quad 
	G'(x) = 1 - \frac{h_{n}^{\alpha}}{\Gamma(\alpha+1)} \frac{\partial f}{\partial x}(t_n, x).
\end{equation}
Here, $ y_n^{(s)} $ is the $ s $-th approximation of $ y_n $, with the initial guess $ y_n^{(0)} = y_{n-1} $. This is a common approach---it is often effective to begin Newton iterations with the most recently calculated approximation. We then terminate our iterations when two successive values differ by less than a specified tolerance.

\begin{remark}
	If the Newton iteration does not converge (which may happen, e.g., if the initial guess $y_n^{(0)}$ is not good enough), or if the computation of the function $G'$ is infeasible or too expensive, then one can also use a fixed point iteration directly based on eq.~\eqref{FODE19}, i.e.\ one (as above) starts with $y_n^{(0)} = y_{n-1}$ and then proceeds according to 
	\[
		y_n^{(s+1)} = y_0 + \frac{h_{n}^{\alpha}}{\Gamma(\alpha+1)}f(t_n, y_n^{(s)}) + \sum_{l=1}^L \Phi_l^n, \qquad s=0, 1, 2, \ldots,
	\]
	again terminating when two successive values are sufficiently close to each other and then accepting the last one as $y_n$. Typically, this procedure will require more iterations to converge than the Newton-Raphson method, but each step is computationally cheaper.
\end {remark}

It should be noted that we can easily compute $K_{l,n,n-1}$ formally using (\ref{FODE15}) which yields
\begin{equation}\label{FODE20}
K_{l,n,n-1} = \frac{c_{\alpha}w_{l}}{b_{l}} \big[\mathrm e^{b_{l}(t_{n}-t_{n-2})}-\mathrm e^{b_{l}(t_{n}-t_{n-1})}\big].
\end{equation}
As noted in Section \ref{Sec2}, before applying Prony's method, we are likely to have exponents $b_l$ near zero, which causes a cancellation of significant digits or even a division by zero when using (\ref{FODE20}) naively and thus may be fatal when implementing the algorithm in finite-precision arithmetic. However, Prony's method can generate negative exponents, which are not very close to zero, particularly for smaller intervals, enabling the use of equation (\ref{FODE20}) without causing any problems, and for very large intervals, where very small exponents may exist, exact values of $K_{l,n,n-1}$ by (\ref{FODE20}) can then be replaced by numerical computation of $K_{l,n,n-1}$ using integration methods (see Table \ref{tab:4} in Appendix).

\subsection{Method of First Order ODE}
Here, for the local part, we use a linear interpolant as the straight line between the points $(t_{n},f_{n})$ and $(t_{n-1},f_{n-1})$. So, we get the approximation  $f( \tau,y(\tau)) \approx (f_{n}-f_{n-1})(\tau-t_{n-1})/h_n + f_{n-1}$ for all $\tau \in (t_{n},t_{n-1})$. By substituting this into (\ref{HL2}) and after some calculation, we get
\begin{equation}\label{Local2}
	L^\alpha_n [f(\cdot, y(\cdot))]  \approx  \frac{h_{n}^{\alpha}}{\Gamma(\alpha+2)}(\alpha f_{n-1}+f_{n}), \quad n=1, \ldots, N.
\end{equation}
For the history part, after considering the approximation of the kernel function, we have
\begin{equation}\label{FODE21}
	H^\alpha_n [f(\cdot, y(\cdot))]  \approx c_{\alpha}\sum_{l = 1}^{L}w_{l}\mu(t_{n},h_{n},b_{l}),
\end{equation}
where
\begin{equation}\label{FODE22}
	\mu(t,h,b_{l}) = \int_{0}^{t-h}\mathrm e^{b_{l}(t-\tau)}f(\tau,y(\tau))\, \mathrm{d}\tau. 
\end{equation}
It is easy to show that for any $b_{l}$ with $l=1,\ldots,L$, the function $\mu(\cdot,h,b_{l})$ is the unique solution to the first order differential equation
\begin{equation}\label{FODE23}
	\frac{\partial \mu}{\partial t}(t,h,b_{l}) = b_{l}\mu(t,h,b_{l}) +\mathrm e^{b_{l}h}f(t-h,y(t-h))  \qquad \text{for }  t \in (h, T]
\end{equation}
subject to the initial condition $\mu(h,h,b_{l}) = 0$.

Now, we require a method for the approximate solution to the initial value problem (\ref{FODE23}). Since all the $b_{l}$ are negative and some of them have large absolute values, i.e.\ $b_{l}\ll 0$, the differential equation (\ref{FODE23}) becomes very stiff, posing challenges for numerical solutions. Therefore, an A-stable scheme is necessary, and preferably one should even choose an L-stable method.  Implicit methods are preferred to achieve reliable results \cite{plato2003concise}. Here, we will test both the backward Euler method (which is L-stable)  and the trapezoidal rule (that is not L-stable) for this purpose. Notably, equation (\ref{FODE23}) indicates that the initial value problem must be solved numerically $L$ times, corresponding to each parameter value $b_l$ ($l = 1, 2, \dots, L$).

Thus, for $n \geq 1$, we get the backward Euler scheme, which is of first order
\begin{equation}\label{FODE24}
	\mu(t_{n},h_{n},b_{l}) = \frac{\mu(t_{n-1},h_{n-1},b_{l}) + h_{n}\mathrm e^{b_{l}h_{n}}f_{n-1}}{1-h_{n}b_{l}}, 
\end{equation}
and the trapezoidal rule, which is of second order
\begin{equation}\label{FODE25}
	\mu(t_{n},h_{n},b_{l}) = \frac{\mu(t_{n-1},h_{n-1},b_{l})\Big[1+\frac{h_{n}b_{l}}{2}\Big]+\frac{h_{n}}{2}\Big[\mathrm e^{b_{l}h_{n}}f_{n-1}+\mathrm e^{b_{l}h_{n-1}}f_{n-2}\Big]}{1-\frac{h_{n}b_{l}}{2}}, 
\end{equation}
with $\mu(t_{1},h_{1},b_{l}) = \mu(h_{1},h_{1},b_{l}) = 0$. Thus, the approximate solution to the fractional differential equation (\ref{FODE1}) at \( t_n \) is
\begin{eqnarray}
	\nonumber
	y_{n} &=& y_{0} +  \frac{h_{n}^{\alpha}}{\Gamma(\alpha+2)}(\alpha f_{n-1}+f_{n}) + c_{\alpha}\sum_{l = 1}^{L}w_{l}\mu(t_{n},h_{n},b_{l}) \\ 
	\label{FODE26}
	 &=& y_{0} +  \frac{h_{n}^{\alpha}}{\Gamma(\alpha+2)}(\alpha f(t_{n-1}, y_{n-1})+f(t_{n}, y_n)) + c_{\alpha}\sum_{l = 1}^{L}w_{l}\mu(t_{n},h_{n},b_{l})
	\qquad \text{for } n \geq 1.
\end{eqnarray}
Similar to the method of constant interpolation, solving (\ref{FODE26}) for $y_{n}$ is here done by Newton's iteration method or a direct fixed point iteration.

\section{Numerical Results}
\label{sec:num}

Here we present some numerical examples of fractional differential equations and solve them using the methods introduced in Section \ref{Sec3} and our new optimization algorithm for approximating the kernel function by exponential sums. Our algorithm has been implemented in MATLAB R2023b on a notebook with an AMD Ryzen 7 PRO 7840U CPU clocked at 3.30 GHz running Windows 11. In all examples, the values of $\epsilon$ and the tolerance for the Newton iterations are set to $10^{-10}$ and for the kernel approximation, $\delta$ is set to $10^{-5}$. For Examples \ref{ex1} and \ref{ex2}, we assume that the discretization of the interval $[0,T]$ is uniform. i.e., $t_{n} = nh,$ where $h = T/N$ is the step size and $N \in \mathbb{N}$.

\begin{example} 
	\label{ex1}
Our first example for testing our numerical method is the nonlinear problem
\begin{equation}\label{example_1}
D^{\alpha}y(t)=\frac{40320}{\Gamma(9-\alpha)}t^{8-\alpha} - 3 \frac{\Gamma(5+\alpha/2)}{\Gamma(5-\alpha/2)}t^{4-\frac{\alpha}{2}}+\frac{9}{4}\Gamma(\alpha+1)+(\frac{3}{2}t^{\alpha/2}-t^{4})^3-[y(t)]^{\frac{3}{2}},
\quad t \in [0,1],
\end{equation}
with the initial condition $y(0) = 0$ originally  proposed in \cite{diethelm2004detailed}. The exact solution of this equation is $y(t) = t^{8} - 3t^{4+\frac{\alpha}{2}}+\frac{9}{4}t^{\alpha}$. 
\end{example}

\begin{example}
	\label{ex2}
For our second example, we consider the following fractional differential equation
\begin{equation}\label{example_2}
	D^{\alpha}y(t)=\lambda y(t), \qquad y(0) = 1,
\end{equation}
whose exact solution is $y(t) = E_{\alpha}(\lambda t^{\alpha})$, with $E_{\alpha}$ being the Mittag-Leffler function of order $\alpha$, which can be numerically evaluated by an algorithm presented in \cite{garrappa2015numerical}. Here, we numerically solved the problem in the interval $[0, 10]$ with various values of $\alpha $ and $\lambda = -1$.
\end{example}

\begin{figure}[H]
    \centering
    \begin{subfigure}{0.33\textwidth}
        \centering
        \includegraphics[width=\textwidth]{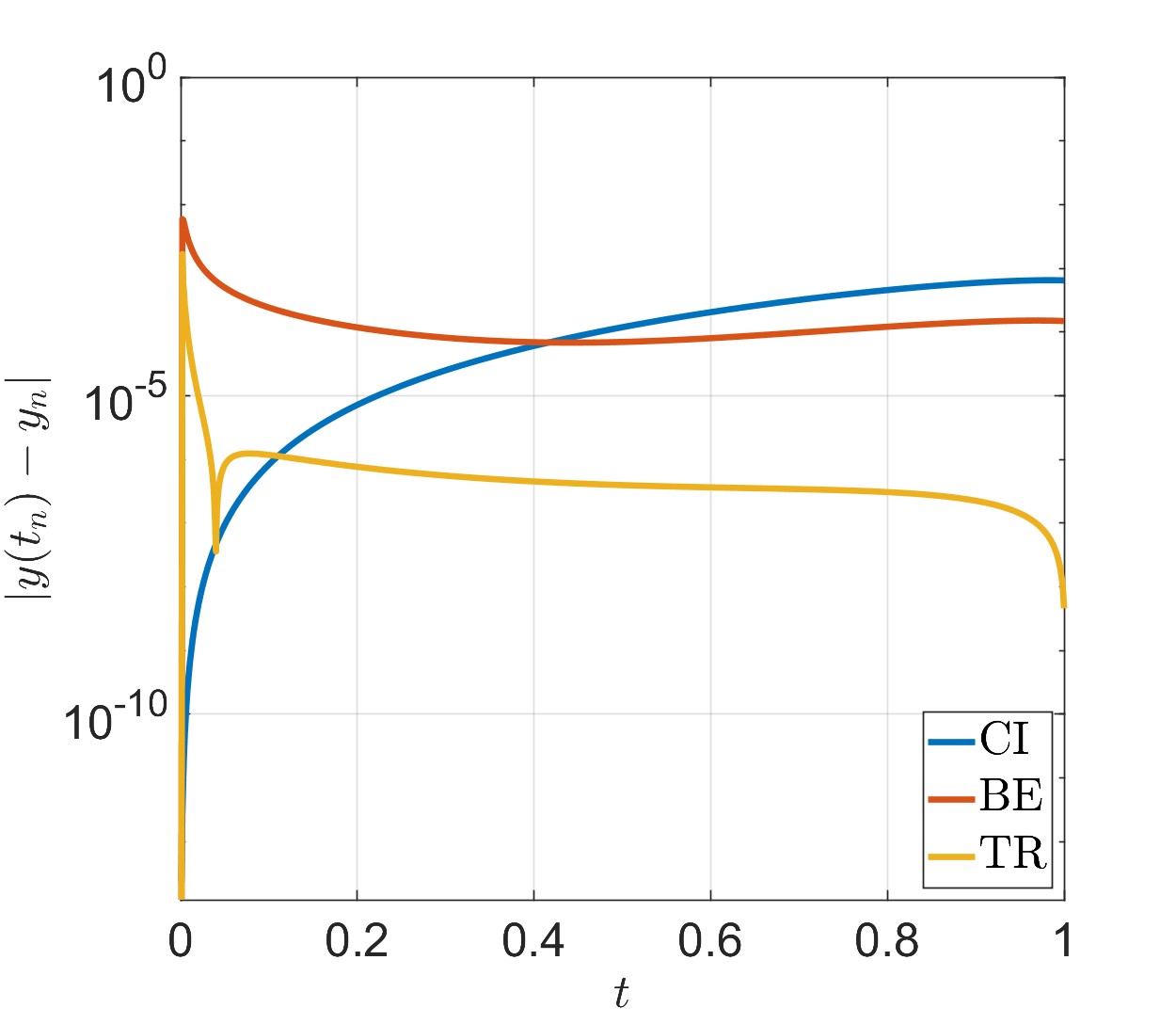}
        \caption{$\alpha = 0.1, L=128, L_{\text f} = 50$}
        \label{fig:example_1_01}
    \end{subfigure}\hfill
    \begin{subfigure}{0.33\textwidth}
        \centering
        \includegraphics[width=\textwidth]{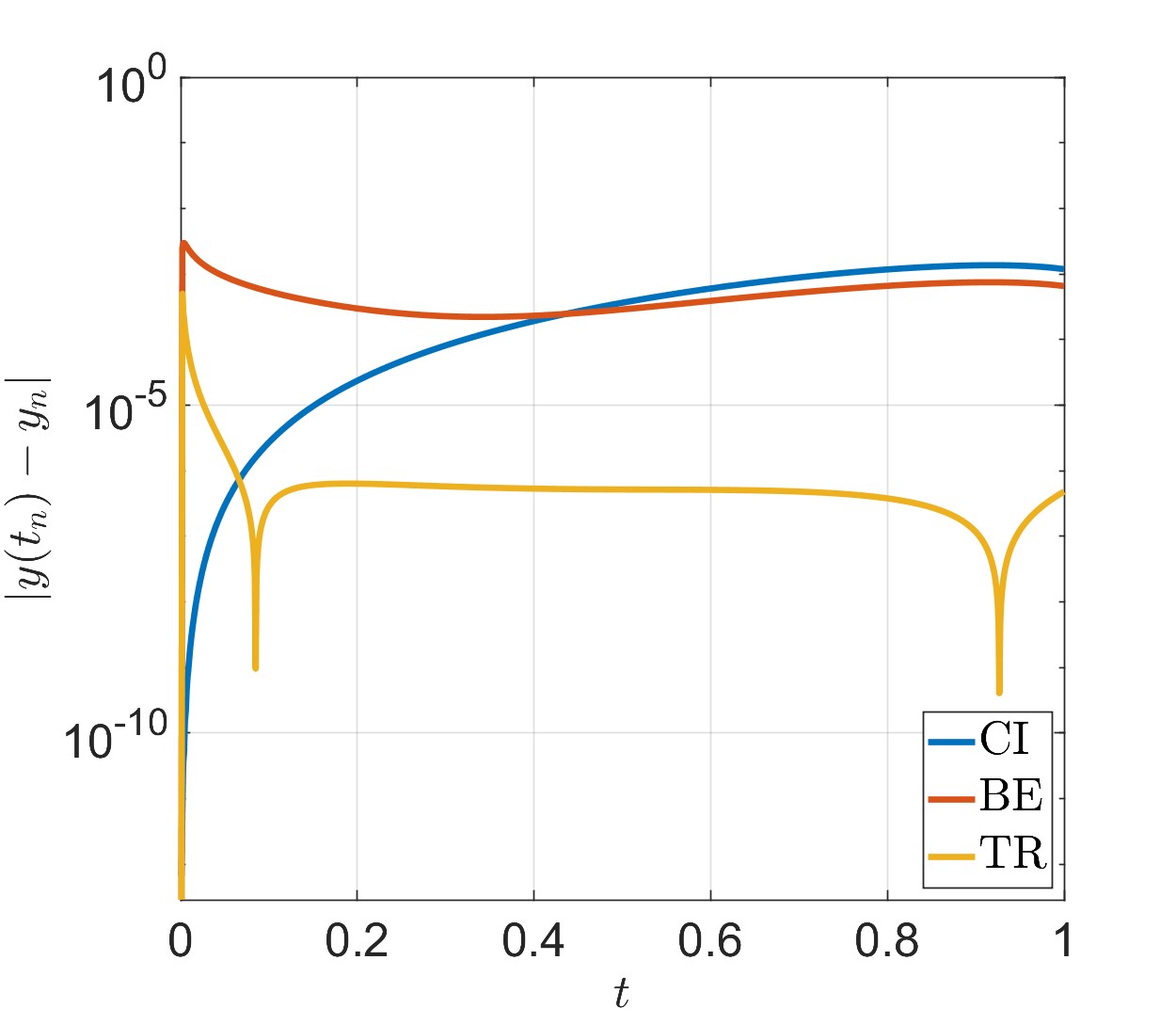}
        \caption{$\alpha = 0.5, L=128, L_{\text f} = 33$}
        \label{fig:example_1_05}
    \end{subfigure}\hfill
    \begin{subfigure}{0.33\textwidth}
        \centering
        \includegraphics[width=\textwidth]{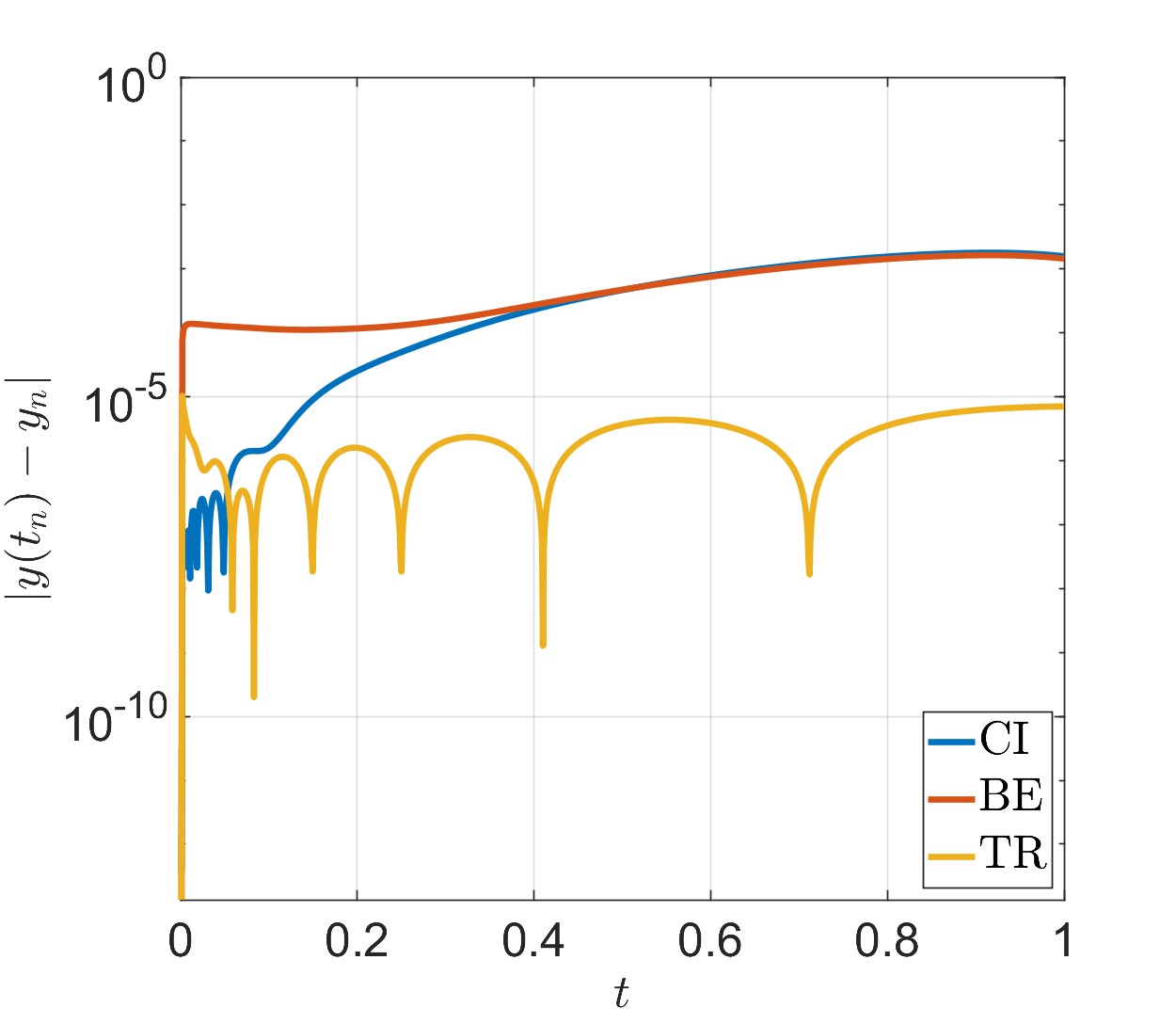}
        \caption{$\alpha = 0.9, L=256, L_{\text f} = 17$}
        \label{fig:example_1_09}
    \end{subfigure}\hfill
    \caption{Comparison of the absolute errors for (\ref{example_1}) as a function of $t$, which are obtained by the methods of constant interpolation (CI) and first order ODE by backward Euler (BE) scheme and trapezoidal rule (TR) for different values of the fractional order $\alpha$. In all cases, the time step size is set to $2^{-10}$. The $y$-axis is plotted in a logarithmic scale.}
    \label{fig:5}
\end{figure}

\begin{figure}[H]
    \centering
    \begin{subfigure}{0.33\textwidth}
        \centering
        \includegraphics[width=\textwidth]{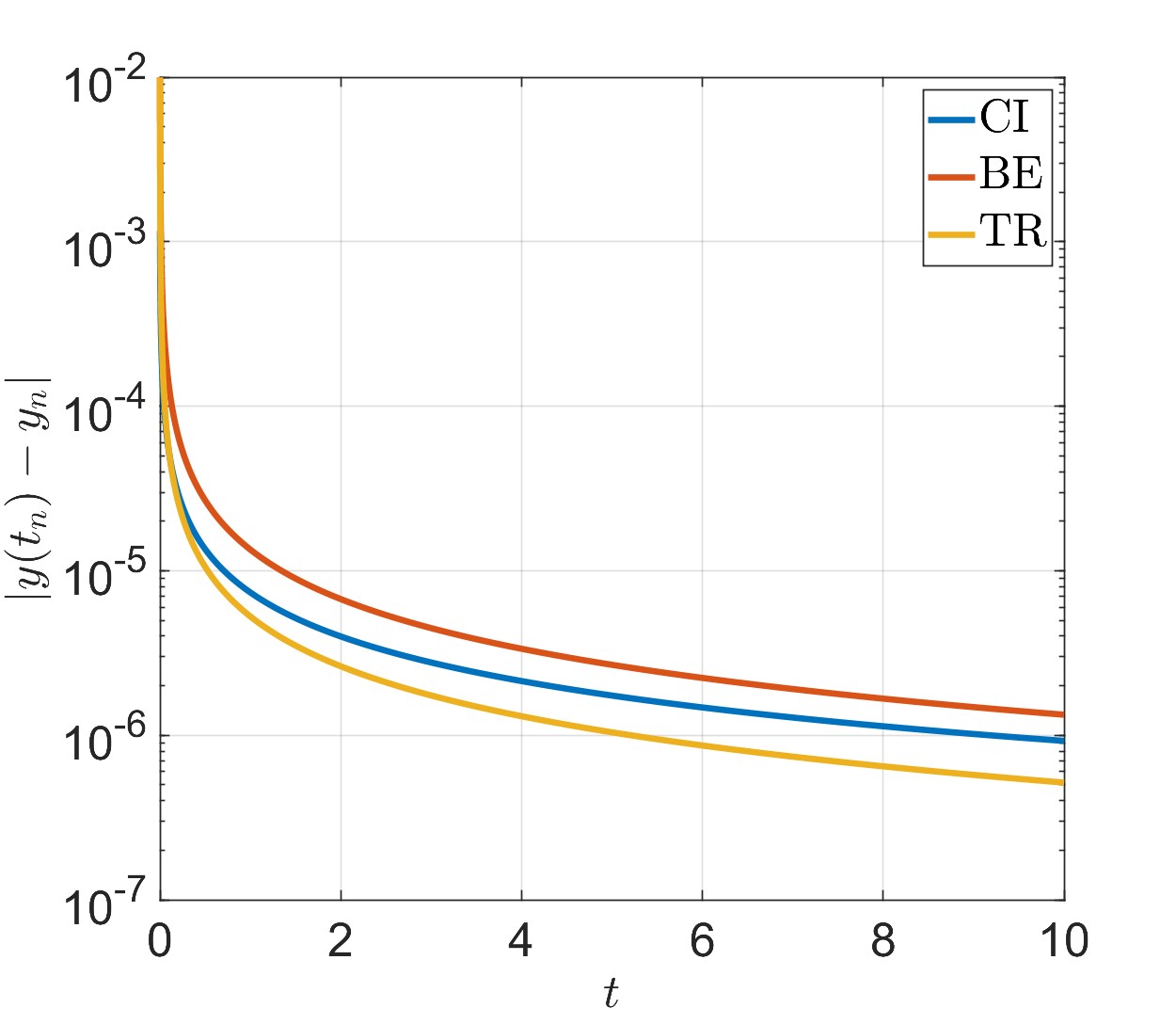}
        \caption{$\alpha = 0.1, L=128, L_{\text f} = 54$}
        \label{fig:example_2_01}
    \end{subfigure}\hfill
    \begin{subfigure}{0.33\textwidth}
        \centering
        \includegraphics[width=\textwidth]{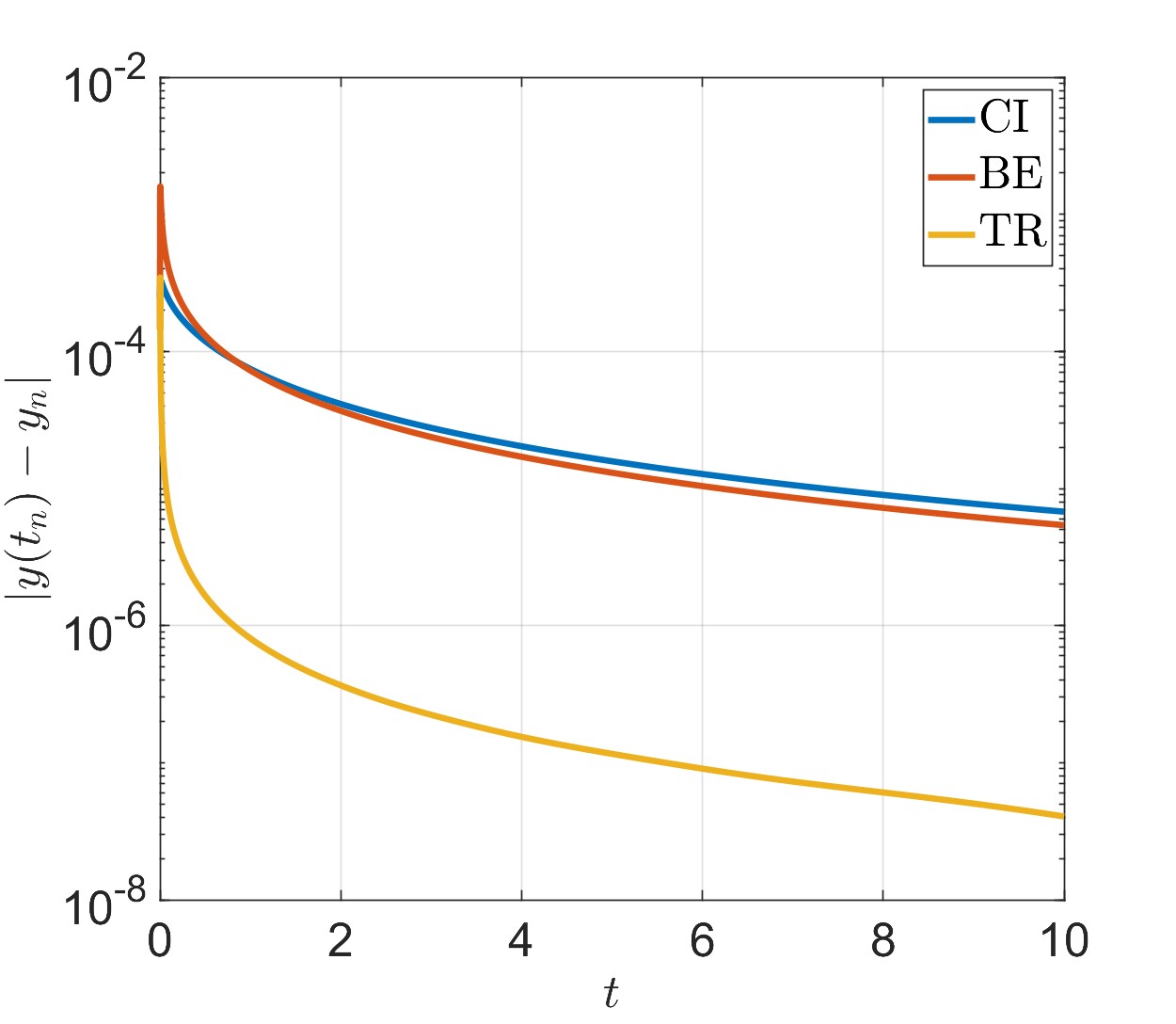}
        \caption{$\alpha = 0.5, L=128, L_{\text f} = 37$}
        \label{fig:example_2_05}
    \end{subfigure}\hfill
    \begin{subfigure}{0.33\textwidth}
        \centering
        \includegraphics[width=\textwidth]{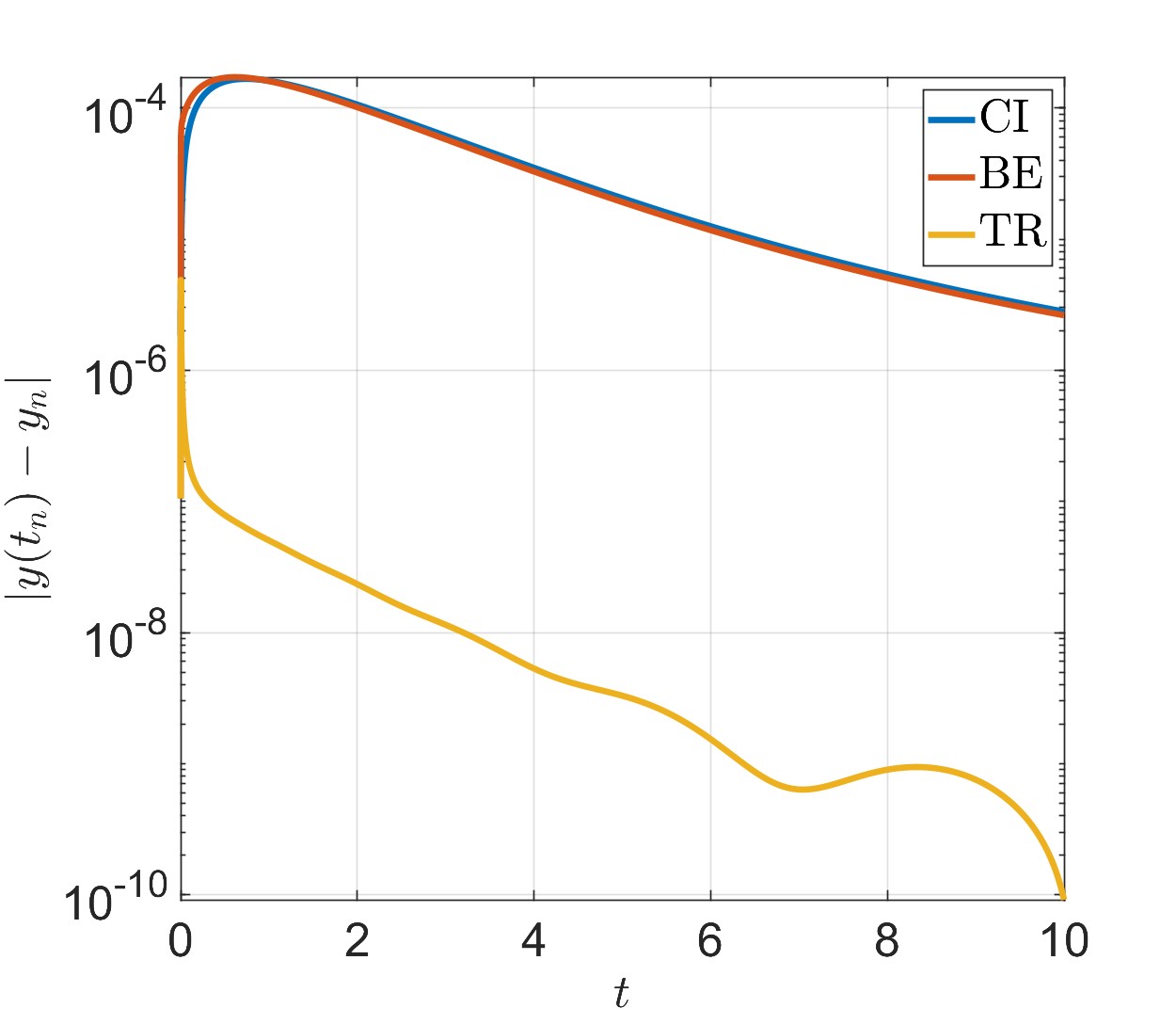}
        \caption{$\alpha = 0.9, L=512, L_{\text f} = 37$}
        \label{fig:example_2_09}
    \end{subfigure}\hfill
    \caption{Comparison of the absolute errors for (\ref{example_2}) as a function of $t$, which are obtained by the methods of constant interpolation (CI) and first order ODE by backward Euler (BE) scheme and trapezoidal rule (TR) for different values of the fractional order $\alpha$. In all cases, the time step size is set to $2^{-10}$. The $y$-axis is plotted in a logarithmic scale.}
    \label{fig:6}
\end{figure}

\begin{figure}[H]
    \centering
    \begin{subfigure}{0.33\textwidth}
        \centering
        \includegraphics[width=\textwidth]{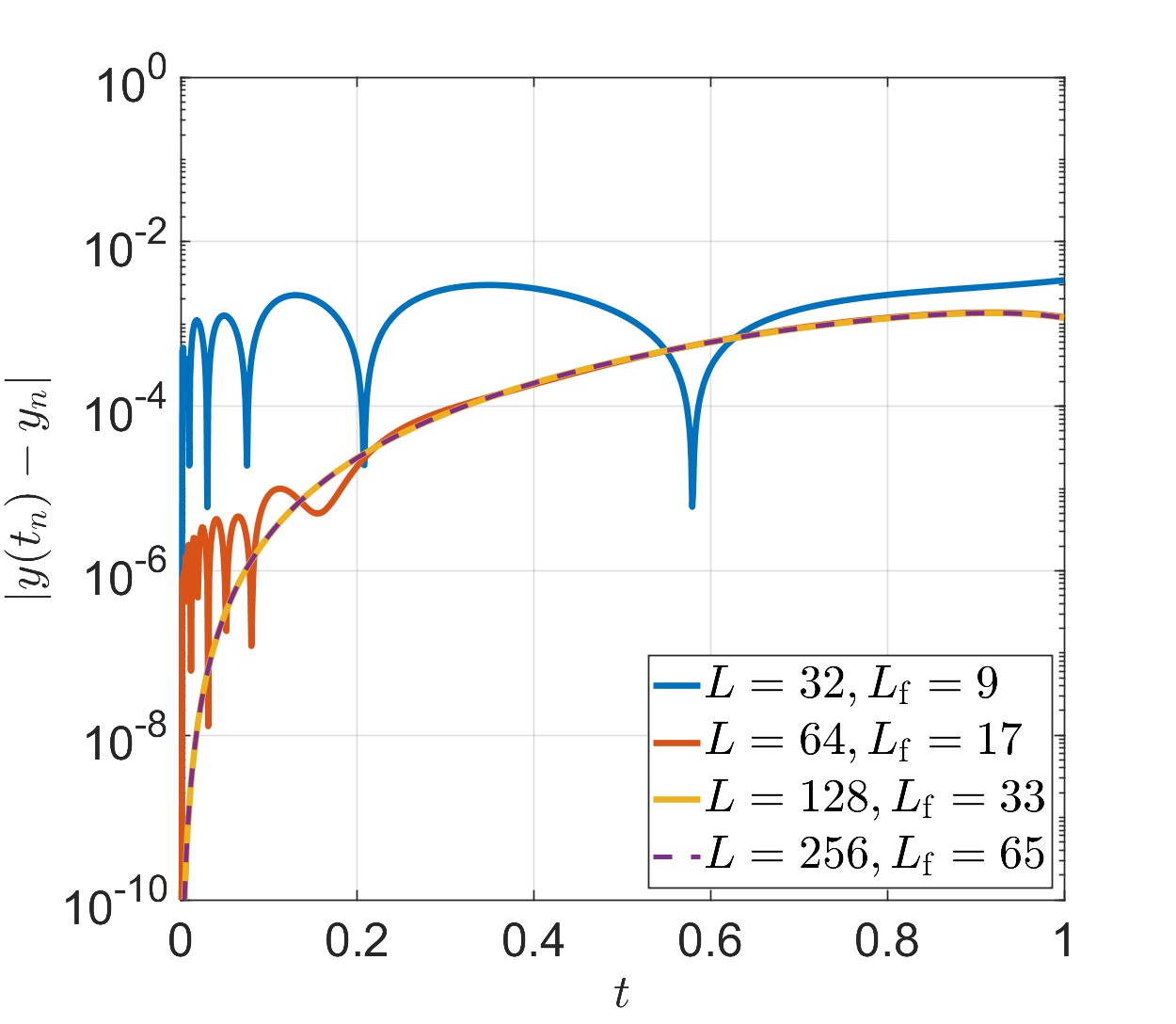}
        \caption{CI}
        \label{fig:example_1_L_ci}
    \end{subfigure}\hfill
    \begin{subfigure}{0.33\textwidth}
        \centering
        \includegraphics[width=\textwidth]{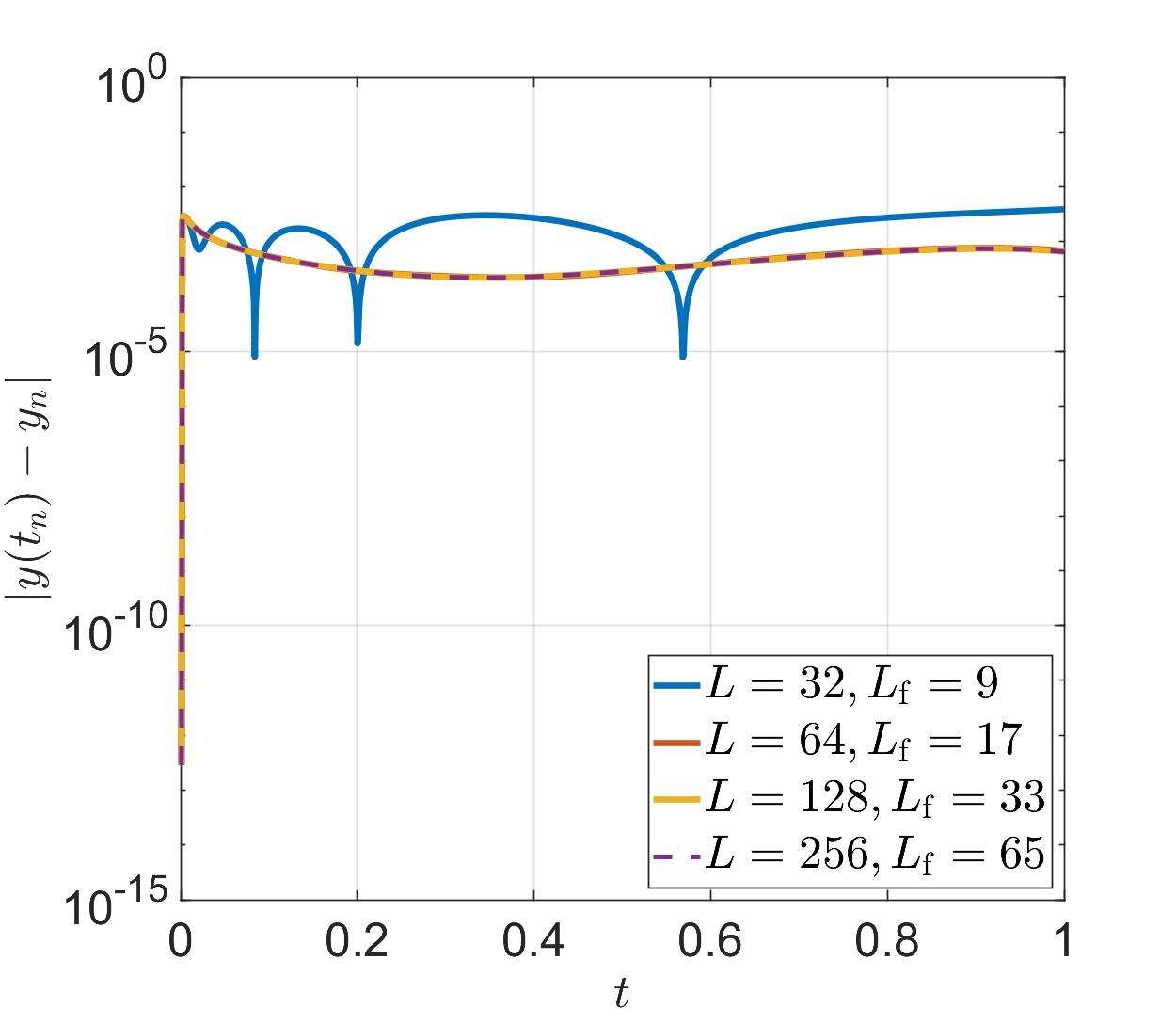}
        \caption{BE}
        \label{fig:example_1_L_be}
    \end{subfigure}\hfill
    \begin{subfigure}{0.33\textwidth}
        \centering
        \includegraphics[width=\textwidth]{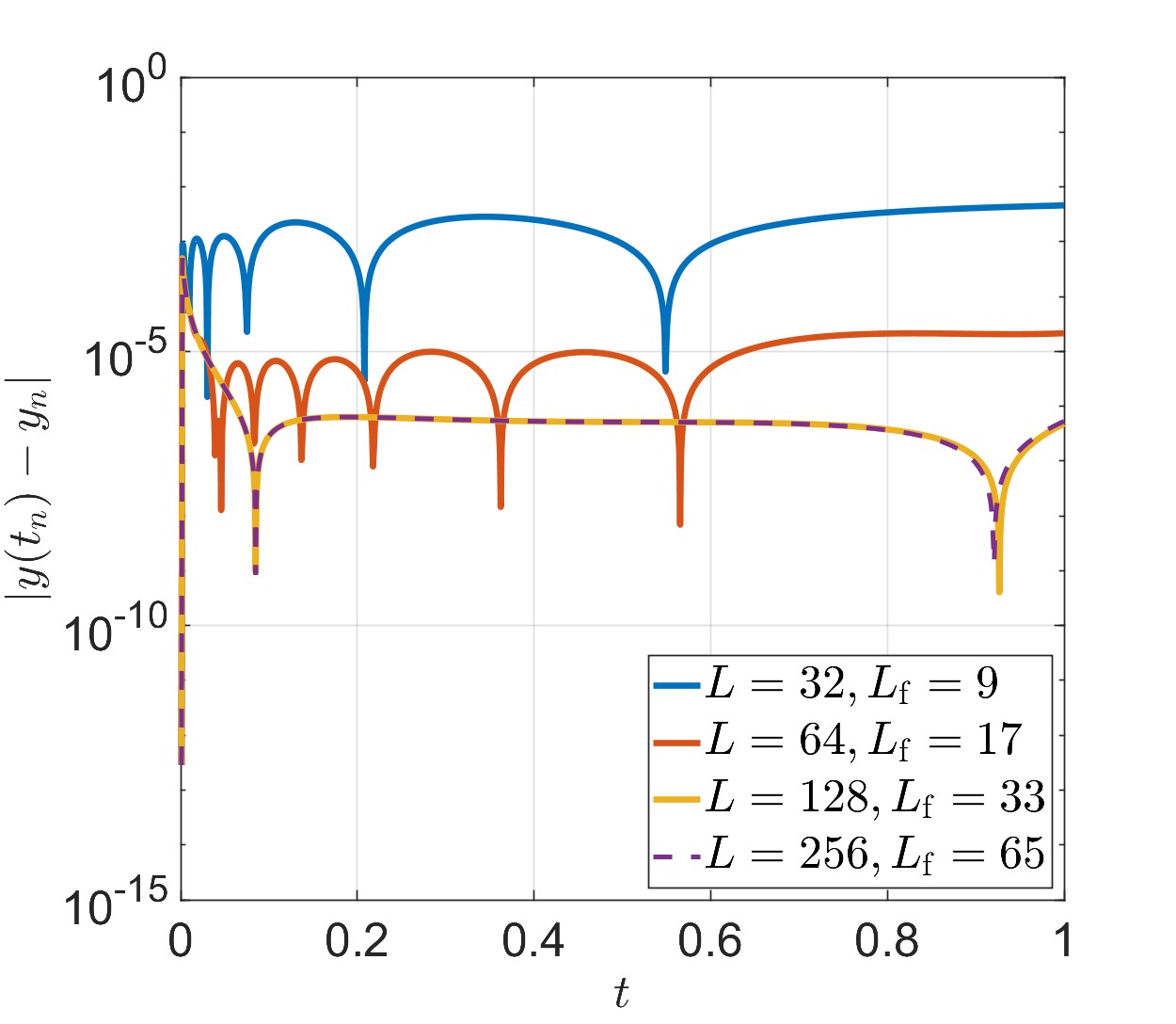}
        \caption{TR}
        \label{fig:example_1_L_trp}
    \end{subfigure}\hfill
    \caption{Comparison of the absolute errors for (\ref{example_1}) as a function of $t$, which are obtained by the methods of constant interpolation (CI) and first order ODE by backward Euler (BE) scheme and trapezoidal rule (TR) for $\alpha =0.5$ and different numbers of exponentials. In all cases, the time step size is set to $2^{-10}$. The $y$-axis is plotted in a logarithmic scale.}
    \label{fig:7}
\end{figure}

\begin{figure}[H]
    \centering
    \begin{subfigure}{0.33\textwidth}
        \centering
        \includegraphics[width=\textwidth]{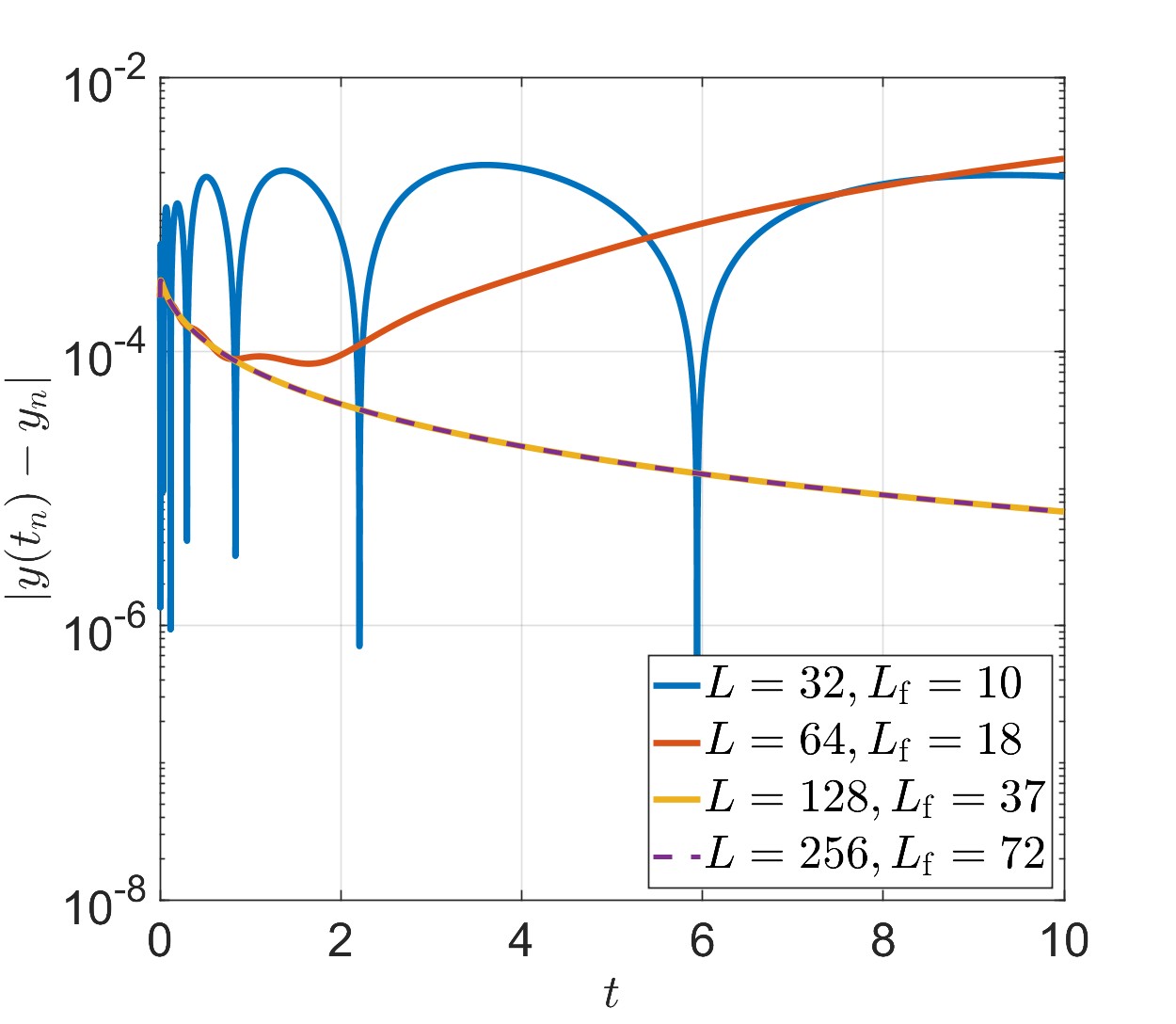}
        \caption{CI}
        \label{fig:example_2_L_ci}
    \end{subfigure}\hfill
    \begin{subfigure}{0.33\textwidth}
        \centering
        \includegraphics[width=\textwidth]{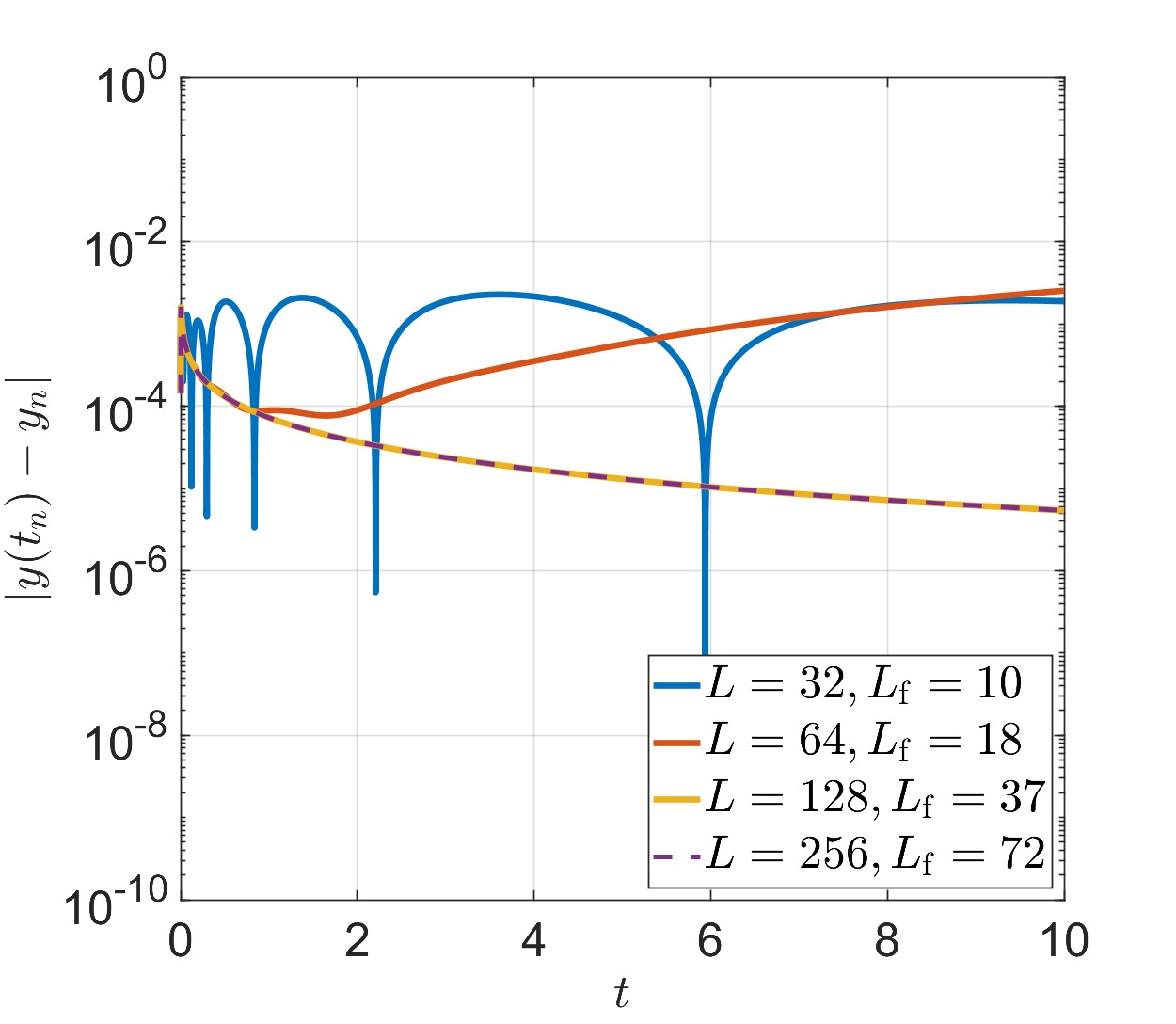}
        \caption{BE}
        \label{fig:example_2_L_be}
    \end{subfigure}\hfill
    \begin{subfigure}{0.33\textwidth}
        \centering
        \includegraphics[width=\textwidth]{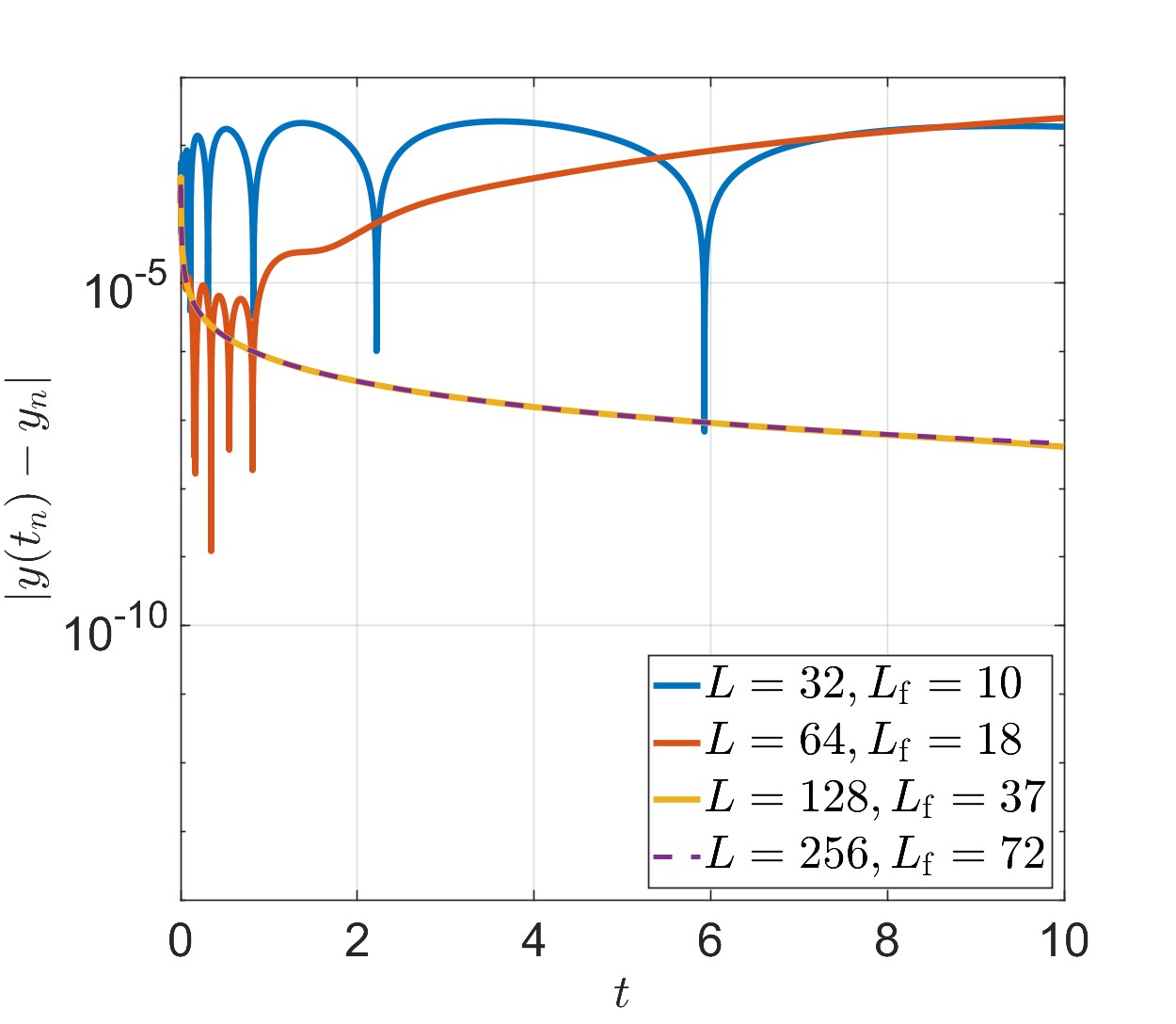}
        \caption{TR}
        \label{fig:example_2_L_trp}
    \end{subfigure}\hfill
    \caption{Comparison of the absolute errors for (\ref{example_2}) as a function of $t$, which are obtained by the methods of constant interpolation (CI) and first order ODE by backward Euler (BE) scheme and trapezoidal rule (TR) for $\alpha =0.5$ and different numbers of exponentials. In all cases, the time step size is set to $2^{-10}$. The $y$-axis is plotted in a logarithmic scale.}
    \label{fig:8}
\end{figure}

For these examples, we compare the performance of the three numerical methods (Constant Interpolation (CI), Backward Euler (BE), and Trapezoidal Rule (TR)) for solving the fractional differential equations (\ref{example_1}) and (\ref{example_2}), respectively, with fractional orders $\alpha \in \{0.1, 0.5, 0.9\}$. Tables \ref{tab:2}--\ref{tab:7} show the absolute errors at the final time $T$ for varying time step sizes $h = 2^{-2}, 2^{-3}, \ldots, 2^{-10}$, different numbers of exponential terms $L$, and corresponding $L_{\text f}$ values obtained after using Prony's method. In these tables, the column marked `EOC' states the experimentally determined order of convergence for each value of $L$ which is calculated using the error at the final time by $\log_{2}(E(h)/E(h/2))$, where $E(h) = |y(T) - y_{N}|$ is the absolute error corresponding to the step size $h$, and the symbol `***' indicates that EOC is not computed because the errors exhibit unexpected behavior, i.e., increasing instead of decreasing. The results of these tables indicate the following facts:
\begin{enumerate}
\item For all methods and all values of $\alpha$, absolute errors decrease as the time step size $h$ decreases (i.e., from $2^{-2}$ to $2^{-10}$). This is expected, as smaller time steps typically improve numerical accuracy.
\item The TR method consistently achieves the lowest errors and highest EOC compared to CI and BE methods for sufficiently large $L$. In Example \ref{ex1}, for all values of $\alpha$ and sufficiently large $L$ ($128$ for $\alpha = 0.1$ and $0.5$, and $512$ for $\alpha = 0.9$), EOC well matches the full second order convergence as we expect for the trapezoidal rule. In Example \ref{ex2}, for the TR method, the EOC is always close to $1+\alpha$. But here, increasing $L$ has not improved the EOC, see columns with $L=256$ and $1024$. Figures \ref{fig:5} and \ref{fig:6} further illustrate the high accuracy of TR compared to CI and BE methods over the entire interval for all values of $\alpha$. 
\item The BE and CI methods generally exhibit values of errors and EOC that are relatively close to each other, particularly for larger values of $\alpha$, as shown in Figures \ref{fig:5} and \ref{fig:6}. However, both CI and BE methods are limited to an EOC of approximately $1$ for sufficiently large $L$, indicating first-order convergence. This matches the well known fact the backward Euler is of first order.
\item Increasing $L$ often reduces errors and improves EOC for all methods and $\alpha$, indicating that more exponential terms enhance accuracy and convergence. Thus, the choice of $L$ is critical to achieve optimal convergence.  However, beyond a certain $L$, further increases do not significantly reduce the absolute errors of the solution anymore. This demonstrates the efficiency of our schemes which require only a modest number of exponential terms to achieve optimal convergence.  This is evident in Figures \ref{fig:7} and \ref{fig:8}, where for $\alpha = 0.5$ and all methods, no accuracy improvement occurs when increasing $L$ from $128$ to $256$. 
\end{enumerate}

\begin{table}[H]
\setlength{\tabcolsep}{4pt}
\renewcommand{\arraystretch}{1.}
\small
\centering
\caption{Errors and EOC at $T = 1$ for (\ref{example_1}) with $\alpha=0.1$ obtained by different methods.}
\label{tab:2}
\begin{tabular}{l *{4}{c} *{4}{c} *{4}{c}}
\toprule
& \multicolumn{4}{c}{\textbf{CI}} & \multicolumn{4}{c}{\textbf{BE}} & \multicolumn{4}{c}{\textbf{TR}} \\
\cmidrule(lr){2-5} \cmidrule(lr){6-9} \cmidrule(lr){10-13}
 & \multicolumn{2}{c}{$\boldsymbol{L = 64}$} & \multicolumn{2}{c}{$\boldsymbol{L = 128}$} & \multicolumn{2}{c}{$\boldsymbol{L = 64}$} & \multicolumn{2}{c}{$\boldsymbol{L = 128}$} & \multicolumn{2}{c}{$\boldsymbol{L = 64}$} & \multicolumn{2}{c}{$\boldsymbol{L = 128}$} \\
\cmidrule(lr){2-3} \cmidrule(lr){4-5} \cmidrule(lr){6-7} \cmidrule(lr){8-9} \cmidrule(lr){10-11} \cmidrule(lr){12-13}
$\boldsymbol{h}$ & \textbf{Error} & \textbf{EOC} & \textbf{Error} & \textbf{EOC} & \textbf{Error} & \textbf{EOC} & \textbf{Error} & \textbf{EOC} & \textbf{Error} & \textbf{EOC} & \textbf{Error} & \textbf{EOC} \\
\midrule
$2^{-2}$ & 8.08(-2) &          & 7.89(-2) &         & 2.30(-2) &        & 2.05(-2)  &             & 2.14(-4) &   & 1.73(-3)&  \\
$2^{-3}$ & 4.77(-2) & 0.76  & 4.57(-2) & 0.79  & 1.78(-2) & 0.37  & 1.55(-2) & 0.40    & 1.71(-4) & 0.32& 2.13(-3) &{***} \\
$2^{-4}$ & 2.80(-2) & 0.77  & 2.60(-2) & 0.81  & 1.17(-2) & 0.60  & 9.60(-3) & 0.70    & 1.10(-3)&  {***}& 8.61(-4) & 1.31\\
$2^{-5}$ & 1.65(-2) & 0.76  & 1.45(-2) & 0.84  & 7.21(-3) & 0.70  & 5.14(-3) & 0.90    & 1.70(-3)&  {***}& 2.63(-4) & 1.71 \\
$2^{-6}$ & 9.95(-3) & 0.73  & 7.98(-3) & 0.87  & 4.60(-3) & 0.65  & 2.58(-3) & 1.00    &  1.89(-3)&  {***}&7.00(-5) & 1.91 \\
$2^{-7}$ & 6.29(-3) & 0.66  & 4.31(-3) & 0.89  & 3.26(-3) & 0.50  & 1.26(-3) & 1.03    &  1.95(-3)&  {***}& 1.68(-5) & 2.06 \\
$2^{-8}$ & 4.28(-3) & 0.56 & 2.31(-3) & 0.90  & 2.60(-3) & 0.33  & 6.17(-4) & 1.03     &  1.96(-3)&  {***}& 3.55(-6) & 2.24 \\
$2^{-9}$ & 3.19(-3) & 0.42  & 1.22(-3) & 0.92  & 2.27(-3) & 0.19  & 3.02(-4) & 1.03    &  1.96(-3)&  {***}& 5.67(-7) & 2.65 \\
$2^{-10}$ & 2.61(-3) & 0.29  & 6.45(-4) & 0.93  & 2.12(-3) & 0.10  & 1.48(-4) & 1.02  &  1.96(-3)&  {***}& 4.58(-9)& 6.95\\
\bottomrule
$\boldsymbol{L_{\text f}}$ & \multicolumn{2}{c}{24} & \multicolumn{2}{c}{50} & \multicolumn{2}{c}{24} & \multicolumn{2}{c}{50} & \multicolumn{2}{c}{24} & \multicolumn{2}{c}{50} \\
\bottomrule
\end{tabular}
\end{table}

\begin{table}[H]
\setlength{\tabcolsep}{4pt}
\renewcommand{\arraystretch}{1.}
\small
\centering
\caption{Errors and EOC at $T = 1$ for (\ref{example_1}) with $\alpha=0.5$ obtained by different methods.}
\label{tab:3}
\begin{tabular}{l *{4}{c} *{4}{c} *{4}{c}}
\toprule
& \multicolumn{4}{c}{\textbf{CI}} & \multicolumn{4}{c}{\textbf{BE}} & \multicolumn{4}{c}{\textbf{TR}} \\
\cmidrule(lr){2-5} \cmidrule(lr){6-9} \cmidrule(lr){10-13}
 & \multicolumn{2}{c}{$\boldsymbol{L = 64}$} & \multicolumn{2}{c}{$\boldsymbol{L = 128}$} & \multicolumn{2}{c}{$\boldsymbol{L = 64}$} & \multicolumn{2}{c}{$\boldsymbol{L = 128}$} & \multicolumn{2}{c}{$\boldsymbol{L = 64}$} & \multicolumn{2}{c}{$\boldsymbol{L = 128}$} \\
\cmidrule(lr){2-3} \cmidrule(lr){4-5} \cmidrule(lr){6-7} \cmidrule(lr){8-9} \cmidrule(lr){10-11} \cmidrule(lr){12-13}
$\boldsymbol{h}$ & \textbf{Error} & \textbf{EOC} & \textbf{Error} & \textbf{EOC} & \textbf{Error} & \textbf{EOC} & \textbf{Error} & \textbf{EOC} & \textbf{Error} & \textbf{EOC} & \textbf{Error} & \textbf{EOC} \\
\midrule
$2^{-2}$ & 3.30(-1) &          & 3.30(-1) &         & 1.01(-1) &        & 1.01(-1)  &             & 3.34(-2) &   & 3.34(-2)&  \\
$2^{-3}$ & 1.51(-1) & 1.12  & 1.51(-1) & 1.12  & 8.48(-2) & 0.26  & 8.47(-2) & 0.26    & 1.38(-2) & 1.27& 1.38(-2) &1.27 \\
$2^{-4}$ & 7.55(-2) & 1.01  & 7.55(-2) & 1.01  & 4.86(-2) & 0.80  & 4.85(-2) & 0.80    & 3.96(-3)&  1.80& 3.98(-3) & 1.80\\
$2^{-5}$ & 3.78(-2) & 1.00  & 3.78(-2) & 1.00  & 2.42(-2) & 1.01  & 2.41(-2) & 1.01    & 9.73(-4)&  2.03& 9.95(-4) & 2.00 \\
$2^{-6}$ & 1.89(-2) & 1.00  &1.89(-2) & 1.00   & 1.16(-2) & 1.06  & 1.15(-2) & 1.06    &  2.08(-4)&  2.22&2.32(-4) & 2.10 \\
$2^{-7}$ & 9.50(-3) & 1.00  & 9.47(-3) & 1.00   & 5.56(-3) & 1.06  & 5.54(-3) & 1.06    &  3.10(-5)& 2.75& 5.22(-5) & 2.15 \\
$2^{-8}$ & 4.75(-3) & 1.00 & 4.73(-3) & 1.00   & 2.71(-3) & 1.04  & 2.68(-3) & 1.04     &  9.71(-6)& 1.68& 1.14(-5) & 2.19 \\
$2^{-9}$ & 2.38(-3) & 0.99  & 2.36(-3) & 1.00   & 1.34(-3) & 1.01  & 1.32(-3) & 1.03    &  1.99(-5)&  {***}& 2.44(-6) & 2.23 \\
$2^{-10}$ & 1.20(-3) &0.99  & 1.18(-3) & 1.00   & 6.74(-4) & 0.99  & 6.52(-4) & 1.02  &  2.13(-5)&  {***}& 4.78(-7)& 2.35\\
\bottomrule
$\boldsymbol{L_{\text f}}$ & \multicolumn{2}{c}{17} & \multicolumn{2}{c}{33} & \multicolumn{2}{c}{17} & \multicolumn{2}{c}{33} & \multicolumn{2}{c}{17} & \multicolumn{2}{c}{33} \\
\bottomrule
\end{tabular}
\end{table}

\begin{table}[H]
\setlength{\tabcolsep}{4pt}
\renewcommand{\arraystretch}{1.}
\small
\centering
\caption{Errors and EOC at $T = 1$ for (\ref{example_1}) with $\alpha=0.9$ obtained by different methods.}
\label{tab:4}
\begin{tabular}{l *{4}{c} *{4}{c} *{4}{c}}
\toprule
& \multicolumn{4}{c}{\textbf{CI}} & \multicolumn{4}{c}{\textbf{BE}} & \multicolumn{4}{c}{\textbf{TR}} \\
\cmidrule(lr){2-5} \cmidrule(lr){6-9} \cmidrule(lr){10-13}
 & \multicolumn{2}{c}{$\boldsymbol{L = 128}$} & \multicolumn{2}{c}{$\boldsymbol{L = 256}$} & \multicolumn{2}{c}{$\boldsymbol{L = 128}$} & \multicolumn{2}{c}{$\boldsymbol{L = 256}$} & \multicolumn{2}{c}{$\boldsymbol{L = 256}$} & \multicolumn{2}{c}{$\boldsymbol{L = 512}$} \\
\cmidrule(lr){2-3} \cmidrule(lr){4-5} \cmidrule(lr){6-7} \cmidrule(lr){8-9} \cmidrule(lr){10-11} \cmidrule(lr){12-13}
$\boldsymbol{h}$ & \textbf{Error} & \textbf{EOC} & \textbf{Error} & \textbf{EOC} & \textbf{Error} & \textbf{EOC} & \textbf{Error} & \textbf{EOC} & \textbf{Error} & \textbf{EOC} & \textbf{Error} & \textbf{EOC} \\
\midrule
$2^{-2}$ & 4.22(-1) &          & 4.25(-1) &         & 2.64(-1) &        & 2.66(-1)  &             & 6.11(-2) &   & 6.11(-2)&  \\
$2^{-3}$ & 1.98(-1) & 1.09  & 1.99(-1) & 1.09  & 1.94(-1) & 0.44  & 1.95(-1) & 0.44    & 1.83(-2) & 1.74& 1.83(-2) &1.74 \\
$2^{-4}$ & 9.86(-2) & 1.01  & 9.99(-2) & 1.00  & 1.03(-1) & 0.91  & 1.04(-1) & 0.91    & 4.71(-3)&  1.96& 4.70(-3) & 1.96\\
$2^{-5}$ & 4.86(-2) & 1.02  & 5.02(-2) & 0.99  & 4.94(-2) & 1.06  & 5.08(-2) & 1.04    & 1.16(-3)&  2.01& 1.16(-3) & 2.02 \\
$2^{-6}$ & 2.36(-2) & 1.04  &2.52(-2) & 0.99   & 2.30(-2) & 1.10  & 2.45(-2) & 1.05    &  2.89(-4)&  2.01&2.83(-4) & 2.03 \\
$2^{-7}$ & 1.11(-2) & 1.09  & 1.26(-2) & 1.00   & 1.05(-2) & 1.13  & 1.19(-2) & 1.04    &  7.53(-5)& 1.94& 6.94(-5) & 2.03 \\
$2^{-8}$ & 4.80(-3) & 1.22 & 6.31(-3) & 1.00   & 4.39(-3) & 1.26  & 5.89(-3) & 1.02     &  2.31(-5)& 1.70& 1.70(-5) & 2.02 \\
$2^{-9}$ & 1.64(-3) & 1.55  & 3.15(-3) & 1.00   & 1.41(-3) & 1.63  & 2.91(-3) & 1.01    &  1.02(-5)&  1.17& 4.22(-6) & 2.02 \\
$2^{-10}$ & 7.42(-5) &4.47  & 1.57(-3) & 1.00   & 4.67(-5) & 4.92  & 1.44(-3) & 1.01  &  7.06(-6)&  0.54& 1.04(-6)& 2.01\\
\bottomrule
$\boldsymbol{L_{\text f}}$ & \multicolumn{2}{c}{9} & \multicolumn{2}{c}{17} & \multicolumn{2}{c}{9} & \multicolumn{2}{c}{17} & \multicolumn{2}{c}{17} & \multicolumn{2}{c}{32} \\
\bottomrule
\end{tabular}
\end{table}

\begin{table}[H]
\setlength{\tabcolsep}{4pt}
\renewcommand{\arraystretch}{1.}
\small
\centering
\caption{Errors and EOC at $T = 10$ for (\ref{example_2}) with $\alpha=0.1$ obtained by different methods.}
\label{tab:5}
\begin{tabular}{l *{4}{c} *{4}{c} *{4}{c}}
\toprule
& \multicolumn{4}{c}{\textbf{CI}} & \multicolumn{4}{c}{\textbf{BE}} & \multicolumn{4}{c}{\textbf{TR}} \\
\cmidrule(lr){2-5} \cmidrule(lr){6-9} \cmidrule(lr){10-13}
 & \multicolumn{2}{c}{$\boldsymbol{L = 64}$} & \multicolumn{2}{c}{$\boldsymbol{L = 128}$} & \multicolumn{2}{c}{$\boldsymbol{L = 64}$} & \multicolumn{2}{c}{$\boldsymbol{L = 128}$} & \multicolumn{2}{c}{$\boldsymbol{L = 128}$} & \multicolumn{2}{c}{$\boldsymbol{L = 256}$} \\
\cmidrule(lr){2-3} \cmidrule(lr){4-5} \cmidrule(lr){6-7} \cmidrule(lr){8-9} \cmidrule(lr){10-11} \cmidrule(lr){12-13}
$\boldsymbol{h}$ & \textbf{Error} & \textbf{EOC} & \textbf{Error} & \textbf{EOC} & \textbf{Error} & \textbf{EOC} & \textbf{Error} & \textbf{EOC} & \textbf{Error} & \textbf{EOC} & \textbf{Error} & \textbf{EOC} \\
\midrule
$2^{-2}$ & 2.58(-4) &          & 1.50(-4) &         & 4.74(-4) &        & 3.59(-4)  &             & 2.35(-4) &   & 2.35(-4)&  \\
$2^{-3}$ & 1.90(-4) & 0.44  & 8.20(-5) & 0.87  & 2.91(-4) & 0.70  & 1.79(-4) & 1.01    & 1.08(-4) & 1.12& 1.08(-4) &1.12 \\
$2^{-4}$ & 1.53(-4) & 0.32  & 4.41(-5) & 0.89  & 1.99(-4) & 0.54  & 8.89(-5) & 1.01    & 5.04(-5)&  1.11& 5.04(-5) & 1.11\\
$2^{-5}$ & 1.32(-4) & 0.21  & 2.35(-5) & 0.91  & 1.54(-4) & 0.37  & 4.41(-5) & 1.01    & 2.34(-5)& 1.10& 2.34(-5) & 1.10 \\
$2^{-6}$ & 1.21(-4) & 0.12  & 1.24(-5) & 0.92  & 1.31(-4) & 0.23  & 2.19(-5) & 1.01    &  1.09(-5)&  1.10&1.09(-5) & 1.10 \\
$2^{-7}$ & 1.15(-4) & 0.07  & 6.55(-6) & 0.93  & 1.20(-4) & 0.13  & 1.08(-5) & 1.01    &  5.08(-6)&  1.10& 5.08(-6) & 1.10 \\
$2^{-8}$ & 1.12(-4) & 0.04 & 3.42(-6) & 0.94  & 1.14(-4) & 0.07  & 5.39(-6) & 1.01     &  2.37(-6)&  1.10& 2.37(-6) & 1.10 \\
$2^{-9}$ & 1.11(-4) & 0.02  & 1.78(-6) & 0.94  & 1.12(-4) & 0.04  & 2.68(-6) & 1.01    &  1.10(-6)&  1.10& 1.10(-6) & 1.10 \\
$2^{-10}$ & 1.10(-4) & 0.01  & 9.24(-7) & 0.95  & 1.10(-4) & 0.02  & 1.33(-6) & 1.01  &  5.15(-7)&  1.10& 5.15(-7)& 1.10\\
\bottomrule
$\boldsymbol{L_{\text f}}$ & \multicolumn{2}{c}{27} & \multicolumn{2}{c}{54} & \multicolumn{2}{c}{27} & \multicolumn{2}{c}{54} & \multicolumn{2}{c}{54} & \multicolumn{2}{c}{105} \\
\bottomrule
\end{tabular}
\end{table}

\begin{table}[H]
\setlength{\tabcolsep}{4pt}
\renewcommand{\arraystretch}{1.}
\small
\centering
\caption{Errors and EOC at $T = 10$ for (\ref{example_2}) with $\alpha=0.5$ obtained by different methods.}
\label{tab:6}
\begin{tabular}{l *{4}{c} *{4}{c} *{4}{c}}
\toprule
& \multicolumn{4}{c}{\textbf{CI}} & \multicolumn{4}{c}{\textbf{BE}} & \multicolumn{4}{c}{\textbf{TR}} \\
\cmidrule(lr){2-5} \cmidrule(lr){6-9} \cmidrule(lr){10-13}
 & \multicolumn{2}{c}{$\boldsymbol{L = 64}$} & \multicolumn{2}{c}{$\boldsymbol{L = 128}$} & \multicolumn{2}{c}{$\boldsymbol{L = 64}$} & \multicolumn{2}{c}{$\boldsymbol{L = 128}$} & \multicolumn{2}{c}{$\boldsymbol{L = 128}$} & \multicolumn{2}{c}{$\boldsymbol{L = 256}$} \\
\cmidrule(lr){2-3} \cmidrule(lr){4-5} \cmidrule(lr){6-7} \cmidrule(lr){8-9} \cmidrule(lr){10-11} \cmidrule(lr){12-13}
$\boldsymbol{h}$ & \textbf{Error} & \textbf{EOC} & \textbf{Error} & \textbf{EOC} & \textbf{Error} & \textbf{EOC} & \textbf{Error} & \textbf{EOC} & \textbf{Error} & \textbf{EOC} & \textbf{Error} & \textbf{EOC} \\
\midrule
$2^{-2}$ & 4.14(-3) &          & 1.69(-3) &         & 4.33(-3) &        & 1.77(-3)  &             & 2.01(-4) &   & 2.01(-4)&  \\
$2^{-3}$ & 3.34(-3) & 0.31  & 8.52(-4) & 0.99  & 3.37(-3) & 0.36  & 8.19(-4) & 1.11    & 6.91(-5) & 1.54& 6.91(-5) &1.54 \\
$2^{-4}$ & 2.94(-3) & 0.18  & 4.28(-4) & 0.99  & 2.93(-3) & 0.20  & 3.87(-4) & 1.08    & 2.39(-5)&  1.53& 2.39(-5) & 1.53\\
$2^{-5}$ & 2.74(-3) & 0.10  & 2.14(-4) & 0.99  & 2.73(-3) & 0.10  & 1.86(-4) & 1.06   & 8.36(-6)& 1.52& 8.36(-6) & 1.52\\
$2^{-6}$ & 2.64(-3) & 0.05  &1.07(-4) & 1.00   & 2.63(-3) & 0.05  & 9.07(-5) & 1.04    &  2.92(-6)&  1.51&2.93(-6) & 1.51 \\
$2^{-7}$ & 2.58(-3) & 0.03  & 5.39(-5) & 1.00   & 2.58(-3) & 0.03  & 4.45(-5) & 1.03    &  1.02(-6)& 1.51& 1.03(-6) & 1.51 \\
$2^{-8}$ & 2.56(-3) & 0.01 & 2.70(-5) & 1.00   & 2.56(-3) & 0.01  & 2.19(-5) & 1.02     &  3.58(-7)& 1.52& 3.63(-7) & 1.51 \\
$2^{-9}$ & 2.55(-3) & 0.01  & 1.35(-5) & 1.00   & 2.55(-3) & 0.01  & 1.08(-5) & 1.01    &  1.23(-7)&  1.54& 1.28(-7) & 1.50 \\
$2^{-10}$ & 2.54(-3) &0.01  & 6.77(-6) & 1.00   & 2.54(-3) & 0.01  & 5.39(-6) & 1.01  &  4.05(-8)&  1.61& 4.51(-8)& 1.50\\
\bottomrule
$\boldsymbol{L_{\text f}}$ & \multicolumn{2}{c}{18} & \multicolumn{2}{c}{37} & \multicolumn{2}{c}{18} & \multicolumn{2}{c}{37} & \multicolumn{2}{c}{37} & \multicolumn{2}{c}{72} \\
\bottomrule
\end{tabular}
\end{table}

\begin{table}[H]
\setlength{\tabcolsep}{4pt}
\renewcommand{\arraystretch}{1.}
\small
\centering
\caption{Errors and EOC at $T = 10$ for (\ref{example_2}) with $\alpha=0.9$ obtained by different methods.}
\label{tab:7}
\begin{tabular}{l *{4}{c} *{4}{c} *{4}{c}}
\toprule
& \multicolumn{4}{c}{\textbf{CI}} & \multicolumn{4}{c}{\textbf{BE}} & \multicolumn{4}{c}{\textbf{TR}} \\
\cmidrule(lr){2-5} \cmidrule(lr){6-9} \cmidrule(lr){10-13}
 & \multicolumn{2}{c}{$\boldsymbol{L = 256}$} & \multicolumn{2}{c}{$\boldsymbol{L = 512}$} & \multicolumn{2}{c}{$\boldsymbol{L = 256}$} & \multicolumn{2}{c}{$\boldsymbol{L = 512}$} & \multicolumn{2}{c}{$\boldsymbol{L = 512}$} & \multicolumn{2}{c}{$\boldsymbol{L = 1024}$} \\
\cmidrule(lr){2-3} \cmidrule(lr){4-5} \cmidrule(lr){6-7} \cmidrule(lr){8-9} \cmidrule(lr){10-11} \cmidrule(lr){12-13}
$\boldsymbol{h}$ & \textbf{Error} & \textbf{EOC} & \textbf{Error} & \textbf{EOC} & \textbf{Error} & \textbf{EOC} & \textbf{Error} & \textbf{EOC} & \textbf{Error} & \textbf{EOC} & \textbf{Error} & \textbf{EOC} \\
\midrule
$2^{-2}$ & 7.75(-4) &          & 7.70(-4) &         & 8.47(-4) &        & 8.44(-4)  &             &1.47(-5) &   & 1.47(-5)\phantom{0}&  \\
$2^{-3}$ & 3.76(-4) & 1.04  & 3.71(-4) & 1.05  & 3.80(-4) & 1.15  & 3.77(-4) & 1.16    &3.78(-6) & 1.97& 3.78(-6)\phantom{0} &1.97 \\
$2^{-4}$ & 1.87(-4) & 1.01  & 1.82(-4) & 1.03  & 1.82(-4) & 1.06  & 1.78(-4) & 1.08    & 9.68(-7)&  1.97& 9.68(-7)\phantom{0} & 1.97\\
$2^{-5}$ & 9.56(-5) & 0.97  & 9.04(-5) & 1.01  & 9.14(-5) & 1.00  & 8.65(-5) & 1.04    &2.47(-7)&  1.97& 2.47(-7)\phantom{0} & 1.97 \\
$2^{-6}$ & 5.02(-5) & 0.93  &4.50(-5) & 1.01   & 4.77(-5) & 0.94  & 4.26(-5) & 1.02    &  6.32(-8)&  1.97&6.34(-8)\phantom{0} & 1.97 \\
$2^{-7}$ & 2.77(-5) & 0.86  & 2.24(-5) & 1.00   & 2.63(-5) & 0.86  & 2.11(-5) & 1.01    & 1.60(-8)& 1.98& 1.62(-8)\phantom{0} & 1.97 \\
$2^{-8}$ & 1.65(-5) & 0.75 & 1.12(-5) & 1.00   & 1.57(-5) & 0.74 & 1.05(-5) & 1.01     &  3.98(-9)& 2.01& 4.16(-9)\phantom{0} & 1.96 \\
$2^{-9}$ & 1.09(-5) & 0.60  & 5.60(-6) & 1.00   & 1.05(-5) & 0.58  & 5.25(-6) & 1.00    &8.85(-10)& 2.17& 1.06(-9)\phantom{0} & 1.96 \\
$2^{-10}$ & 8.09(-6) &0.43  & 2.80(-6) & 1.00   & 7.91(-6) & 0.41  & 2.62(-6) & 1.00  & 9.16(-11)& 3.27& 2.74(-10)& 1.96\\
\bottomrule
$\boldsymbol{L_{\text f}}$ & \multicolumn{2}{c}{19} & \multicolumn{2}{c}{37} & \multicolumn{2}{c}{19} & \multicolumn{2}{c}{37} & \multicolumn{2}{c}{37} & \multicolumn{2}{c}{70} \\
\bottomrule
\end{tabular}
\end{table}

\begin{figure}[H]
    \centering
    \begin{subfigure}{0.33\textwidth}
        \centering
        \includegraphics[width=\textwidth]{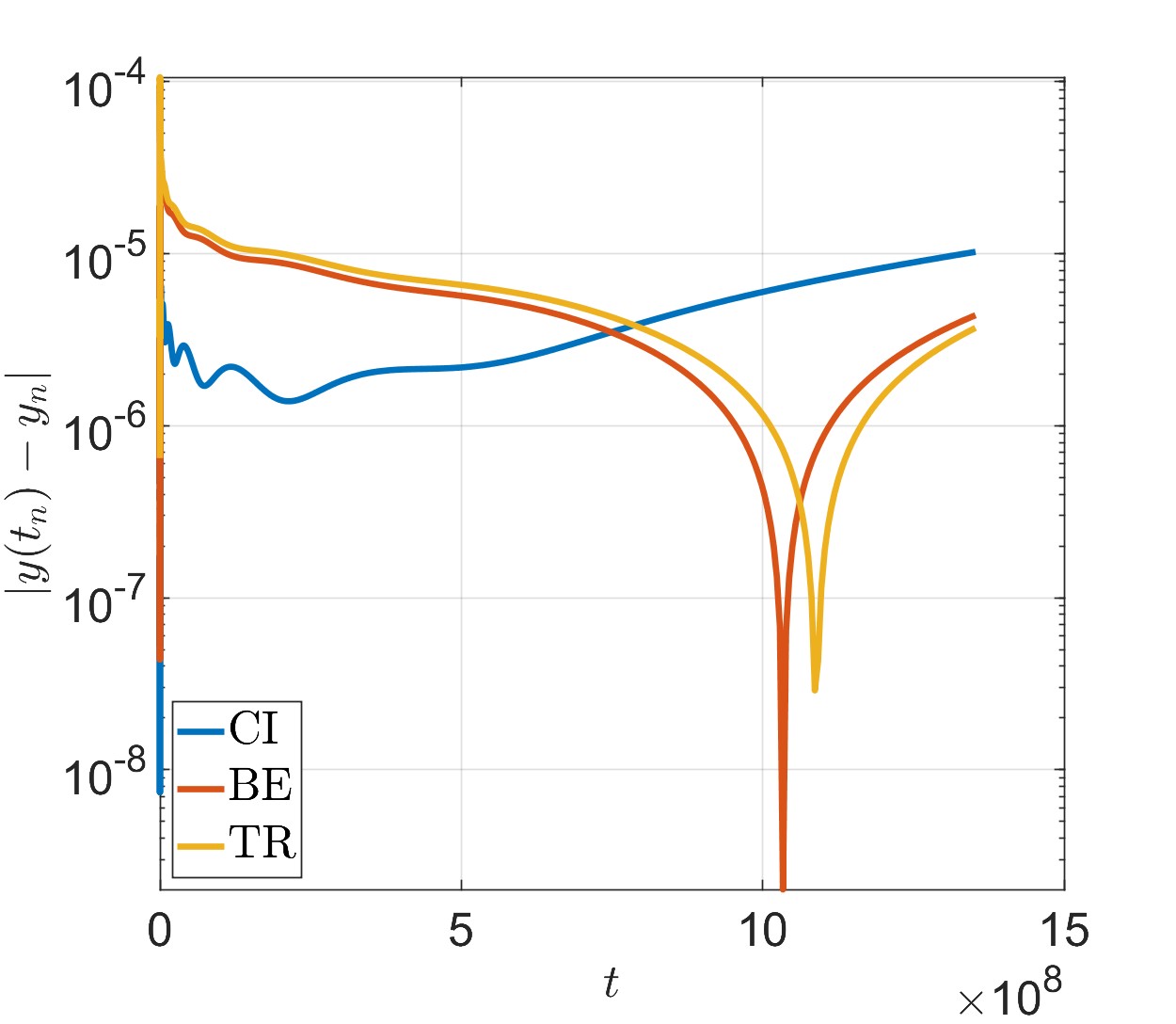}
        \caption{ $L =64, L_{\text f} = 34$}
        \label{fig:err_kelvin_ci_br_tr}
    \end{subfigure}\hfill
    \begin{subfigure}{0.33\textwidth}
        \centering
        \includegraphics[width=\textwidth]{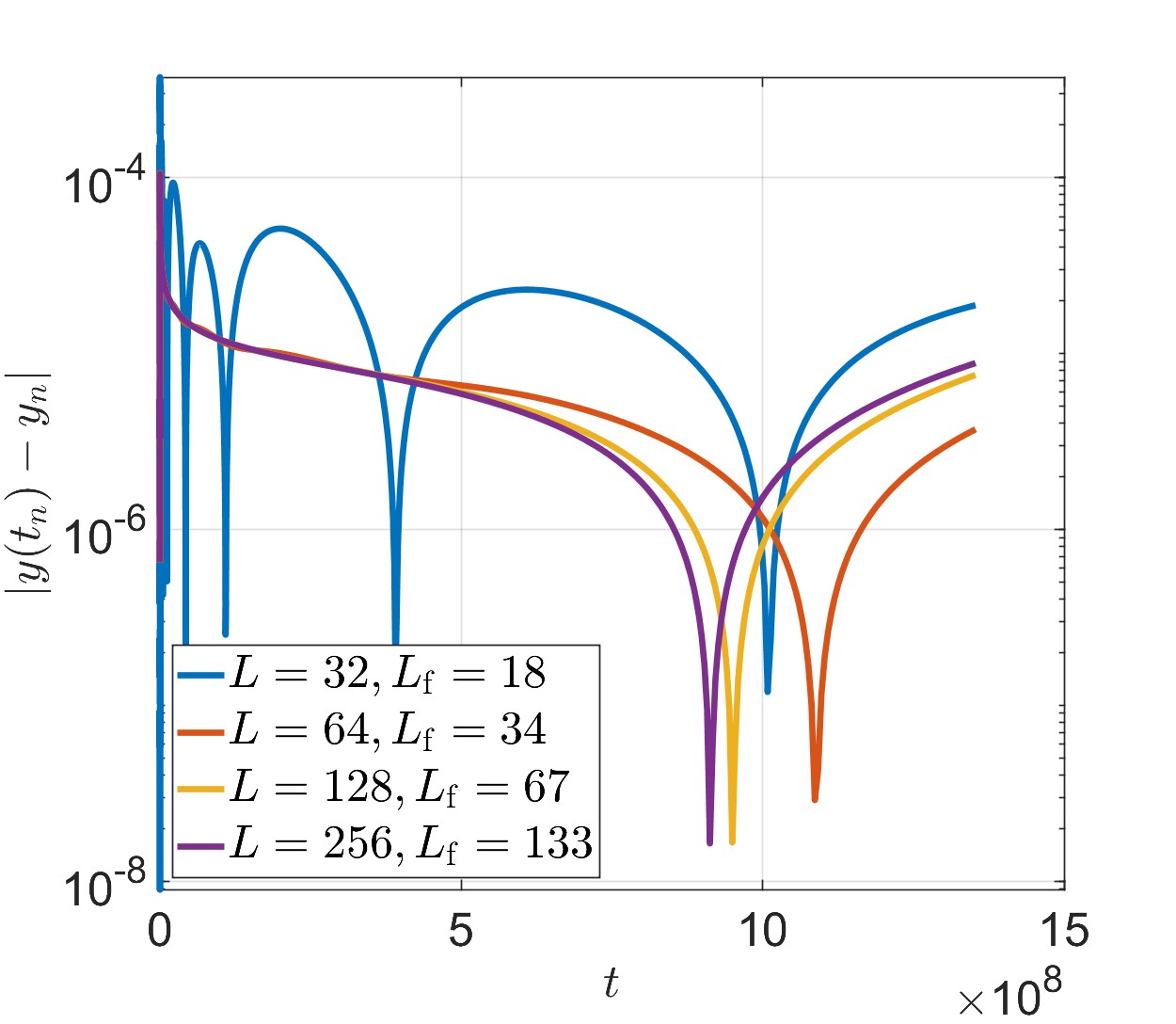}
        \caption{TR}
        \label{fig:err_kelvin_tr_L}
    \end{subfigure}\hfill
    \begin{subfigure}{0.33\textwidth}
        \centering
        \includegraphics[width=\textwidth]{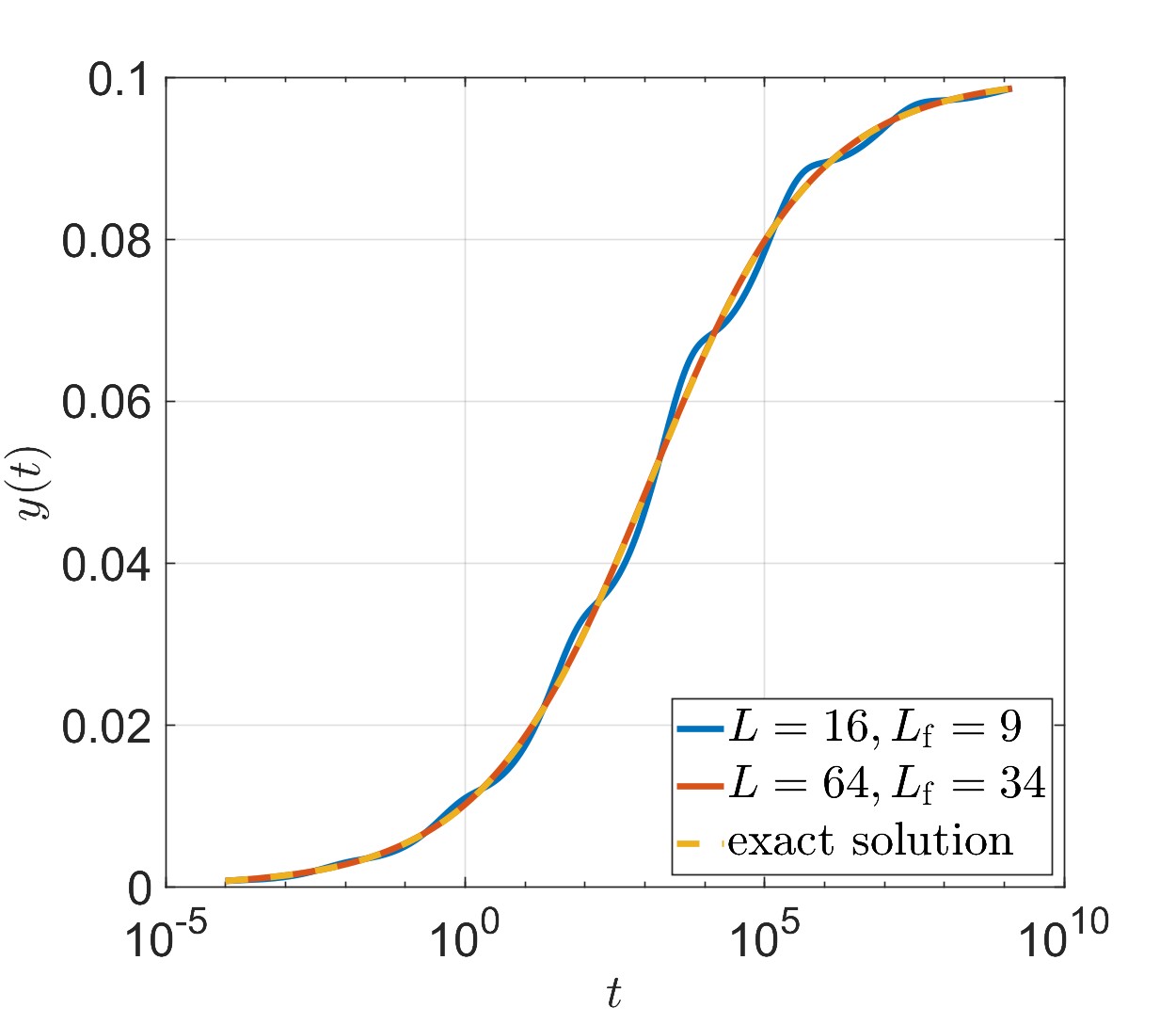}
        \caption{}
        \label{fig:err_kelvin_tr_exact}
    \end{subfigure}\hfill
    \caption{Accuracy test for (\ref{example_3}): (a) Comparison of the absolute errors as a function of $t$, which are obtained by the method of constant interpolation and first order ODE by backward Euler scheme and trapezoidal rule for $L =64$. (b) Comparison of the absolute errors as a function of $t$, which are obtained by the trapezoidal rule for various values of $L$. (c) Comparison of the numerical solution obtained by the trapezoidal rule for $L = 16$ and $64$ and the exact solution given by (\ref{81}). In (a) and (b), the $y$-axis is shown on a logarithmic scale, while in (c), the $x$-axis is displayed on a logarithmic scale.}
    \label{fig:9}
\end{figure}

\begin{example}
	\label{ex3}
Here we consider the differential equation of motion of a massless one-dimensional fractional Kelvin-Voigt model subject to the external force $f(t)=1$ for $t\neq 0$,
\begin{equation}\label{example_3}
	D^{\alpha}y(t) + ky(t) = f(t), \qquad y(0)=0,
\end{equation}
with $c=100, k=10$ and $\alpha=0.3$, see \cite{schmidt2006critique}. We use the same discretization of the variable $t$ used in \cite{schmidt2006critique} by the mesh points
\begin{equation}\label{76}
t_{0}=0, \qquad t_{j}=t_{j-1}+h_{j} \qquad j=1,2,\ldots,5000
\end{equation}
with $h_{1}=10^{-4}$ and $h_{j}=1.005h_{j-1}$ for $j=2,\dots,5000$. The exact solution of this equation is
\begin{equation}\label{81}
	y(t) = \dfrac{1}{k} \left( 1 - E_{\alpha}(-\dfrac{k}{c}t^{\alpha}) \right), \quad t>0,
\end{equation}
where $E_{\alpha}$ is the Mittag-Leffler function of order $\alpha$. The numerical methods CI, BE, and TR are applied to solve problem (\ref{example_3}). A comparison of absolute errors over time, as shown in Figure \ref{fig:err_kelvin_ci_br_tr}, reveals that the BE and TR methods exhibit similar behavior but do not significantly improve compared to the CI method. These findings align with the results from previous examples, which suggest that the TR method can offer considerable improvements over the CI and BE methods mainly for larger $\alpha$ values. Nevertheless, the trapezoidal rule (TR) can be considered as a more efficient option compared to the other two methods for solving fractional differential equations. To verify the impact of the parameter $L$ on solution accuracy, Figure \ref{fig:err_kelvin_tr_L} illustrates that for the trapezoidal rule, increasing $L$ enhances accuracy, but beyond $L \geq 64$, further increases yield no significant improvement. Figure \ref{fig:err_kelvin_tr_exact} compares the numerical solutions using the trapezoidal rule for $L = 16$ and $L = 64$ with the exact solution from (\ref{81}), demonstrating that a considerable error reduction requires a sufficiently large $L$.
\end{example}

\section{Conclusions}

In this work, we proposed a sum of exponential approximation for power functions $t^{\alpha-1}, 0 < \alpha < 1$ with pre-specified number of exponential terms. We then combined this with Prony's algorithm to efficiently decrease the number of exponentials while preserving or even enhancing the maximum approximation error. The findings highlight the importance of choosing the appropriate method for solving fractional differential equations and number of exponential terms $L$ based on the fractional order $\alpha$ to achieve optimal accuracy and convergence. Here, the Trapezoidal Rule (TR) is the most effective method for solving fractional differential equations, offering the lowest errors and highest EOC. The CI and BE methods exhibit slower convergence and lower accuracy comparatively, rendering them suboptimal for high‑precision applications. In all methods, increasing the number of exponential terms $L$ enhances the solution accuracy.

\section*{Appendix}

Table \ref{tab:8} presents the same parameters as Table \ref{tab:1}, but evaluated over the interval $[10^{-2}, 10^{3}]$. It includes the initial maximum error, $\max_{t \in [10^{-2},10^{3}]} |e(t,\theta_0)|$, and the maximum error after applying Prony's method, $\max_{t \in [10^{-2},10^{3}]} |e_{\text p}(t)|$, both of which are computed using rescaling of the weights and exponents according to Remark \ref{rem:rescaling} for the interval $[10^{-5},1]$. Table \ref{tab:8} confirms that Prony's method, with the rescaling approach, effectively approximates the kernel function over the large interval $[10^{-2}, 10^3]$ without loss in accuracy and produces a significant term reduction ($L_{\text f} \ll L$) across all cases. As in Table \ref{tab:1}, the parameter $K$ ranges from $1$ to $5$, and in all cases, the obtained $L_{\text p}$ is equal to the initialized value $M $.  Compared to Table \ref{tab:1}, the value of $L_{\text f}$ increases for fixed $\alpha$ and $L$ because reducing $\delta$ from $10^{-2}$ to $10^{-5}$ widens the interval, resulting in a greater number of positive nodes. 

Therefore,  Table \ref{tab:8} demonstrates the robustness of Prony's method for larger intervals, making it suitable for practical applications where computational efficiency and accuracy over wide ranges are critical.
\begin{table}[H]
\centering
\renewcommand{\arraystretch}{1.2} 
\caption{Max. error values obtained before and after the application of Prony's method for various values of $L$ and $\alpha \in \{0.1, 0.5, 0.9\}$ over the interval $[10^{-2}, 10^{3}]$. In all cases, the numbers of non-positive nodes $M$, the values of $L_{\text p}$ and $K$ in Prony's method and the final number of exponential terms after reduction $L_{\text f}$ are provided. Here, the value of $\epsilon$ is set to $10^{-10}$.}
\label{tab:8}
\begin{tabular}{ccccccccc}
\toprule
$\boldsymbol{\alpha}$ & $\boldsymbol{L}$ & $\boldsymbol{M}$ & $\boldsymbol{L_{\text p}}$ & $\boldsymbol{K}$ & $\boldsymbol{L_{\text f}}$ & $\boldsymbol{\max_{t \in [10^{-2},10^3]} |e(t,\theta_{0})|}$ & $\boldsymbol{\max_{t \in [10^{-2}, 10^3]} |e_{\text p}(t)|}$ \\
\midrule
\multirow{4}{*}{$0.1$} & $32$  & 20 &20 &  1 & 13  & 2.472386(-1)\phantom{0} & 2.472386(-1)\phantom{0} \\
& $64$  & 41  &  41&1 & 24  & 1.335328(-4)\phantom{0} & 1.335343(-4)\phantom{0} \\
& $128$  & 81  & 81  &  3 & 50  & 1.582185(-8)\phantom{0} & 1.542839(-10) \\
& $256$  & 163 & 163  &3& 96  & 8.214528(-9)\phantom{0} & 6.045368(-10) \\
\midrule
\multirow{4}{*}{$0.5$} & $32$  &24  & 24 &1 & 9 & 2.050205(-1)\phantom{0} & 2.050205(-1)\phantom{0} \\
& $64$  & 49& 49   & 2 & 17  & 1.159256(-3)\phantom{0} & 1.159256(-3)\phantom{0} \\
& $128$  & 98 & 98  &  3 & 33  & 4.535111(-8)\phantom{0} & 4.532144(-8)\phantom{0} \\
& $256$  & 195 & 195  & 4 & 65 & 3.326726(-10) & 4.989634(-12)\\
\midrule
\multirow{4}{*}{$0.9$} & $128$  & 121  & 121& 2 & 9  & 4.834038(-3)\phantom{0} & 4.834038(-3)\phantom{0} \\
& $256$  & 242  & 242   & 3 & 17   & 3.344168(-5)\phantom{0} & 3.344168(-5)\phantom{0} \\
& $512$  & 484 & 484  & 4 & 32 & 2.046576(-9)\phantom{0} & 2.046172(-9)\phantom{0} \\
& $1024$  &968 & 968  & 5 & 61  & 8.260503(-12) & 5.237699(-12) \\
\bottomrule
\end{tabular}
\end{table}

In Table \ref{tab:9}, we focus on fixed $L$ values ($L = 256$ for $\alpha = 0.1, 0.5$, and $L = 1024$ for $\alpha = 0.9$) and varying the parameter $\epsilon$ from $10^{-10}$ to $10^{-13}$. The table shows that the final number of terms ($L_{\text f}$) decreases as $\epsilon$ becomes stricter (smaller). For $\alpha = 0.1$, $L_{\text f}$ decreases from $96$ ($\epsilon = 10^{-10}$) to $83$ ($\epsilon = 10^{-13}$). For $\alpha = 0.5$, $L_{\text f}$ decreases from $65$ to $55$, and for $\alpha = 0.9$, $L_{\text f}$ decreases from $61$ to $51$. This reduction in $L_{\text f}$ demonstrates that Prony's method effectively compresses the number of terms while maintaining accuracy, especially for stricter $\epsilon$ values. Furthermore, decreasing $\epsilon$ improves both initial and Prony errors, with Prony errors remaining close to or better than initial errors in most cases.

In Table \ref{tab:10}, we have shown the values of the parameters $\rho_k$ (weights) and $\eta_k$ (exponents) after the application of Prony's method in the cases of $L = 256$ for $\alpha = 0.1, 0.5$ and $L = 1024$ for $\alpha = 0.9$ over three intervals, $[10^{-2}, 1]$, $[10^{-5}, 1]$, and $[10^{-2}, 10^3]$. In comparison with what we explained about Figure~\ref{fig:3}, before using Prony's method, there are very small (very close to zero) weights $w_{l}$ and exponents $b_{l}$ obtained by the nodes $\omega_{i}$ on the interval $[l_{\min},0]$. However, with Prony's method we get weights $\rho_{k}$ and exponents $\eta_{k}$ which are orders of magnitude larger. The weights $\rho_k$ are generally larger for smaller intervals $[10^{-2}, 1]$, $[10^{-5}, 1]$ and decrease significantly for the larger interval $[10^{-2}, 10^3]$ due to the rescaling described in Remark \ref{rem:rescaling} ($w_l = \tilde{w}_l / T^{1-\alpha}$). As $\alpha$ increases, the weights are significantly larger (especially for smaller intervals). The exponents $\eta_k$ for the large interval $[10^{-2}, 10^3]$ are orders of magnitude smaller (closer to zero) than those for $[10^{-2}, 1]$ or $[10^{-5}, 1]$ due to the rescaling factor $b_l = \tilde{b}_l / T$. 
\begin{table}[H]
\centering
\renewcommand{\arraystretch}{1.2} 
\caption{Max. error values obtained before and after the application of Prony's method for fixed values of $L$ and $\alpha \in \{0.1, 0.5, 0.9\}$ over the interval $[10^{-2}, 10^{3}]$. In all cases, the numbers of non-popsitive nodes, $M$, the values of $L_{\text p}$ and $K$ in Prony's method, and the final number of exponential terms after reduction $L_{\text f}$ are provided.}
\label{tab:9}
\begin{tabular}{ccccccccc}
\toprule
$\boldsymbol{\alpha}$ & $\boldsymbol{\epsilon}$ & $\boldsymbol{M}$ & $\boldsymbol{L_{\text p}}$ & $\boldsymbol{K}$ & $\boldsymbol{L_{\text f}}$ & $\boldsymbol{\max_{t \in [10^{-2},10^{3}]} |e(t,\theta_{0})|}$ & $\boldsymbol{\max_{t \in [10^{-2},10^{3}]} |e_{\text p}(t)|}$ \\
\midrule
\multirow{4}{*}{ $0.1$  $(L=256)$} & $10^{-10}$ & 163 & 163 & 3 &96 & 8.214528(-9)\phantom{0} & 6.045368(-10) \\
& $10^{-11}$ & 168& 168 & 3 & 91 & 9.331558(-10) & 5.064520(-10) \\
& $10^{-12}$ & 173 &173 & 4 & 87 & 1.058709(-10) & 1.094236(-12) \\
& $10^{-13}$ & 177 & 177& 4 & 83 & 1.196554(-11) & 6.685971(-13) \\
\midrule
\multirow{4}{*}{ $0.5$ $(L=256)$} & $10^{-10}$ & 195 & 195 & 4 & 65 & 3.326726(-10) & 4.989634(-12) \\
& $10^{-11}$ & 199 & 199 & 4 & 61 & 3.741008(-11) & 3.887078(-12) \\
& $10^{-12}$ &203 & 203 & 5 & 58 & 4.185097(-12) & 3.197442(-14) \\
& $10^{-13}$ & 206 & 206 & 5 & 55 & 5.879741(-13) & 1.278977(-13) \\
\midrule
\multirow{4}{*}{$0.9$ $(L=1024)$} & $10^{-10}$ & 968 & 968 & 5 & 61 & 8.260503(-12) & 5.237699(-12)\\
& $10^{-11}$ & 972 & 972 & 5 & 57 & 8.537615(-13) & 5.442313(-13) \\
& $10^{-12}$ & 975 & 975 & 5 & 54 & 8.926193(-14) & 5.584422(-14) \\
& $10^{-13}$ & 978 & 978 & 5 & 51 & 3.264056(-14) & 2.842171(-14) \\
\bottomrule
\end{tabular}
\end{table}
\begin{table}[H]
\centering
\renewcommand{\arraystretch}{1.2}
\caption{Values of the weights $\rho_k$ and exponents $\eta_k$ obtained by Prony's method for various values of $\alpha \in \{0.1, 0.5, 0.9\}$ over the intervals $[10^{-2}, 1]$, $[10^{-5}, 1]$, and $[10^{-2}, 10^3]$. Here, the value of $\epsilon$ is set to $10^{-10}$.}
\label{tab:10}
\begin{tabular}{ccccccc}
\toprule
\multirow{2}{*}{} & \multicolumn{2}{c}{$\boldsymbol{[10^{-2}, 1]}$} & \multicolumn{2}{c}{$\boldsymbol{[10^{-5}, 1]}$} & \multicolumn{2}{c}{$\boldsymbol{[10^{-2}, 10^3]}$} \\
\cmidrule(lr){2-3} \cmidrule(lr){4-5} \cmidrule(lr){6-7}
{$\boldsymbol{\alpha}$} & {$\boldsymbol{\rho_k}$} & {$\boldsymbol{\eta_k}$} & {$\boldsymbol{\rho_k}$} & {$\boldsymbol{\eta_k}$} & {$\boldsymbol{\rho_k}$} & {$\boldsymbol{\eta_k}$} \\
\midrule
\multirow{4}{*}{$0.1$  $(L=256)$} & 0.1887  & $-0.8580$ & 0.3076 & $-0.8810$ & $0.6138(-3)$ & $-0.8810(-3)$ \\
& 0.3202  & $-0.6074$ & 0.4804 & $-0.4834$ & $0.9586(-3)$ & $-0.4834(-3)$ \\
& 0.3384  & $-0.2926$ & 0.3347 & $-0.1006$ & $0.6677(-3)$ & $-0.1006(-3)$ \\
& 0.2033  & $-0.0569$ &  &  &  &  \\
\midrule
\multirow{5}{*}{$0.5$  $(L=256)$} & 0.2239 & $-0.9500$ & 0.2600 & $-0.8070$ & 0.0082 & $-0.8070(-3)$ \\
& 0.3026 & $-0.7184$ & 0.4149 & $-0.5384$ & 0.0131 & $-0.5384(-3)$ \\
& 0.4290 & $-0.4413$ & 0.5775 & $-0.2336$ & 0.0183 & $-0.2336(-3)$ \\
& 0.5265 & $-0.1795$ & 0.6670 & $-0.0284$ & 0.0211 & $-0.0284(-3)$ \\
& 0.5778 & $-0.0212$ &  &  &  &  \\
\midrule
\multirow{5}{*}{$0.9$  $(L=1024)$} & 0.2618 & $-0.8433$ & 0.2715 & $-0.9807$ & 0.1361 & $-0.9807(-3)$ \\
& 0.3766 & $-0.6154$ & 0.3836 & $-0.7131$ & 01922 & $-0.7131(-3)$ \\
& 0.6551 & $-0.3579$ & 0.6648 & $-0.4145$ & 0.3332 & $-0.4145(-3)$ \\
& 1.2235 & $-0.1258$ & 1.2416 & $-0.1458$ & 0.6223 & $-0.1458(-3)$ \\
& 7.4423 & $-0.0034$ & 7.5525 & $-0.0040$ & 3.7852 & $-0.0040(-3)$ \\
\bottomrule
\end{tabular}
\end{table}

\bibliographystyle{siam}
\bibliography{sample.bib}

\end{document}